%% file: draft.tex
\documentclass{siamart1116}

\usepackage{braket,amsfonts}

\usepackage{array}

\usepackage{graphicx,epstopdf}

\usepackage{booktabs}
\usepackage[caption=false]{subfig}
\usepackage{pgfplots}
\usetikzlibrary{calc}

\newsiamthm{claim}{Claim}
\newsiamremark{rem}{Remark}
\newsiamremark{expl}{Example}
\newsiamremark{hypothesis}{Hypothesis}
\crefname{hypothesis}{Hypothesis}{Hypotheses}

\usepackage{algorithmic}

\Crefname{ALC@unique}{Line}{Lines}

\numberwithin{theorem}{section}

\usepackage{xspace}
\usepackage{bold-extra}
\usepackage[most]{tcolorbox}

\colorlet{texcscolor}{blue!50!black}
\colorlet{texemcolor}{red!70!black}
\colorlet{texpreamble}{red!70!black}
\colorlet{codebackground}{black!25!white!25}

\lstdefinestyle{siamlatex}{%
  style=tcblatex,
  texcsstyle=*\color{texcscolor},
  texcsstyle=[2]\color{texemcolor},
  keywordstyle=[2]\color{texemcolor},
  moretexcs={cref,Cref,maketitle,mathcal,text,headers,email,url},
}

\tcbset{%
  colframe=black!75!white!75,
  coltitle=white,
  colback=codebackground, 
  colbacklower=white, 
  fonttitle=\bfseries,
  arc=0pt,outer arc=0pt,
  top=1pt,bottom=1pt,left=1mm,right=1mm,middle=1mm,boxsep=1mm,
  leftrule=0.3mm,rightrule=0.3mm,toprule=0.3mm,bottomrule=0.3mm,
  listing options={style=siamlatex}
}

\newtcblisting[use counter=example]{example}[2][]{%
  title={Example~\thetcbcounter: #2},#1}

\newtcbinputlisting[use counter=example]{\examplefile}[3][]{%
  title={Example~\thetcbcounter: #2},listing file={#3},#1}

\DeclareTotalTCBox{\code}{ v O{} }
{ 
  fontupper=\ttfamily\color{black},
  nobeforeafter,
  tcbox raise base,
  colback=codebackground,colframe=white,
  top=0pt,bottom=0pt,left=0mm,right=0mm,
  leftrule=0pt,rightrule=0pt,toprule=0mm,bottomrule=0mm,
  boxsep=0.5mm,
  #2}{#1}

\patchcmd\newpage{\vfil}{}{}{}
\graphicspath{{../fig/plots/}}

\newcommand{\dudx}[1]{\frac{d #1}{d x}}
\newcommand{\ddudxx}[1]{\frac{d^2 #1}{d x^2}}

\newtheorem{Lemma}{Lemma}

\begin{tcbverbatimwrite}{tmp_\jobname_header.tex}
\title{Deep Multigrid: learning prolongation and restriction matrices%
  \thanks{Submitted to the editors 10 November 2017.
\funding{This study was supported by the Ministry of Education and Science of the Russian Federation grant 14.756.31.0001 (reformulation of the multigrid method in terms of deep neural network  and training multigrid operators with the backpropagtion approach) and by RSF grant 17-11-01376 (idea to replace the objective functional (spectral radius) by the approximation using Gelfand's formula and stochastic trace estimator)}}}
\author{Alexandr Katrutsa%
  \thanks{Skolkovo Institute of Science and Technology, Nobel St., 3, 143025, Moscow, Russia, and Moscow Institute of Physics and Technology, Institutskii per. 9, Dolgoprudny, 141700, Moscow Region, Russia (\email{aleksandr.katrutsa@phystech.edu})}
  \and
  Talgat Daulbaev%
  \thanks{
  Lomonosov Moscow State University,GSP-1, Leninskie Gory, Moscow, 119991, Russian Federation, and Skolkovo Institute of Science and Technology, Nobel St., 3, 143025, Moscow, Russia, (\email{PVnuRT@gmail.com})}
  \and
Ivan Oseledets%
  \thanks{Skolkovo Institute of Science and Technology, Nobel St., 3, 143025, Moscow, Russia, and Institute of Numerical Mathematics, Russian Academy of Sciences, Gubkina St. 8,
119333 Moscow, Russia (\email{i.oseledets@skoltech.ru})}
}

\headers{Deep Multigrid}
{Alexandr Katrutsa, Talgat Daulbaev, and Ivan Oseledets}
\end{tcbverbatimwrite}
\input{tmp_\jobname_header.tex}

\begin{document}
\maketitle

\begin{tcbverbatimwrite}{tmp_\jobname_abstract.tex}
\begin{abstract}
This paper proposes the method to optimize restriction and prolongation operators in the two-grid method.
The proposed method is straightforwardly extended to the geometric multigrid method (GMM).
GMM is used in solving discretized partial differential equation (PDE) and based on the restriction and prolongation operators.
The operators are crucial for fast convergence of GMM, but they are unknown.
To find them we propose a reformulation of the two-grid method in terms of a deep neural network with a specific architecture.
This architecture is based on the idea that every operation in the two-grid method can be considered as a layer of a deep neural network.
The parameters of layers correspond to the restriction and prolongation operators.
Therefore, we state an optimization problem with respect to these operators and get optimal ones through backpropagation approach.
To illustrate the performance of the proposed approach, we carry out experiments on the discretized Laplace equation, Helmholtz equation and singularly perturbed convection-diffusion equation and demonstrate that proposed approach gives operators, which lead to faster convergence.
\end{abstract}

\begin{keywords}
geometric multigrid method, deep neural network, spectral radius minimization, Helmholtz equation, Poisson equation
\end{keywords}

\begin{AMS}
65N55, 65M55, 35Q93 
\end{AMS}
\end{tcbverbatimwrite}
\input{tmp_\jobname_abstract.tex}

\section{Introduction}

In this paper, we propose a method to optimize the parameters of the geometric multigrid method (GMM). 
GMM is often a method of choice for solving large sparse systems arising from partial differential equation (PDE) discretization~\cite{hackbusch2013multi,briggs2000multigrid}.
The main challenge in GMM is to define \emph{prolongation} and \emph{restriction} operators.
These operators are constructed from the prior knowledge about system structure or different heuristics~\cite{hackbusch2013multi}.

To treat this problem we develop \emph{Deep Multigrid Method}, where we find parameters of GMM directly to minimize the spectral radius of the iteration matrix, which depends on the restriction and prolongation operators.
The method is called ``deep'' because every operation in one iteration of GMM can be considered as a layer in a deep neural network with a specific architecture and every iteration of GMM is equivalent to a forward pass in this neural network.
To use deep neural network methods in optimization of the prolongation and restriction operators, one has to define a loss function, which estimates the quality of a given pair of operators.
Since GMM can be represented as an iterative process, the natural choice of the loss function is the spectral radius of the iteration matrix, which depends on the restriction and prolongation operators and significantly affects the speed of GMM convergence.
However, the minimization problem of the spectral radius is non-smooth and computation of the gradient of the objective function requires left and right leading eigenvectors, which is too costly.
Therefore, instead of using spectral radius as the objective function we propose a stochastic estimate derived from the Gelfand's formula~\cite{kozyakin2009accuracy}.
To compute this estimate, we need only a small number of GMM iterations for different initial vectors.
Gradients of the proposed estimate can be computed with any automatic differentiation tool, like Autograd~\cite{maclaurin2015autograd}, TensorFlow~\cite{abadi2016tensorflow}, Theano~\cite{theano} and others.
We use Autograd in our experiments.
To minimize the introduced approximation, we use Adam optimizer~\cite{kingma2014adam}, which is a modification of the pure stochastic gradient descent.
This optimizer requires a good initial approximation, therefore we use the homotopy optimization approach~\cite{watson1989modern} to find an appropriate initial approximation of restriction and prolongation operators.

In the numerical experiments, we compare the proposed method with standard GMM and Algebraic Multigrig Method (AMG)~\cite{trottenberg2000multigrid,ruge1987algebraic,pyamg} and obtain better convergence for model problems. 

The main contributions of this study are the following.
\begin{itemize}
\item We reformulate GMM in terms of deep neural network with specific architecture
\item We introduce a loss function, which estimates the spectral radius of the iteration matrix
\item We use homotopy initialization approach to find initial approximation for optimization method
\item We demonstrate performance of the Deep Multigrid method in comparison with GMM for considered types of differential equation
\end{itemize}

\section{Multigrid method as a Neural Network}
\label{deepmult::sec::mm_nn}
To demonstrate the performance of the proposed approach we only consider a two-grid method.
The extension from two-grid to multigrid is straightforward.
The general form of the classical two-grid method for solving a linear system
\begin{equation}
Au = f, \quad A \in \mathbb{R}^{n \times n}
\label{deepmult::eq::general_system}
\end{equation}
is provided in the next Section.

\subsection{Two-grid method}
\label{deepmult::sec::2grid}
The idea of the two-grid method comes from the consideration of linear systems arising in PDE discretization on a sequence of grids.
Let $A$ be a matrix of the linear system obtained after PDE discretization on a fine grid.

Let $u^{(k)} \in \mathbb{R}^n$ be a current approximation. 
The \emph{pre-smoothing} step, consisting of $s_1$ iterations, is defined from the splitting of the matrix $A = M_1 - K_1$ and has the form:
\begin{equation}
u^{(k)}_{pre} = (M_1^{-1}K_1)^{s_1} u^{(k)} + s_1M_1^{-1}f.
\label{deepmult::eq::presmoother}
\end{equation}
 

After that, we compute the residual vector 
\begin{equation}
r^{(k)} = Au^{(k)}_{pre} - f.
\label{deepmult::eq::error_vector}
\end{equation}
The next step of the two-grid method is the restriction of the residual vector $r^{(k)}$ to the coarse grid.
To perform this restriction, we introduce the \emph{restriction operator} $R \in \mathbb{R}^{n_c \times n}$, where $n_c$ is the size of the coarse grid, and compute the restricted residual vector $r^{(k)}_{c}$
\begin{equation}
r^{(k)}_{c} = Rr^{(k)}.
\label{deepmult::eq::error_vector_restr}
\end{equation}
Now we assume that we can accurately solve the following restricted linear system:
\begin{equation}
A_{c}u^{(k)}_{c} = r^{(k)}_{c},
\label{deepmult::eq::2grid_course_sys}
\end{equation}
where $A_{c} \in \mathbb{R}^{n_c \times n_c}$ is the matrix on the coarse grid.
The usual choice is the Galerkin projection method, where
\begin{equation}
A_{c} = RAP,
\label{deepmult::eq::A_2h}
\end{equation}
and $P \in \mathbb{R}^{n \times n_c}$ is the \emph{prolongation operator} that maps a vector from the coarse grid to the fine grid.
After solving Equation~\eqref{deepmult::eq::2grid_course_sys}, the solution $u_{c}$ has to be prolongated to the fine grid with the prolongation operator $P$, and
\begin{equation}
\hat{u}^{(k)} = u^{(k)}_{pre} + Pu_{2h}.
\label{deepmult::eq::2grid_update}
\end{equation}
After that we perform $s_2$ \emph{postsmoothing} iterations:
\begin{equation}
u^{(k+1)} = (M^{-1}_2K_2)^{s_2} \hat{u}^{(k)} + s_2M^{-1}_2f,
\label{deepmult::eq::postsmoother}
\end{equation}
where $M_2$ and $K_2$ satisfy $A = M_2 - K_2$ and $M_2$ is invertible.

In our work the \emph{damped Jacobi method}~\cite{briggs2000multigrid} with parameter~$\omega$ performs pre- and postsmoothing.
In this case, the matrix $A$ is split in the form $A = \omega^{-1}D - \omega^{-1}D + A$, where $D$ is the diagonal matrix with diagonal equal to the diagonal of the matrix $A$.
Therefore, the matrices $M_1 = M_2 = \omega^{-1}D$ are invertible if the diagonal of $A$ has no zero elements, and the matrices $K_1 = K_2 = \omega^{-1}D - A$.

Now we summarize the operations which are performed in equations~\cref{deepmult::eq::presmoother,deepmult::eq::error_vector,deepmult::eq::error_vector_restr,deepmult::eq::2grid_course_sys,deepmult::eq::A_2h,deepmult::eq::2grid_update,deepmult::eq::postsmoother}.
Backward substitution procedure gives the following form of one iteration of the two-grid method:


\begin{equation}
u^{(k+1)} = Cu^{(k)} + b,
\label{deepmult::eq::gen_iter_proc}
\end{equation}
where $C$ is the \emph{iteration matrix} given by
\begin{equation}
C = (M^{-1}_2K_2)^{s_2}(I +  P(RAP)^{-1}RA)(M_1^{-1}K_1)^{s_1},
\label{deepmult::eq::2grid_C}
\end{equation}
where $I$ is the identity matrix, and the vector $b$ is computed as
\[
b = ((M^{-1}_2K_2)^{s_2}P(RAP)^{-1}R(s_1AM_1^{-1} - I) + s_2M_2^{-1})f.
\]
Note that, to multiply the matrix $C$ by any vector $y$ one needs to perform one iteration of the two-grid method with $u^{(k)} = y$ and $f \equiv 0$.
The iteration matrix $C$ depends on the operators $R$, $P$ and the damp factor $\omega$: $C = C(R, P, \omega)$.
Further to simplify the notation we write $C$ in place of $C(R, P, \omega)$.

\subsection{Spectral radius optimization problem}
\label{deepmult::sec::opt_problem}

Since we formulate the two-grid method as an iterative process with the iteration matrix $C$, see equations~\eqref{deepmult::eq::gen_iter_proc},~\eqref{deepmult::eq::2grid_C}, convergence of this iterative process is estimated as~\cite{olshanskii2014iterative} 
\[
\|u^{(k)} - u^* \|_2 \leq \rho^k(C) \|u^{(0)} - u^*\|_2,
\]
where $\rho(C)$ is the spectral radius of the iteration matrix $C$, $u^{(0)}$ is an initial approximation, $u^{(k)}$ is the approximation after $k$ iterations and $u^*$ is the solution of the linear system $Au^* = f$.
Thus, it is natural to determine operators $R, P$ and damp factor $\omega$ from the following optimization problem
\begin{equation}
\rho(C) = \max_{i = 1,\dots,n} |\lambda_i(C)| \to \min\limits_{R, P, \omega}.
\label{deepmult::eq::spect_norm_min_init}
\end{equation}

However, this optimization problem is non-smooth~\cite{nesterov2007smoothing} and the computation of the objective subgradient is too costly because it requires left and right leading eigenvectors~\cite{nesterov2007smoothing}. 
Therefore, we need some approximation of this objective which is more appropriate for minimization.
We propose to replace $\rho(C)$ with its approximation
\begin{equation}
\rho(C) \approx \sqrt[K]{\|C^K\|_F},
\label{deepmult::eq::}
\end{equation} 
where $K$ is a positive integer constant and $\|\cdot \|_F$ is the Frobenius matrix norm.
This approximation is inspired by the Gelfand's formula~\cite{kozyakin2009accuracy}
\begin{equation}
\rho(C) = \lim \limits_{k \to \infty} \sqrt[k]{\|C^k \|} \ ,
\label{deepmult::eq::gelfand}
\end{equation}
where $\|\cdot \|$ is any matrix norm.
The introduced approximation is used to bound the spectral radius~\cite{kozyakin2009accuracy}
\begin{equation}
\gamma^{(1 + \ln K)/ K}\| C^K \|^{1/K}_F \leq \rho(C) \leq \| C^K \|^{1/K}_F,
\label{deepmult::eq::rho_bounds}
\end{equation}
for some $\gamma \in (0, 1)$ and any positive $K$.
This bound becomes tighter for larger constant $K$.
Therefore, we fix $K$ and minimize the upper bound in~\eqref{deepmult::eq::rho_bounds}.
This minimization problem is equivalent to the following problem:
\begin{equation}
F_K = \| C^K\|^2_F \to \min_{R, P, \omega}.
\label{deepmult::eq::fro_norm_min}
\end{equation}
We can compute the product of the matrix $C$ by any vector with one iteration of the two-grid method.
Therefore, we use the following stochastic unbiased estimate of the objective in problem~\eqref{deepmult::eq::fro_norm_min}~\cite{avron2011randomized}
\begin{equation}
\| C^K\|^2_F = \mathbb{E}_z \|C^K z \|^2_2,
\label{deepmult::eq::fro_norm_expect}
\end{equation}
where $z = [z_i]$ is a random vector, which elements are generated by Rademacher distribution:
\begin{equation}
\mathbb{P}(z_i = \pm 1) = \frac{1}{2}.
\label{deepmult::eq::rademacher}
\end{equation}
So, we minimize the expectation
\begin{equation}
\mathbb{E}_z \|C^K z \|^2_2 \to \min_{R, P, \omega}.
\label{deepmult::eq::final_opt}
\end{equation}
To minimize~\eqref{deepmult::eq::final_opt}, it is natural to use a stochastic gradient based method.
The stochastic unbiased estimate of the objective function in~\eqref{deepmult::eq::final_opt} is the following:
\begin{equation}
\hat{F}_K = \frac{1}{N} \sum\limits_{i=1}^N \|C^K z^i\|^2_2,
\label{deepmult::eq::F_estimate}
\end{equation}
where $N$ is the batch size and $z^i$ is the $i$-th random vector with elements generated by Rademacher distribution~\eqref{deepmult::eq::rademacher}. 

Following Lemma is a direct consequence of the stochastic trace estimator~\cite{hutchinson1990stochastic} and proves that estimate~\eqref{deepmult::eq::F_estimate} is unbiased.
\begin{Lemma}
The used stochastic estimation~\eqref{deepmult::eq::F_estimate} of the objective in~\eqref{deepmult::eq::fro_norm_min} has the following expectation and variance:
\[
\mathbb{E}_{z^i} \left(\frac{1}{N} \sum\limits_{i=1}^N \|C^K z^i\|^2_2 \right) = \|C^K\|^2_F
\]
\[
\mathrm{Var}_{z^i} \left(\frac{1}{N} \sum\limits_{i=1}^N \|C^K z^i\|^2_2 \right) = \frac{2}{N}\left( \left\|\left(C^K\right)^{\top} C^K \right\|^2_F - \sum^n_{i=1} \left(\left(C^K\right)^{\top} C^K\right)^2_{ii} \right).
\]
\end{Lemma} 

Thus, we have the unbiased estimate~\eqref{deepmult::eq::F_estimate} of the objective function.
Therefore, we can compute the stochastic gradient of the objective as the gradient of this unbiased estimate with respect to operators $R, P$ and damp factor $\omega$.
To compute the stochastic gradient we use automatic differentiation tool Autograd~\cite{maclaurin2015autograd}, which requires a function for computing $\hat{F}_K$ as a superposition of differentiable operations and returns the gradient of $\hat{F}_K$ with respect to operators $R, P$ and damp factor $\omega$.
Computation of $\hat{F}_K$ consists of $K$ matrix $C$ by vector products for every generated vector $z^i$, computing squared 2-norm of the obtained vectors and its averaging. 
So, the most complicated operation here is the multiplication of the matrix $C$ by some vector.
However, in Section~\ref{deepmult::sec::2grid} we have shown that the multiplication of the matrix $C$ by a vector is equivalent to one iteration of the two-grid method with a given vector as initial approximation and zero right-hand side. 
The Algorithm for the computation of the estimate $\hat{F}_K$ is summarized in Algorithm~\ref{deepmult::alg::F_compute}.

\begin{algorithm}[H]
 \caption{Computation of spectral radius estimation $\hat{F}_K$}
 \label{deepmult::alg::F_compute}
 \begin{algorithmic}[1]
\REQUIRE{Batch size $N$,  positive integer $K$, function \texttt{Two-grid}$(u_0, f|A, R, P, \omega)$ that performs one iteration of two-grid method for initial approximation $u_0$ and right-hand side $f$, given matrix $A$, operators $R, P$ and damp factor $\omega$}
\ENSURE{Estimate of the $\rho(C) \approx \hat{F}_K$}
  \STATE {$\hat{F}_K = 0$}
 \FOR{$i = 1,\ldots, N$}
  \STATE {Generate random vector $x$ with elements from Rademacher distribution~\eqref{deepmult::eq::rademacher}}
  \FOR{$k = 1, \ldots, K$}
    \STATE {$x = \texttt{Two-grid}(x, 0)$}
  \ENDFOR
  \STATE {$\hat{F}_K = \hat{F}_K + \|x\|^2_2$}
 \ENDFOR
 \STATE {$\hat{F}_K = \hat{F}_K / N$}
 \RETURN {$\hat{F}_K$}
 \end{algorithmic}
\end{algorithm}

To represent the function $\texttt{Two-grid}$ as a superposition of differentiable operations we reformulate it in the form of a deep neural network.

\subsection{Deep neural network reformulation}

In this section, we present a reformulation of the two-grid method in terms of a deep neural network.
Extension of this reformulation for a general GMM is straightforward.
To provide the required reformulation, we represent operations in the two-grid method as layers in a deep neural network.
The first operation in the two-grid method is presmoothing~\eqref{deepmult::eq::presmoother}, which is performed by a fixed linear transformation of the current approximation $u^{(k)}$ with an unknown damp factor $\omega$.
This operation can be represented as a fully connected layer with a given matrix and learned damp factor $\omega$.
The next operation is the computation of the residual vector $r^{(k)}$~\eqref{deepmult::eq::error_vector} which is also linear and can be represented as a fully connected layer with a fixed matrix (no parameters to be learned).
After that, we have to restrict the residual vector $r^{(k)}$ on the coarse grid by the unknown matrix $R$ and get a vector $r^{(k)}_{c}$.
Therefore, we represent this operation as a fully connected layer with unknown matrix $R$, which is optimized by the proposed approach.
After that, we need to get the projected matrix $A_{c}$, which is computed according to~\eqref{deepmult::eq::A_2h}.
We represent this operation as a \emph{Projection layer} with unknown parameters $P$ and $R$, which are going to be optimized.
After projection we have to solve Equation~\eqref{deepmult::eq::2grid_course_sys} and get solution $u^{(k)}_{c}$.
This operation is performed by \emph{Solver layer}, which gives solution $u^{(k)}_{c}$ of the linear system with matrix $A_{c}$ and right-hand side $r^{(k)}_{c}$.
The matrix $A_c$ depends on the unknown matrices $P$ and $R$~\eqref{deepmult::eq::A_2h}, therefore the output of this layer is differentiable with respect to $P$ and $R$.
After that we prolongate this solution on the fine grid and update presmoothed $u^{(k)}$ according to~\eqref{deepmult::eq::2grid_update}.
Similar to the restriction operation, the prolongation is performed by a fully connected layer with unknown matrix $P$.
The final operation is postsmoothing~\eqref{deepmult::eq::postsmoother}, which is the same as a presmoothing.

\begin{figure}[!ht]
\centering
\scalebox{0.9}{
\begin{tikzpicture}
\coordinate (input_u) at (0, 0);
\coordinate (input_A) at (2.1, 0);
\coordinate (input_f) at (2.6, 0);
\coordinate (presmoother_left_up) at (-1.5, -1);
\coordinate (presmoother_right_bottom) at (1.5, -2);
\coordinate (y_scale) at (0, -2);
\node[scale=1] (u) at (input_u) {$u^{(k)}$};
\node (A) at (input_A) {$A$};
\node (f) at (input_f) {$f$};

\draw (presmoother_left_up) rectangle (presmoother_right_bottom) node[pos=.5] {Presmoothing};
\draw[-latex] ($(input_u) + (0, -0.2)$) -- ($(presmoother_left_up) + (1.5, 0)$);
\draw[-latex] ($(input_f) + (0, -0.2)$) to [bend left=50, bend right=20] ($(presmoother_left_up) + (2.8, 0)$);
\draw[-latex] ($(input_A) + (0, -0.2)$) to [bend left=50, bend right=25] ($(presmoother_left_up) + (2.4, 0)$);

\coordinate (res_vec_left_up) at ($(presmoother_left_up) + (y_scale)$);
\coordinate (res_vec_right_bottom) at ($(presmoother_right_bottom) + (y_scale)$);
\draw[-latex] ($(presmoother_right_bottom) + (-1.5, 0)$) -- ($(res_vec_left_up) + (1.5, 0)$) node[pos=0.5,left,scale=1] {$u^{(k)}_{pre}$};
\draw (res_vec_left_up) rectangle (res_vec_right_bottom) node[pos=0.5] {Residual vector};
\draw[-latex] ($(input_A) + (0, -0.2)$) to [bend left=10, bend right=-25] ($(res_vec_left_up) + (2.4, 0)$);
\draw[-latex] ($(input_f) + (0, -0.2)$) to [bend left=10, bend right=-20] ($(res_vec_left_up) + (2.8, 0)$);

\coordinate (restr_left_up) at ($(res_vec_left_up) + (y_scale)$);
\coordinate (restr_right_bottom) at ($(res_vec_right_bottom) + (y_scale)$);
\draw (restr_left_up) rectangle (restr_right_bottom) node[pos=0.5] {Restriction};
\draw[-latex] ($(res_vec_right_bottom) + (-1.5, 0)$) -- ($(restr_left_up) + (1.5, 0)$) node[scale=1,pos=0.5,left] {$r^{(k)} = Au^{(k)}_{pre} - f$};

\coordinate (galerkin_left_up) at ($(restr_left_up) + (y_scale) + (3, 0)$);
\coordinate (galerkin_right_bottom) at ($(restr_right_bottom) + (y_scale) + (3, 0)$);
\draw (galerkin_left_up) rectangle (galerkin_right_bottom) node[scale=1,pos=0.5] {Matrix projection};
\draw[-latex] ($(input_A) + (0, -0.2)$) to[bend left=10,bend right=-10] ($(galerkin_left_up) + (1.5, 0)$);

\coordinate (solv_lin_sys_left_up) at ($(restr_left_up) + 2*(y_scale)$);
\coordinate (solv_lin_sys_right_bottom) at ($(restr_right_bottom) + 2*(y_scale)$);
\draw (solv_lin_sys_left_up) rectangle (solv_lin_sys_right_bottom) node[scale=0.8,pos=0.5] {Solving linear system};
\draw[-latex] ($(galerkin_right_bottom) + (-1.5, 0)$) to[out=-90,in=0] ($(solv_lin_sys_right_bottom) + (0, 0.5)$);
\node (A_coarse) at ($(galerkin_right_bottom) + (-0.9, -1.3)$) {$A_{c} = RAP$};
\draw[-latex] ($(restr_right_bottom) + (-1.5,0)$) -- ($(solv_lin_sys_left_up) + (1.5, 0)$) node[pos=0.5,left,scale=1] {$r^{(k)}_{c} = Rr^{(k)}$};

\coordinate (prolong_left_up) at ($(solv_lin_sys_left_up) + (y_scale)$);
\coordinate (prolong_right_bottom) at ($(solv_lin_sys_right_bottom) + (y_scale)$);
\draw (prolong_left_up) rectangle (prolong_right_bottom) node[scale=0.7,pos=0.5] {Prolongation and update};
\draw[-latex] ($(solv_lin_sys_right_bottom) + (-1.5, 0)$) -- ($(prolong_left_up) + (1.5, 0)$) node[scale=1,left,pos=0.5] {$u^{(k)}_{c} = A^{-1}_{c}r^{(k)}_{c}$};
\draw[-latex] ($(presmoother_left_up) + (0, -0.5)$) to[in=180,out=180] ($(prolong_left_up) + (0, -0.5)$);
\node (u_to_update) at ($(galerkin_left_up) + (-6.5, -0.5)$) [scale=1] {$u^{(k)}_{pre}$};

\coordinate (postsmooth_left_up) at ($(prolong_left_up) + (y_scale)$);
\coordinate (postsmooth_right_bottom) at ($(prolong_right_bottom) + (y_scale)$);
\draw (postsmooth_left_up) rectangle (postsmooth_right_bottom) node[scale=1,pos=0.5] {Postsmoothing};
\draw[-latex] ($(prolong_right_bottom) + (-1.5,0)$) -- ($(postsmooth_left_up) + (1.5, 0)$) node[pos=0.5,left,scale=1] {$\hat{u}^{(k)} = u^{(k)}_{pre} + Pu^{(k)}_{c}$};
\draw[-latex] ($(postsmooth_right_bottom) + (-1.5, 0)$) -- ($(postsmooth_right_bottom) + (-1.5, 0) + (0, -1)$) node[below, scale=1] {$u^{(k+1)}$};

\coordinate (par_xleft_shift) at (7, 0);
\coordinate (par_xright_shift) at (5, 0);
\coordinate (par_restr_left_up) at ($(restr_left_up) + (par_xleft_shift)$);
\coordinate (par_restr_right_bottom) at ($(restr_right_bottom) + (par_xright_shift)$);
\draw (par_restr_left_up) rectangle (par_restr_right_bottom) node[pos=0.5] {$R$};
\draw[-latex] ($(par_restr_left_up) + (0, -0.5)$) -- ($(restr_right_bottom) + (0, 0.5)$);
\draw[-latex] ($(par_restr_right_bottom) + (-0.5, 0)$) to[out=-90,in=0] ($(galerkin_right_bottom) + (0, 0.67)$);

\coordinate (par_prolong_left_up) at ($(prolong_left_up) + (par_xleft_shift)$);
\coordinate (par_prolong_right_bottom) at ($(prolong_right_bottom) + (par_xright_shift)$);
\draw (par_prolong_left_up) rectangle (par_prolong_right_bottom) node[pos=0.5] {$P$};
\draw[-latex] ($(par_prolong_left_up) + (0, -0.5)$) -- ($(prolong_right_bottom) + (0, 0.5)$);
\draw[-latex] ($(par_prolong_left_up) + (0.5, 0)$) to[out=90,in=0] ($(galerkin_right_bottom) + (0, 0.33)$);

\coordinate (damp_factor_left_up) at ($(par_restr_right_bottom) + (0.2, -1)$);
\coordinate (damp_factor_right_bottom) at ($(damp_factor_left_up) + (1, -1)$);
\draw (damp_factor_left_up) rectangle (damp_factor_right_bottom) node[scale=1,pos=0.5] {$\omega$};
\draw[-latex] ($(damp_factor_left_up) + (0.5, 0)$) to[out=90,in=0] ($(presmoother_right_bottom) + (0, 0.5)$);
\draw[-latex] ($(damp_factor_right_bottom) + (-0.5, 0)$) to[out=-90,in=0,looseness=1.6] ($(postsmooth_right_bottom) + (0, 0.5)$);
\end{tikzpicture}
}
\caption{Deep Multigrid neural network architecture}
\label{deepmult::fig::deep_nn}
\end{figure}
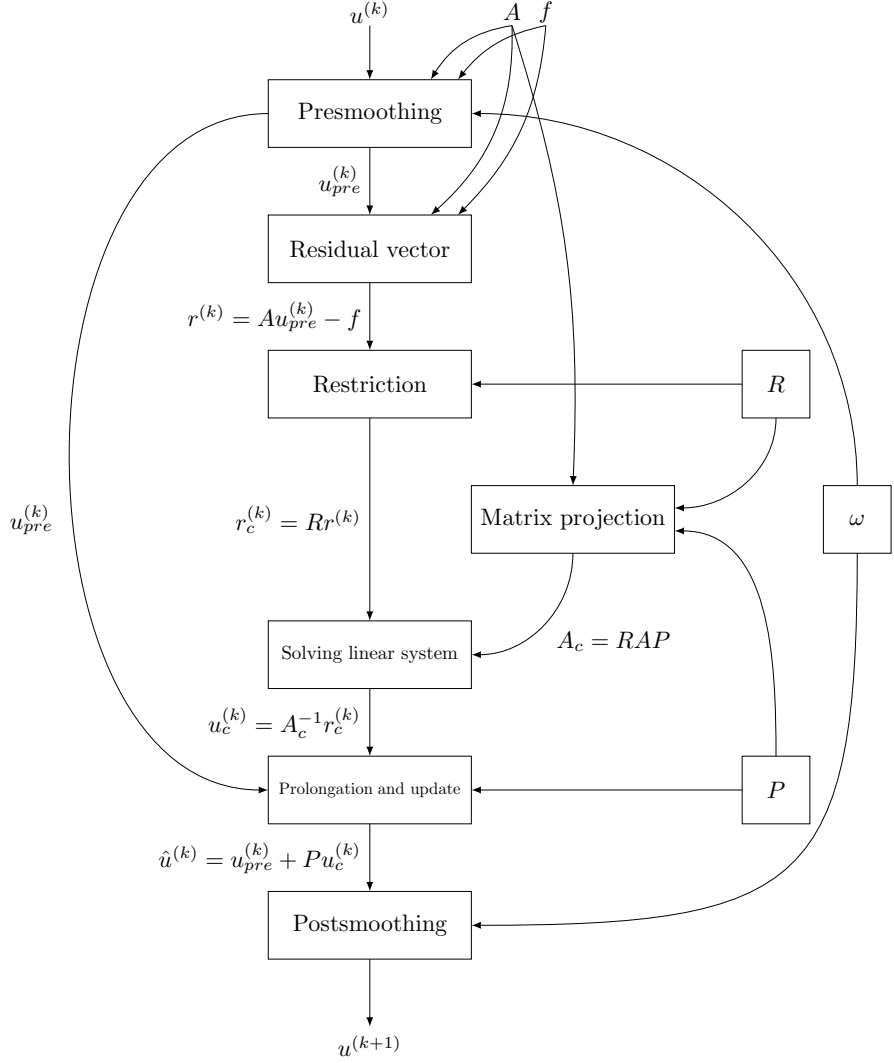

Note that every layer described above is differentiable, therefore backpropagation approach~\cite{hecht1988theory} can be applied to optimize $P, R$ and $\omega$. 
The layers described above and connections between them are summarized in Fig.~\ref{deepmult::fig::deep_nn}.

\section{Initialization approaches}
A good initial approximation is crucial for the convergence of the gradient-type methods.
These methods can converge to different local minima from different initial approximations. 
We observed this in our problem.
Therefore, the choice of initial approximation is important.
We can use an \emph{apriori} knowledge about the problem to choose better initial approximation.

\subsection{Standard initialization}
In this approach, we initialize operators $R, P$ and damp factor $\omega$ with some values that are good for some known problem and perform iterations of the optimization process until the convergence is achieved.
For example, it is well-known~\cite{briggs2000multigrid} that linear interpolation operators $R, P$~\eqref{deepmult::eq::linear_inter} and damp factor $\omega_0 = \frac{2}{3}$ are appropriate for solving Poisson equation.
Therefore, we can use it as initialization for solving optimization problem~\eqref{deepmult::eq::final_opt} for $A$ coming from close problems, i.e. Helmholtz or diffusion-convection equations. 
However, our experiments for the Helmholtz equation in a high-frequency setting show that this initialization leads to the convergence to a local minimum with spectral radius $\rho(C)$ greater than $1$.
Therefore, we use homotopy approach~\cite{watson1989modern} to find better local minimum.

\subsection{Homotopy approach}
Assume we want to solve the optimization problem~\eqref{deepmult::eq::final_opt} for a matrix $A_1$ for which we do not know an appropriate initialization of the parameters, but we know it for the problem corresponding to some matrix $A_0$.
So we introduce a sequence of matrices $M_i$ according to the following equation:
\begin{equation}
\begin{split}
& M_i = \alpha_i A_1 + (1 - \alpha_i)A_0, \\ 
& \alpha_0 = 0, \; 0 < \alpha_1 < \alpha_2 < \ldots < \alpha_{k-1} < 1, \; \alpha_k = 1.
\end{split}
\label{deepmult::eq::homotopy}
\end{equation}
This sequence helps to tune initialization for the optimization problem corresponding to the matrix $A_1$ in the following way. 

Assume we know a good initialization for the problem corresponding to the matrix $M_i$.
Therefore, we can solve it and get parameters $R_i$, $P_i$ and $\omega_i$, which establish an appropriate local minimum.
After that we use this solution as initialization of the problem for the matrix $M_{i+1}$.
If the matrices $M_i$ and $M_{i+1}$ are close enough, then according to the continuity of the problem and homotopy theory~\cite{watson1989modern} we get a good local minimum for the problem corresponding to the matrix $M_{i+1}$.
Further, we use the solution $R_{i+1}$, $P_{i+1}$ and $\omega_{i+1}$ as an initialization for the problem  corresponding to the matrix $M_{i+2}$ and so on.
But how to measure the quality of the local minimum obtained for the matrix $M_{i+1}$?
To measure this quality we introduce the \emph{acceptance rate} $\tau \in (0, 1)$ and number of trials $p = 1$.
The local minimum is good if the spectral radius for the solution $R_{i+1}$, $P_{i+1}$ and $\omega_{i+1}$ is smaller than $\tau$. 
If the local minimum for $M_{i+1}$ is not good, it means that the matrices are not close enough.
Therefore, we recompute $M_{i+1}$ by decreasing the value of $\alpha_{i+1}$, for example according to the following equation 
\[
\alpha_{i+1} = \alpha_i + \frac{\delta}{2^p},
\]
where $\delta \in (0, 1)$ is a fixed \emph{homotopy step size}.
After that we increase $p$ by one: $p = p + 1$.
This process continues until we find a good local minimum and then repeat it for the other $M_i$.
The idea is summarized in Algorithm~\ref{deepmult::alg::homotopy}.

\begin{algorithm}[H]
 \caption{Homotopy initialization}
 \label{deepmult::alg::homotopy}
 \begin{algorithmic}[1]
 \REQUIRE{Acceptance rate $\tau$, number $\alpha_{i}$, homotopy step size $\delta$, matrix $M_i$, matrices $A_0$ and $A_1$}
 \ENSURE{Appropriate initialization $R_{i+1}$, $P_{i+1}$ and $\omega_{i+1}$ for the problem corresponding to the matrix $M_{i+2}$}
 \STATE{Solve optimization problem for the matrix $M_i$ and get local optimal parameters $R_i$, $P_i$ and $\omega_i$}
 \STATE {Compute $\alpha_{i+1} = \alpha_i + \delta$}
 \STATE {Compute $M_{i+1} = (1 - \alpha_{i+1})A_0 + \alpha_{i+1}A_1$}
 \STATE {Solve the optimization problem for the matrix $M_{i+1}$ initialized by $R_i$, $P_i$ and $\omega_i$ and get corresponding local optimal parameters $R_{i+1}$, $P_{i+1}$ and $\omega_{i+1}$}
 \STATE {Compute estimate of the spectral radius $\hat{F}$ for $R_{i+1}$, $P_{i+1}$ and $\omega_{i+1}$}
 \STATE {$p = 1$}
 \WHILE{$\hat{F} \geq \tau$} 
     \STATE {Compute $\alpha_{i+1} = \alpha_i + \frac{\delta}{2^p}$}
     \STATE {$p = p + 1$}
     \STATE {Update $M_{i+1} = (1 - \alpha_{i+1})A_0 + \alpha_{i+1}A_1$}
     \STATE {Update $R_{i+1}$, $P_{i+1}$ and $\omega_{i+1}$ by solving the optimization problem for the matrix $M_{i+1}$ initialized by $R_i$, $P_i$ and $\omega_i$}
     \STATE {Compute estimate of the spectral radius $\hat{F}$ for $R_{i+1}$, $P_{i+1}$ and $\omega_{i+1}$}
 \ENDWHILE
 \RETURN {$R_{i+1}$, $P_{i+1}$ and $\omega_{i+1}$}
\end{algorithmic}
\end{algorithm}

For example, we know that linear interpolation operators and damp factor $\omega_0 = \frac{2}{3}$ are appropriate for the solution of the discretized Poisson equation. 
Therefore, we consider the  matrix from the discretization of the Poisson equation as a matrix $A_0$ and the matrix from the discretization of the Helmholtz equation as matrix $A_1$.
In this case the homotopy helps to find much better starting guess.

\section{Complexity}
The most time consuming operation in one iteration of the proposed Deep Multigrid method is the computation of the gradient of the objective function.
From the paper~\cite{baur1983complexity} it follows that the complexity of the computation of the objective function gradient is upper bounded by $3M$, where $M$ is the complexity of the objective function computation, and auto-differentiation tools implement this approach.
Therefore, the total complexity of the Deep Multigrid method is linearly dependent on the complexity of the one iteration of the two-grid method.
However, in the case of using homotopy approach the Deep Multigrid method becomes more computationally expensive, since there are a lot of intermediate subproblems that have to be solved.

\section{Numerical Experiments}
\label{deepmult::sec::num_exp}

In this section we present numerical experiments, which demonstrate the performance of the proposed approach.
We select some model problems described below to show that Deep Multigrid method (DMG) gives the operators which establish faster convergence of the two-grid method compared to the linear interpolation operators and to the two-grid Algebraic Multigrid method (AMG).

\subsection{Model problems overview}
To demonstrate the performance of the proposed DMG we consider the following model problems and discretizations. 
Grid size is $n = 2^l - 1$ for $l > 2$ and $l \in \mathbb{N}$.
Denote by $h = \frac{1}{n+1}$ the grid step.

\paragraph{Poisson equation}
\begin{equation}
-\ddudxx{u(x)} = f(x), \qquad u(0) = 0, \quad u(1) = 0
\label{deepmult::eq::poisson_eq}
\end{equation}
This equation gives the following matrix $A$: 
\begin{equation*}
A = -\frac{1}{h^2}
\begin{bmatrix}
-2  & 1 & ~ & ~ & ~ \\
1 & -2 & 1 & ~ & ~ \\
~ & \ddots & \ddots & \ddots \\
~ & ~ & 1 & -2 & 1 \\
~ & ~ & ~ & 1 & -2
\end{bmatrix}
\in \mathbb{R}^{n \times n}
\end{equation*}
\paragraph{Helmholtz equation}

\begin{equation}
-\ddudxx{u(x)} - k^2(x) u(x) = f(x), \qquad u(0) = 0, \quad u(1) = 0
\label{deepmult::eq::helm_eq}
\end{equation}
This equation gives the following matrix $A$:
\begin{equation*}
A = -\frac{1}{h^2}
\begin{bmatrix}
-2  & 1 & ~ & ~ & ~ \\
1 & -2 & 1 & ~ & ~ \\
~ & \ddots & \ddots & \ddots \\
~ & ~ & 1 & -2 & 1 \\
~ & ~ & ~ & 1 & -2
\end{bmatrix}
-
\mathrm{diag}(k^2(x_1), \ldots, k^2(x_n))
\in \mathbb{R}^{n \times n},
\end{equation*}
where $\mathrm{diag}(x)$ is a diagonal matrix with $x$ on the diagonal.
\paragraph{Singularly perturbed diffusion-convection equation} 
\begin{equation}
-\varepsilon \ddudxx{u(x)} + \dudx{u(x)} = f(x), \qquad u(0) = 0, \quad u(1) = 0
\label{deepmult::eq::diffusion_eq}
\end{equation}
This equation gives the following matrix $A$:
\begin{equation*}
A = -\frac{\varepsilon}{h^2}
\begin{bmatrix}
-2  & 1 & ~ & ~ & ~ \\
1 & -2 & 1 & ~ & ~ \\
~ & \ddots & \ddots & \ddots \\
~ & ~ & 1 & -2 & 1 \\
~ & ~ & ~ & 1 & -2
\end{bmatrix}
+ 
\frac{1}{h}
\begin{bmatrix}
-1  & 1 & ~ & ~ & ~ \\
0 & -1 & 1 & ~ & ~ \\
~ & \ddots & \ddots & \ddots \\
~ & ~ & \ddots & -1 & 1 \\
~ & ~ & ~ & 0  & -1
\end{bmatrix}
\in \mathbb{R}^{n \times n}
\end{equation*}

In this problem one can observe a boundary layer that has to be covered by the introduced grid~\cite{stynes2007convection}. 
Therefore, we get the constraint on the grid step $h$:
\[
h < \varepsilon.
\]  

The model problems give tridiagonal matrices $A$.
We also consider the restriction and prolongation operators that have the special structure which is shown below in the case of $n = 7$:
\[
R = 
\begin{bmatrix}
\times & \times & \times & & & & \\
& & \times & \times & \times & & \\
& & & & \times & \times & \times
\end{bmatrix} 
\quad
P = 
\begin{bmatrix}
\times & & \\
\times & & \\
\times & \times & \\
& \times & \\
& \times & \times \\
& & \times \\
& & \times
\end{bmatrix},
\]
where crosses $\times$ indicate non-zero elements. 
The standard choice for $R$ and $P$ comes from the linear interpolation between points in the fine grid and points in the coarse grid.
For example, for $n = 7$: 
\begin{equation}
R_{\mathrm{lin}} = 
\frac{1}{4}
\begin{bmatrix}
1 & 2 & 1\\
  &   & 1 & 2 & 1\\
  &   &   &   & 1 & 2 & 1
\end{bmatrix}
\quad
P_{\mathrm{lin}} = \frac{1}{2}
\begin{bmatrix}
1 & & \\
2 & & \\
1 & 1 & \\
  & 2 & \\
  & 1 & 1 \\
  &   & 2 \\
  &   & 1
\end{bmatrix}
\label{deepmult::eq::linear_inter}
\end{equation}

Note that the matrix $A$ is stored as an $n \times 3$ array, which stores only non-zero elements.
In the same way operators $R$ and $P$ are stored as $\frac{n-1}{2} \times 3$ arrays.
 

\subsection{Poisson equation}
We begin the study of the DMG performance with optimization of operators $R, P$ and damp factor $\omega$ for the Poisson equation~\eqref{deepmult::eq::poisson_eq}.
Initial approximations of optimized parameters are the linear interpolation in the form~\eqref{deepmult::eq::linear_inter} and $\omega_0 = \frac{2}{3}$.
We call linear interpolation operators and damp factor $\omega_0$ as \emph{linear parameters}.
In the same way, optimized operators and damp factor is called \emph{DMG parameters}.

Comparison $\rho(C)$ for linear parameters, AMG method and optimized ones is given in Table~\ref{deepmult::tab::poisson_rho}.
Parameters of the optimizer: $K = 10$, number of iterations $T=1000$, step size $s= 10^{-4}$ and batch size $N = 10$.
This table shows that DMG parameters provide smaller spectral radius $\rho(C)$ compared to linear parameters and to the AMG method~\cite{pyamg}.
We get better $\rho(C)$ but this gain becomes smaller with larger $N$.
 
\begin{table}[!ht]
\centering
\caption{Spectral radii $\rho$ for DMG, AMG and linear parameters in the case of the Poisson equation}
\begin{tabular}{cccc}
\toprule
 Grid size &    Linear &       AMG &       DMG \\
\midrule
         7 &  0.061728 &  0.182358 &  0.015088 \\
        15 &  0.061728 &  0.193726 &  0.018481 \\
        31 &  0.061728 &  0.196578 &  0.027819 \\
        63 &  0.061728 &  0.197207 &  0.045068 \\
       127 &  0.061728 &  0.195878 &  0.045400 \\
\bottomrule
\end{tabular}
\label{deepmult::tab::poisson_rho}
\end{table}

\subsection{Helmholtz equation}

For the Helmholtz equation we apply the proposed approach in three settings: low frequency ($k \approx 10$), high frequency ($k \approx 100$) and piecewise constant $k$:
\begin{equation}
k(x) = 
\begin{cases}
1, & 0 \leq x < 0.5\\
k_{\max}, & 0.5 \leq x \leq 1.
\end{cases}
\label{deepmult::eq::kx}
\end{equation}

In Table~\ref{deepmult::tab::helm_lowfreq} we present the performance of DMG in the low frequency case. 
It is shown that DMG parameters provide smaller spectral radius for different wavenumbers and different grid size.
In this case standard initialization with the linear parameters gives spectral radius less than one, therefore we do not need to use homotopy approach.
We use the following parameters for the optimizer: $K = 10$, number of iterations $T=1000$, step size $s= 10^{-4}$ and batch size $N = 10$.

\begin{table}[!ht]
\caption{Spectral radii $\rho$ in the case of the low frequency setting}
\centering
\begin{tabular}{ccccc}
\toprule
 Grid size &  $k$ &    Linear &       AMG &       DMG \\
\midrule
         7 &    5 &  0.226356 &  0.226214 &  0.012505 \\
        13 &   10 &  1.808608 &  0.255912 &  0.044337 \\
        17 &   15 &  0.826753 &  0.406821 &  0.062037 \\
        23 &   20 &  3.388036 &  0.418464 &  0.067183 \\
\bottomrule
\end{tabular}
\label{deepmult::tab::helm_lowfreq}
\end{table}

In Table~\ref{deepmult::tab::helm_highfreq} we provide results for the high frequency case of the Helmhotz equation.
In this experiment we have to use the homotopy approach for every considered $k$ with the homotopy step size $\delta = 0.1$.
Since we use the solution of the optimization problem corresponding to $k_i$ as an initialization for the optimization problem corresponding to $k_{i+1}$, we have to use the fixed grid size $n=1115$.
The parameters of the optimizer are the following: $K = 10$, the number of iterations $T=100$, step size $s= 10^{-4}$, batch size $N = 10$ and acceptance rate $\tau = 0.5$. 

\begin{table}[!ht]
\centering
\caption{Spectral radii $\rho$ in the case of the high frequency setting, grid size $n = 1115$}
\begin{tabular}{cccc}
\toprule
  $k$ &       Linear &       AMG &       DMG \\
\midrule
  100 &     0.180680 &  0.198093 &  0.061088 \\
  300 &    13.389492 &  0.203956 &  0.053827 \\
  500 &    14.608550 &  0.218872 &  0.066820 \\
  700 &    99.555631 &  0.243871 &  0.060205 \\
  900 &    62.940589 &  0.377024 &  0.091268 \\
 1000 &  4789.842424 &  0.607620 &  0.116077 \\
\bottomrule
\end{tabular}
\label{deepmult::tab::helm_highfreq}
\end{table}

In Table~\ref{deepmult::tab::helm_piecewise} we consider piecewise constant $k(x)$~\eqref{deepmult::eq::kx} for different grid sizes.
For all considered grid sizes DMG gives significantly better results.
Similarly to the Poisson equation, with larger grid size the difference between DMG and linear parameters is smaller.  

\begin{table}[!ht]
\centering
\caption{Spectral radii $\rho$ in the case of the piecewise constant $k(x)$~\eqref{deepmult::eq::kx}}
\begin{tabular}{ccccc}
\toprule
 Grid size &  $k_{\max}$ &    Linear &       AMG &       DMG \\
\midrule
       127 &         100 &  3.147622 &  0.330212 &  0.078162 \\
       255 &         100 &  1.642432 &  0.212405 &  0.047063 \\
       511 &         100 &  0.194238 &  0.200955 &  0.055769 \\
\bottomrule
\end{tabular}
\label{deepmult::tab::helm_piecewise}
\end{table}

\begin{figure}[!ht]
\centering
\resizebox{0.6\textwidth}{!}{\input{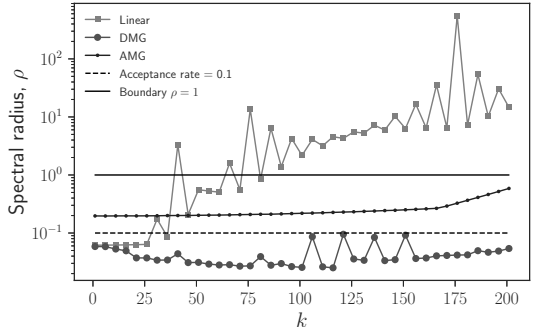}}
\caption{Spectral radius $\rho$ for DMG, AMG and linear parameters}
\label{deepmult::fig::helm_krange}
\end{figure}

Moreover, to observe how DMG performance changes from the low frequency setting to the high frequency one, we compute the spectral radius for $k \in [1, 200]$, the grid size $n = 225$ and plot it in Fig.~\ref{deepmult::fig::helm_krange}.
This plot shows that the proposed method provides parameters, which give significantly lower spectral radii in a wide range of frequencies.
To obtain this result we used the homotopy approach with the  acceptance rate $\tau = 0.1$ and the homotopy step size $\delta = 0.1$. 
Parameters of the optimizer are the following: $K = 10$, $N = 10$, $s = 10^{-4}$ and $T = 200$.
Initialization of parameters for the $k = 1$ is linear one.
The solution of the problem for $k = 1$ uses as initialization for the next $k$ and so on for the other~$k$'s.

\begin{figure}[!ht]
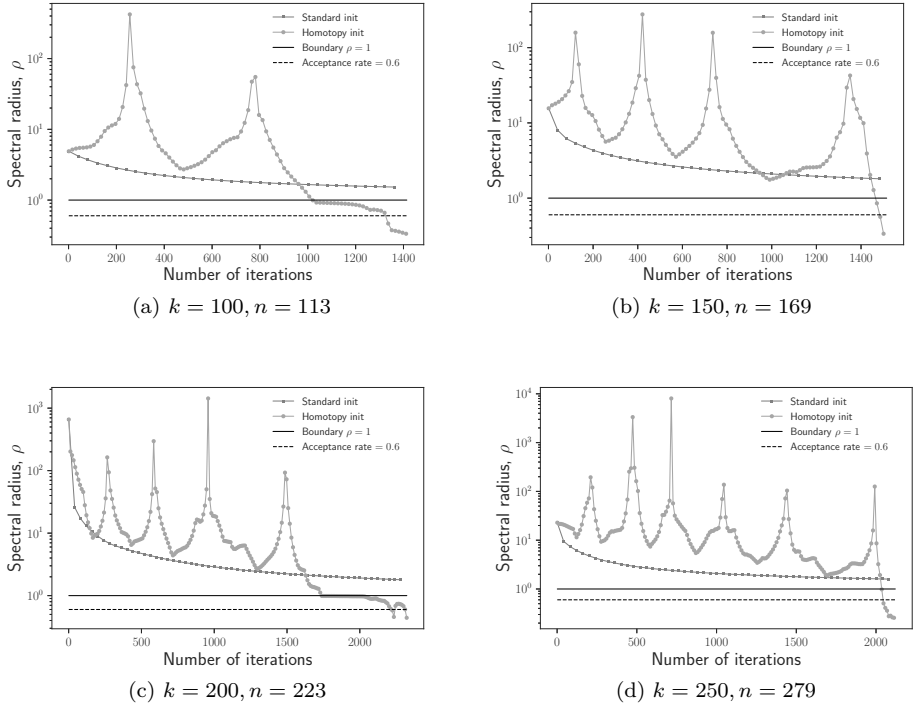

\centering
\subfloat[$k=100, n = 113$]{\label{deepmult::fig::helm_standard_hom_init_k_100}\resizebox{0.49\textwidth}{!}{\input{helm_standard_hom_init_m3points_k_100_seq_hom_bw.pgf}}}
\subfloat[$k=150, n = 169$]{\label{deepmult::fig::helm_standard_hom_init_k_150}\resizebox{0.49\textwidth}{!}{\input{helm_standard_hom_init_m3points_k_150_seq_hom_bw.pgf}}}
\\
\subfloat[$k=200, n = 223$]{\label{deepmult::fig::helm_standard_hom_init_k_200} \resizebox{0.49\textwidth}{!}{\input{helm_standard_hom_init_m3points_k_200_seq_hom_bw.pgf}}}
\subfloat[$k=250, n = 279$]{\label{deepmult::fig::laplace_random_linear_init_k_250} \resizebox{0.49\textwidth}{!}{\input{helm_standard_hom_init_m3points_k_250_seq_hom_bw.pgf}}}
\caption{Convergence of the DMG method for standard and homotopy initialization for different $k$}
\label{deepmult::fig::helm_standard_hom_init}
\end{figure}

Now we compare the convergence of the problem corresponding to the high frequency Helmholtz equation for homotopy and standard initialization.
In the case of the homotopy we use the matrix from the discretization of the Poisson equation as the matrix for which we know good starting guess.
This starting guess is the linear parameters.
We use the following parameters: the homotopy step $\delta = 0.1$, $K = 10$, $N = 10$, $s = 10^{-4}$.
The total number of iterations for both considered initialization approaches is chosen equally.
For every considered $k$ we set the grid size $n$ in a way that the number of points per wavelength equals $7$.
Fig.~\ref{deepmult::fig::helm_standard_hom_init} shows that homotopy gives the better local minimum than standard initialization, although the homotopy convergence is not monotone.
The reason of this non-monotonicity is that in the intermediate steps of the homotopy we solve optimization problem not for the target matrix, but for some auxiliary matrices.
However, the result of the homotopy approach is much better than for the standard initialization.

\subsection{Diffusion-convection equation}

\begin{figure}[!ht]
\centering
\subfloat[$n=63$]{\label{deepmult::fig::diff_epsrange_63}\resizebox{0.485\textwidth}{!}{\input{diff_m3points_epsrange_n_63_bw.pgf}}}
\subfloat[$n=127$]{\label{deepmult::fig::diff_epsrange_127}\resizebox{0.485\textwidth}{!}{\input{diff_m3points_epsrange_n_127_bw.pgf}}}
\\
\subfloat[$n=255$]{\label{deepmult::fig::diff_epsrange_255} \resizebox{0.485\textwidth}{!}{\input{diff_m3points_epsrange_n_255_bw.pgf}}}
\subfloat[$n=511$]{\label{deepmult::fig::diff_epsrange_511} \resizebox{0.485\textwidth}{!}{\input{diff_m3points_epsrange_n_511_bw.pgf}}}
\caption{Spectral radii for the DMG, AMG and linear parameters in the case of convection-diffusion equation}
\label{deepmult::fig::diff_epsrange}
\end{figure}
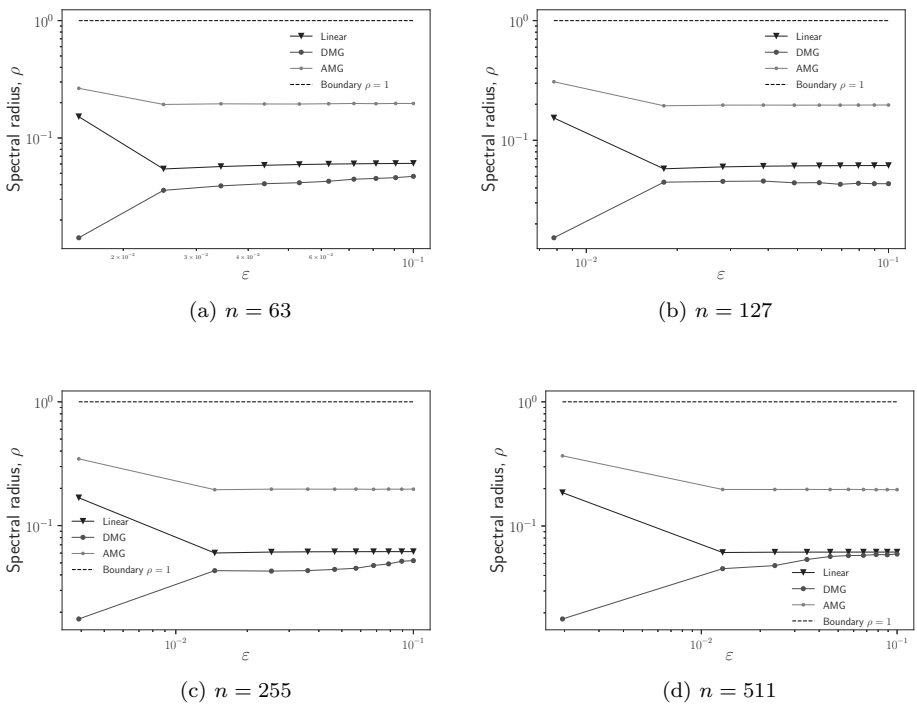

We consider Equation~\eqref{deepmult::eq::diffusion_eq} to demonstrate the performance of the proposed approach in the case of a non-symmetric matrix $A$.
For every grid size $n$ we study the spectral radii for various $\varepsilon \in [h, 10^{-1}]$.
The comparison of DMG parameters and linear parameters is shown in Fig.~\ref{deepmult::fig::diff_epsrange}.
This plot demonstrates that the DMG approach gives better results.
In this experiment we use standard initialization and initialize parameters with linear ones. 
The parameters of the optimizer are the following: $K = 10$, $T = 500$, $s = 10^{-4}$ and $N = 5$.

\section{Related work}
The problem of finding good restriction and prolongation operators in the multigrid method is studied a lot in the literature~\cite{yserentant1993old,trottenberg2000multigrid,hackbusch2013multi,hemker1990order,xu1996auxiliary}.
One of the most efficient ``black-box'' approaches is the black-box multigrid method, which is proposed in~\cite{de1990matrix,dendy1982black} and works for matrices coming from the discretization of 2D PDEs with $3 \times 3$ stencil.
Classical approach to the automatic selection of the operators in the multigrid method is the AMG 
~\cite{ruge1987algebraic,stuben2001review}, which we already considered in our numerical experiments. 
Although it often provides good convergence, geometric-type multigrid can give better results (as we also demonstrated for our problems).

Many methods were proposed to find appropriate operators in the multigrid method for particular types of PDEs.
For example,~\cite{stolk2014multigrid,livshits2014scalable} describe how to construct operators in the multigrid method for the Helmholtz equation.
The problem of finding operators of the multigrid method for the convection-diffusion is considered in~\cite{khelifi2014hybrid,dai2013fast}.
Although the proposed method is studied for structured grids, it can be extended to unstructured grids.
The multigrid method for these grids was studied in~\cite{grasedyck2016nearly}. 
All these methods are designed for specific problems, but clearly show, that operators that provide small spectral radius exist. 
Thus, our approach to optimize the spectral radius directly has a potential to find operators of the multigrid method for different problems.
Previously, the problem of spectral radius minimization was studied in~\cite{mengi2014numerical,nesterov2013optimizing,kressner2014subspace}, but the proposed methods require computation of the leading eigenvector, which is too costly.

The proposed DMG method is based on modern advances in deep neural networks~\cite{lecun2015deep,goodfellow2016deep} to fit very complex functions for classification~\cite{krizhevsky2012imagenet,socher2012convolutional}, computer vision~\cite{sun2014deep,kawulok2016advances}, and other tasks~\cite{yi2014deep,sosnovik2017neural,park2015deep,socher2013recursive}.
One of the papers devoted to image recognition problem~\cite{he2016deep} explores the idea of the multigrid method to motivate the deep neural network architecture.
Therefore, it is very natural to use the neural network approach to find local optimal parameters of the multigrid method itself.

\section{Conclusion and future work}
In this paper we show that the proposed DMG approach outperforms the linear interpolation operators for the one-dimensional Poisson, Helmholtz and singularly perturbed convection-diffusion equations.
However, the 1D case is purely academic and is aimed to prove the concept of the operators optimization to speed-up the  convergence.
More interesting and useful from practical viewpoint 2D and 3D cases will be considered in future work.
One of the promising ideas is to adopt the early stopping strategy in the optimization process, which will give trade off between the optimization of operators and the quality of the obtained local minimum.
Also the DMG approach solves a non-convex problem, therefore there are a lot of empirical parameters like batch size or power $K$ in Gelfand's formula approximation that affect the convergence and the quality of the local minimum.
Optimization of these parameters is an open question for further research.
Finally, the current approach requires learning of DMG parameters for every new matrix, which maybe impractical.
It would be very interesting to learn a mapping from the input parameters, i.e. coefficients of PDE, to the multigrid operators in a ``feed-forward'' manner.

\bibliographystyle{siamplain}
\bibliography{lib}

\end{document}

%% file: diff_m3points_epsrange_n_63_bw.pgf
\begingroup%
\makeatletter%
\begin{pgfpicture}%
\pgfpathrectangle{\pgfpointorigin}{\pgfqpoint{9.000000in}{6.000000in}}%
\pgfusepath{use as bounding box, clip}%
\begin{pgfscope}%
\pgfsetbuttcap%
\pgfsetmiterjoin%
\definecolor{currentfill}{rgb}{1.000000,1.000000,1.000000}%
\pgfsetfillcolor{currentfill}%
\pgfsetlinewidth{0.000000pt}%
\definecolor{currentstroke}{rgb}{1.000000,1.000000,1.000000}%
\pgfsetstrokecolor{currentstroke}%
\pgfsetdash{}{0pt}%
\pgfpathmoveto{\pgfqpoint{0.000000in}{0.000000in}}%
\pgfpathlineto{\pgfqpoint{9.000000in}{0.000000in}}%
\pgfpathlineto{\pgfqpoint{9.000000in}{6.000000in}}%
\pgfpathlineto{\pgfqpoint{0.000000in}{6.000000in}}%
\pgfpathclose%
\pgfusepath{fill}%
\end{pgfscope}%
\begin{pgfscope}%
\pgfsetbuttcap%
\pgfsetmiterjoin%
\definecolor{currentfill}{rgb}{1.000000,1.000000,1.000000}%
\pgfsetfillcolor{currentfill}%
\pgfsetlinewidth{0.000000pt}%
\definecolor{currentstroke}{rgb}{0.000000,0.000000,0.000000}%
\pgfsetstrokecolor{currentstroke}%
\pgfsetstrokeopacity{0.000000}%
\pgfsetdash{}{0pt}%
\pgfpathmoveto{\pgfqpoint{1.125000in}{0.750000in}}%
\pgfpathlineto{\pgfqpoint{8.100000in}{0.750000in}}%
\pgfpathlineto{\pgfqpoint{8.100000in}{5.280000in}}%
\pgfpathlineto{\pgfqpoint{1.125000in}{5.280000in}}%
\pgfpathclose%
\pgfusepath{fill}%
\end{pgfscope}%
\begin{pgfscope}%
\pgfsetbuttcap%
\pgfsetroundjoin%
\definecolor{currentfill}{rgb}{0.150000,0.150000,0.150000}%
\pgfsetfillcolor{currentfill}%
\pgfsetlinewidth{0.803000pt}%
\definecolor{currentstroke}{rgb}{0.150000,0.150000,0.150000}%
\pgfsetstrokecolor{currentstroke}%
\pgfsetdash{}{0pt}%
\pgfsys@defobject{currentmarker}{\pgfqpoint{0.000000in}{-0.083333in}}{\pgfqpoint{0.000000in}{0.000000in}}{%
\pgfpathmoveto{\pgfqpoint{0.000000in}{0.000000in}}%
\pgfpathlineto{\pgfqpoint{0.000000in}{-0.083333in}}%
\pgfusepath{stroke,fill}%
}%
\begin{pgfscope}%
\pgfsys@transformshift{7.782955in}{0.750000in}%
\pgfsys@useobject{currentmarker}{}%
\end{pgfscope}%
\end{pgfscope}%
\begin{pgfscope}%
\definecolor{textcolor}{rgb}{0.150000,0.150000,0.150000}%
\pgfsetstrokecolor{textcolor}%
\pgfsetfillcolor{textcolor}%
\pgftext[x=7.782955in,y=0.588889in,,top]{\color{textcolor}\sffamily\fontsize{16.000000}{19.200000}\selectfont \(\displaystyle 10^{-1}\)}%
\end{pgfscope}%
\begin{pgfscope}%
\pgfsetbuttcap%
\pgfsetroundjoin%
\definecolor{currentfill}{rgb}{0.150000,0.150000,0.150000}%
\pgfsetfillcolor{currentfill}%
\pgfsetlinewidth{0.401500pt}%
\definecolor{currentstroke}{rgb}{0.150000,0.150000,0.150000}%
\pgfsetstrokecolor{currentstroke}%
\pgfsetdash{}{0pt}%
\pgfsys@defobject{currentmarker}{\pgfqpoint{0.000000in}{-0.041667in}}{\pgfqpoint{0.000000in}{0.000000in}}{%
\pgfpathmoveto{\pgfqpoint{0.000000in}{0.000000in}}%
\pgfpathlineto{\pgfqpoint{0.000000in}{-0.041667in}}%
\pgfusepath{stroke,fill}%
}%
\begin{pgfscope}%
\pgfsys@transformshift{2.285292in}{0.750000in}%
\pgfsys@useobject{currentmarker}{}%
\end{pgfscope}%
\end{pgfscope}%
\begin{pgfscope}%
\definecolor{textcolor}{rgb}{0.150000,0.150000,0.150000}%
\pgfsetstrokecolor{textcolor}%
\pgfsetfillcolor{textcolor}%
\pgftext[x=2.285292in,y=0.661111in,,top]{\color{textcolor}\sffamily\fontsize{8.000000}{9.600000}\selectfont \(\displaystyle 2\times10^{-2}\)}%
\end{pgfscope}%
\begin{pgfscope}%
\pgfsetbuttcap%
\pgfsetroundjoin%
\definecolor{currentfill}{rgb}{0.150000,0.150000,0.150000}%
\pgfsetfillcolor{currentfill}%
\pgfsetlinewidth{0.401500pt}%
\definecolor{currentstroke}{rgb}{0.150000,0.150000,0.150000}%
\pgfsetstrokecolor{currentstroke}%
\pgfsetdash{}{0pt}%
\pgfsys@defobject{currentmarker}{\pgfqpoint{0.000000in}{-0.041667in}}{\pgfqpoint{0.000000in}{0.000000in}}{%
\pgfpathmoveto{\pgfqpoint{0.000000in}{0.000000in}}%
\pgfpathlineto{\pgfqpoint{0.000000in}{-0.041667in}}%
\pgfusepath{stroke,fill}%
}%
\begin{pgfscope}%
\pgfsys@transformshift{3.670316in}{0.750000in}%
\pgfsys@useobject{currentmarker}{}%
\end{pgfscope}%
\end{pgfscope}%
\begin{pgfscope}%
\definecolor{textcolor}{rgb}{0.150000,0.150000,0.150000}%
\pgfsetstrokecolor{textcolor}%
\pgfsetfillcolor{textcolor}%
\pgftext[x=3.670316in,y=0.661111in,,top]{\color{textcolor}\sffamily\fontsize{8.000000}{9.600000}\selectfont \(\displaystyle 3\times10^{-2}\)}%
\end{pgfscope}%
\begin{pgfscope}%
\pgfsetbuttcap%
\pgfsetroundjoin%
\definecolor{currentfill}{rgb}{0.150000,0.150000,0.150000}%
\pgfsetfillcolor{currentfill}%
\pgfsetlinewidth{0.401500pt}%
\definecolor{currentstroke}{rgb}{0.150000,0.150000,0.150000}%
\pgfsetstrokecolor{currentstroke}%
\pgfsetdash{}{0pt}%
\pgfsys@defobject{currentmarker}{\pgfqpoint{0.000000in}{-0.041667in}}{\pgfqpoint{0.000000in}{0.000000in}}{%
\pgfpathmoveto{\pgfqpoint{0.000000in}{0.000000in}}%
\pgfpathlineto{\pgfqpoint{0.000000in}{-0.041667in}}%
\pgfusepath{stroke,fill}%
}%
\begin{pgfscope}%
\pgfsys@transformshift{4.653007in}{0.750000in}%
\pgfsys@useobject{currentmarker}{}%
\end{pgfscope}%
\end{pgfscope}%
\begin{pgfscope}%
\definecolor{textcolor}{rgb}{0.150000,0.150000,0.150000}%
\pgfsetstrokecolor{textcolor}%
\pgfsetfillcolor{textcolor}%
\pgftext[x=4.653007in,y=0.661111in,,top]{\color{textcolor}\sffamily\fontsize{8.000000}{9.600000}\selectfont \(\displaystyle 4\times10^{-2}\)}%
\end{pgfscope}%
\begin{pgfscope}%
\pgfsetbuttcap%
\pgfsetroundjoin%
\definecolor{currentfill}{rgb}{0.150000,0.150000,0.150000}%
\pgfsetfillcolor{currentfill}%
\pgfsetlinewidth{0.401500pt}%
\definecolor{currentstroke}{rgb}{0.150000,0.150000,0.150000}%
\pgfsetstrokecolor{currentstroke}%
\pgfsetdash{}{0pt}%
\pgfsys@defobject{currentmarker}{\pgfqpoint{0.000000in}{-0.041667in}}{\pgfqpoint{0.000000in}{0.000000in}}{%
\pgfpathmoveto{\pgfqpoint{0.000000in}{0.000000in}}%
\pgfpathlineto{\pgfqpoint{0.000000in}{-0.041667in}}%
\pgfusepath{stroke,fill}%
}%
\begin{pgfscope}%
\pgfsys@transformshift{5.415240in}{0.750000in}%
\pgfsys@useobject{currentmarker}{}%
\end{pgfscope}%
\end{pgfscope}%
\begin{pgfscope}%
\pgfsetbuttcap%
\pgfsetroundjoin%
\definecolor{currentfill}{rgb}{0.150000,0.150000,0.150000}%
\pgfsetfillcolor{currentfill}%
\pgfsetlinewidth{0.401500pt}%
\definecolor{currentstroke}{rgb}{0.150000,0.150000,0.150000}%
\pgfsetstrokecolor{currentstroke}%
\pgfsetdash{}{0pt}%
\pgfsys@defobject{currentmarker}{\pgfqpoint{0.000000in}{-0.041667in}}{\pgfqpoint{0.000000in}{0.000000in}}{%
\pgfpathmoveto{\pgfqpoint{0.000000in}{0.000000in}}%
\pgfpathlineto{\pgfqpoint{0.000000in}{-0.041667in}}%
\pgfusepath{stroke,fill}%
}%
\begin{pgfscope}%
\pgfsys@transformshift{6.038031in}{0.750000in}%
\pgfsys@useobject{currentmarker}{}%
\end{pgfscope}%
\end{pgfscope}%
\begin{pgfscope}%
\definecolor{textcolor}{rgb}{0.150000,0.150000,0.150000}%
\pgfsetstrokecolor{textcolor}%
\pgfsetfillcolor{textcolor}%
\pgftext[x=6.038031in,y=0.661111in,,top]{\color{textcolor}\sffamily\fontsize{8.000000}{9.600000}\selectfont \(\displaystyle 6\times10^{-2}\)}%
\end{pgfscope}%
\begin{pgfscope}%
\pgfsetbuttcap%
\pgfsetroundjoin%
\definecolor{currentfill}{rgb}{0.150000,0.150000,0.150000}%
\pgfsetfillcolor{currentfill}%
\pgfsetlinewidth{0.401500pt}%
\definecolor{currentstroke}{rgb}{0.150000,0.150000,0.150000}%
\pgfsetstrokecolor{currentstroke}%
\pgfsetdash{}{0pt}%
\pgfsys@defobject{currentmarker}{\pgfqpoint{0.000000in}{-0.041667in}}{\pgfqpoint{0.000000in}{0.000000in}}{%
\pgfpathmoveto{\pgfqpoint{0.000000in}{0.000000in}}%
\pgfpathlineto{\pgfqpoint{0.000000in}{-0.041667in}}%
\pgfusepath{stroke,fill}%
}%
\begin{pgfscope}%
\pgfsys@transformshift{6.564592in}{0.750000in}%
\pgfsys@useobject{currentmarker}{}%
\end{pgfscope}%
\end{pgfscope}%
\begin{pgfscope}%
\pgfsetbuttcap%
\pgfsetroundjoin%
\definecolor{currentfill}{rgb}{0.150000,0.150000,0.150000}%
\pgfsetfillcolor{currentfill}%
\pgfsetlinewidth{0.401500pt}%
\definecolor{currentstroke}{rgb}{0.150000,0.150000,0.150000}%
\pgfsetstrokecolor{currentstroke}%
\pgfsetdash{}{0pt}%
\pgfsys@defobject{currentmarker}{\pgfqpoint{0.000000in}{-0.041667in}}{\pgfqpoint{0.000000in}{0.000000in}}{%
\pgfpathmoveto{\pgfqpoint{0.000000in}{0.000000in}}%
\pgfpathlineto{\pgfqpoint{0.000000in}{-0.041667in}}%
\pgfusepath{stroke,fill}%
}%
\begin{pgfscope}%
\pgfsys@transformshift{7.020721in}{0.750000in}%
\pgfsys@useobject{currentmarker}{}%
\end{pgfscope}%
\end{pgfscope}%
\begin{pgfscope}%
\pgfsetbuttcap%
\pgfsetroundjoin%
\definecolor{currentfill}{rgb}{0.150000,0.150000,0.150000}%
\pgfsetfillcolor{currentfill}%
\pgfsetlinewidth{0.401500pt}%
\definecolor{currentstroke}{rgb}{0.150000,0.150000,0.150000}%
\pgfsetstrokecolor{currentstroke}%
\pgfsetdash{}{0pt}%
\pgfsys@defobject{currentmarker}{\pgfqpoint{0.000000in}{-0.041667in}}{\pgfqpoint{0.000000in}{0.000000in}}{%
\pgfpathmoveto{\pgfqpoint{0.000000in}{0.000000in}}%
\pgfpathlineto{\pgfqpoint{0.000000in}{-0.041667in}}%
\pgfusepath{stroke,fill}%
}%
\begin{pgfscope}%
\pgfsys@transformshift{7.423055in}{0.750000in}%
\pgfsys@useobject{currentmarker}{}%
\end{pgfscope}%
\end{pgfscope}%
\begin{pgfscope}%
\definecolor{textcolor}{rgb}{0.150000,0.150000,0.150000}%
\pgfsetstrokecolor{textcolor}%
\pgfsetfillcolor{textcolor}%
\pgftext[x=4.612500in,y=0.318273in,,top]{\color{textcolor}\sffamily\fontsize{24.000000}{28.800000}\selectfont \(\displaystyle \varepsilon\)}%
\end{pgfscope}%
\begin{pgfscope}%
\pgfsetbuttcap%
\pgfsetroundjoin%
\definecolor{currentfill}{rgb}{0.150000,0.150000,0.150000}%
\pgfsetfillcolor{currentfill}%
\pgfsetlinewidth{0.803000pt}%
\definecolor{currentstroke}{rgb}{0.150000,0.150000,0.150000}%
\pgfsetstrokecolor{currentstroke}%
\pgfsetdash{}{0pt}%
\pgfsys@defobject{currentmarker}{\pgfqpoint{-0.083333in}{0.000000in}}{\pgfqpoint{0.000000in}{0.000000in}}{%
\pgfpathmoveto{\pgfqpoint{0.000000in}{0.000000in}}%
\pgfpathlineto{\pgfqpoint{-0.083333in}{0.000000in}}%
\pgfusepath{stroke,fill}%
}%
\begin{pgfscope}%
\pgfsys@transformshift{1.125000in}{2.848437in}%
\pgfsys@useobject{currentmarker}{}%
\end{pgfscope}%
\end{pgfscope}%
\begin{pgfscope}%
\definecolor{textcolor}{rgb}{0.150000,0.150000,0.150000}%
\pgfsetstrokecolor{textcolor}%
\pgfsetfillcolor{textcolor}%
\pgftext[x=0.439258in,y=2.732361in,left,base]{\color{textcolor}\sffamily\fontsize{22.000000}{26.400000}\selectfont \(\displaystyle 10^{-1}\)}%
\end{pgfscope}%
\begin{pgfscope}%
\pgfsetbuttcap%
\pgfsetroundjoin%
\definecolor{currentfill}{rgb}{0.150000,0.150000,0.150000}%
\pgfsetfillcolor{currentfill}%
\pgfsetlinewidth{0.803000pt}%
\definecolor{currentstroke}{rgb}{0.150000,0.150000,0.150000}%
\pgfsetstrokecolor{currentstroke}%
\pgfsetdash{}{0pt}%
\pgfsys@defobject{currentmarker}{\pgfqpoint{-0.083333in}{0.000000in}}{\pgfqpoint{0.000000in}{0.000000in}}{%
\pgfpathmoveto{\pgfqpoint{0.000000in}{0.000000in}}%
\pgfpathlineto{\pgfqpoint{-0.083333in}{0.000000in}}%
\pgfusepath{stroke,fill}%
}%
\begin{pgfscope}%
\pgfsys@transformshift{1.125000in}{5.074091in}%
\pgfsys@useobject{currentmarker}{}%
\end{pgfscope}%
\end{pgfscope}%
\begin{pgfscope}%
\definecolor{textcolor}{rgb}{0.150000,0.150000,0.150000}%
\pgfsetstrokecolor{textcolor}%
\pgfsetfillcolor{textcolor}%
\pgftext[x=0.594814in,y=4.958016in,left,base]{\color{textcolor}\sffamily\fontsize{22.000000}{26.400000}\selectfont \(\displaystyle 10^{0}\)}%
\end{pgfscope}%
\begin{pgfscope}%
\pgfsetbuttcap%
\pgfsetroundjoin%
\definecolor{currentfill}{rgb}{0.150000,0.150000,0.150000}%
\pgfsetfillcolor{currentfill}%
\pgfsetlinewidth{0.401500pt}%
\definecolor{currentstroke}{rgb}{0.150000,0.150000,0.150000}%
\pgfsetstrokecolor{currentstroke}%
\pgfsetdash{}{0pt}%
\pgfsys@defobject{currentmarker}{\pgfqpoint{-0.041667in}{0.000000in}}{\pgfqpoint{0.000000in}{0.000000in}}{%
\pgfpathmoveto{\pgfqpoint{0.000000in}{0.000000in}}%
\pgfpathlineto{\pgfqpoint{-0.041667in}{0.000000in}}%
\pgfusepath{stroke,fill}%
}%
\begin{pgfscope}%
\pgfsys@transformshift{1.125000in}{1.292771in}%
\pgfsys@useobject{currentmarker}{}%
\end{pgfscope}%
\end{pgfscope}%
\begin{pgfscope}%
\pgfsetbuttcap%
\pgfsetroundjoin%
\definecolor{currentfill}{rgb}{0.150000,0.150000,0.150000}%
\pgfsetfillcolor{currentfill}%
\pgfsetlinewidth{0.401500pt}%
\definecolor{currentstroke}{rgb}{0.150000,0.150000,0.150000}%
\pgfsetstrokecolor{currentstroke}%
\pgfsetdash{}{0pt}%
\pgfsys@defobject{currentmarker}{\pgfqpoint{-0.041667in}{0.000000in}}{\pgfqpoint{0.000000in}{0.000000in}}{%
\pgfpathmoveto{\pgfqpoint{0.000000in}{0.000000in}}%
\pgfpathlineto{\pgfqpoint{-0.041667in}{0.000000in}}%
\pgfusepath{stroke,fill}%
}%
\begin{pgfscope}%
\pgfsys@transformshift{1.125000in}{1.684689in}%
\pgfsys@useobject{currentmarker}{}%
\end{pgfscope}%
\end{pgfscope}%
\begin{pgfscope}%
\pgfsetbuttcap%
\pgfsetroundjoin%
\definecolor{currentfill}{rgb}{0.150000,0.150000,0.150000}%
\pgfsetfillcolor{currentfill}%
\pgfsetlinewidth{0.401500pt}%
\definecolor{currentstroke}{rgb}{0.150000,0.150000,0.150000}%
\pgfsetstrokecolor{currentstroke}%
\pgfsetdash{}{0pt}%
\pgfsys@defobject{currentmarker}{\pgfqpoint{-0.041667in}{0.000000in}}{\pgfqpoint{0.000000in}{0.000000in}}{%
\pgfpathmoveto{\pgfqpoint{0.000000in}{0.000000in}}%
\pgfpathlineto{\pgfqpoint{-0.041667in}{0.000000in}}%
\pgfusepath{stroke,fill}%
}%
\begin{pgfscope}%
\pgfsys@transformshift{1.125000in}{1.962760in}%
\pgfsys@useobject{currentmarker}{}%
\end{pgfscope}%
\end{pgfscope}%
\begin{pgfscope}%
\pgfsetbuttcap%
\pgfsetroundjoin%
\definecolor{currentfill}{rgb}{0.150000,0.150000,0.150000}%
\pgfsetfillcolor{currentfill}%
\pgfsetlinewidth{0.401500pt}%
\definecolor{currentstroke}{rgb}{0.150000,0.150000,0.150000}%
\pgfsetstrokecolor{currentstroke}%
\pgfsetdash{}{0pt}%
\pgfsys@defobject{currentmarker}{\pgfqpoint{-0.041667in}{0.000000in}}{\pgfqpoint{0.000000in}{0.000000in}}{%
\pgfpathmoveto{\pgfqpoint{0.000000in}{0.000000in}}%
\pgfpathlineto{\pgfqpoint{-0.041667in}{0.000000in}}%
\pgfusepath{stroke,fill}%
}%
\begin{pgfscope}%
\pgfsys@transformshift{1.125000in}{2.178448in}%
\pgfsys@useobject{currentmarker}{}%
\end{pgfscope}%
\end{pgfscope}%
\begin{pgfscope}%
\pgfsetbuttcap%
\pgfsetroundjoin%
\definecolor{currentfill}{rgb}{0.150000,0.150000,0.150000}%
\pgfsetfillcolor{currentfill}%
\pgfsetlinewidth{0.401500pt}%
\definecolor{currentstroke}{rgb}{0.150000,0.150000,0.150000}%
\pgfsetstrokecolor{currentstroke}%
\pgfsetdash{}{0pt}%
\pgfsys@defobject{currentmarker}{\pgfqpoint{-0.041667in}{0.000000in}}{\pgfqpoint{0.000000in}{0.000000in}}{%
\pgfpathmoveto{\pgfqpoint{0.000000in}{0.000000in}}%
\pgfpathlineto{\pgfqpoint{-0.041667in}{0.000000in}}%
\pgfusepath{stroke,fill}%
}%
\begin{pgfscope}%
\pgfsys@transformshift{1.125000in}{2.354678in}%
\pgfsys@useobject{currentmarker}{}%
\end{pgfscope}%
\end{pgfscope}%
\begin{pgfscope}%
\pgfsetbuttcap%
\pgfsetroundjoin%
\definecolor{currentfill}{rgb}{0.150000,0.150000,0.150000}%
\pgfsetfillcolor{currentfill}%
\pgfsetlinewidth{0.401500pt}%
\definecolor{currentstroke}{rgb}{0.150000,0.150000,0.150000}%
\pgfsetstrokecolor{currentstroke}%
\pgfsetdash{}{0pt}%
\pgfsys@defobject{currentmarker}{\pgfqpoint{-0.041667in}{0.000000in}}{\pgfqpoint{0.000000in}{0.000000in}}{%
\pgfpathmoveto{\pgfqpoint{0.000000in}{0.000000in}}%
\pgfpathlineto{\pgfqpoint{-0.041667in}{0.000000in}}%
\pgfusepath{stroke,fill}%
}%
\begin{pgfscope}%
\pgfsys@transformshift{1.125000in}{2.503678in}%
\pgfsys@useobject{currentmarker}{}%
\end{pgfscope}%
\end{pgfscope}%
\begin{pgfscope}%
\pgfsetbuttcap%
\pgfsetroundjoin%
\definecolor{currentfill}{rgb}{0.150000,0.150000,0.150000}%
\pgfsetfillcolor{currentfill}%
\pgfsetlinewidth{0.401500pt}%
\definecolor{currentstroke}{rgb}{0.150000,0.150000,0.150000}%
\pgfsetstrokecolor{currentstroke}%
\pgfsetdash{}{0pt}%
\pgfsys@defobject{currentmarker}{\pgfqpoint{-0.041667in}{0.000000in}}{\pgfqpoint{0.000000in}{0.000000in}}{%
\pgfpathmoveto{\pgfqpoint{0.000000in}{0.000000in}}%
\pgfpathlineto{\pgfqpoint{-0.041667in}{0.000000in}}%
\pgfusepath{stroke,fill}%
}%
\begin{pgfscope}%
\pgfsys@transformshift{1.125000in}{2.632748in}%
\pgfsys@useobject{currentmarker}{}%
\end{pgfscope}%
\end{pgfscope}%
\begin{pgfscope}%
\pgfsetbuttcap%
\pgfsetroundjoin%
\definecolor{currentfill}{rgb}{0.150000,0.150000,0.150000}%
\pgfsetfillcolor{currentfill}%
\pgfsetlinewidth{0.401500pt}%
\definecolor{currentstroke}{rgb}{0.150000,0.150000,0.150000}%
\pgfsetstrokecolor{currentstroke}%
\pgfsetdash{}{0pt}%
\pgfsys@defobject{currentmarker}{\pgfqpoint{-0.041667in}{0.000000in}}{\pgfqpoint{0.000000in}{0.000000in}}{%
\pgfpathmoveto{\pgfqpoint{0.000000in}{0.000000in}}%
\pgfpathlineto{\pgfqpoint{-0.041667in}{0.000000in}}%
\pgfusepath{stroke,fill}%
}%
\begin{pgfscope}%
\pgfsys@transformshift{1.125000in}{2.746596in}%
\pgfsys@useobject{currentmarker}{}%
\end{pgfscope}%
\end{pgfscope}%
\begin{pgfscope}%
\pgfsetbuttcap%
\pgfsetroundjoin%
\definecolor{currentfill}{rgb}{0.150000,0.150000,0.150000}%
\pgfsetfillcolor{currentfill}%
\pgfsetlinewidth{0.401500pt}%
\definecolor{currentstroke}{rgb}{0.150000,0.150000,0.150000}%
\pgfsetstrokecolor{currentstroke}%
\pgfsetdash{}{0pt}%
\pgfsys@defobject{currentmarker}{\pgfqpoint{-0.041667in}{0.000000in}}{\pgfqpoint{0.000000in}{0.000000in}}{%
\pgfpathmoveto{\pgfqpoint{0.000000in}{0.000000in}}%
\pgfpathlineto{\pgfqpoint{-0.041667in}{0.000000in}}%
\pgfusepath{stroke,fill}%
}%
\begin{pgfscope}%
\pgfsys@transformshift{1.125000in}{3.518425in}%
\pgfsys@useobject{currentmarker}{}%
\end{pgfscope}%
\end{pgfscope}%
\begin{pgfscope}%
\pgfsetbuttcap%
\pgfsetroundjoin%
\definecolor{currentfill}{rgb}{0.150000,0.150000,0.150000}%
\pgfsetfillcolor{currentfill}%
\pgfsetlinewidth{0.401500pt}%
\definecolor{currentstroke}{rgb}{0.150000,0.150000,0.150000}%
\pgfsetstrokecolor{currentstroke}%
\pgfsetdash{}{0pt}%
\pgfsys@defobject{currentmarker}{\pgfqpoint{-0.041667in}{0.000000in}}{\pgfqpoint{0.000000in}{0.000000in}}{%
\pgfpathmoveto{\pgfqpoint{0.000000in}{0.000000in}}%
\pgfpathlineto{\pgfqpoint{-0.041667in}{0.000000in}}%
\pgfusepath{stroke,fill}%
}%
\begin{pgfscope}%
\pgfsys@transformshift{1.125000in}{3.910344in}%
\pgfsys@useobject{currentmarker}{}%
\end{pgfscope}%
\end{pgfscope}%
\begin{pgfscope}%
\pgfsetbuttcap%
\pgfsetroundjoin%
\definecolor{currentfill}{rgb}{0.150000,0.150000,0.150000}%
\pgfsetfillcolor{currentfill}%
\pgfsetlinewidth{0.401500pt}%
\definecolor{currentstroke}{rgb}{0.150000,0.150000,0.150000}%
\pgfsetstrokecolor{currentstroke}%
\pgfsetdash{}{0pt}%
\pgfsys@defobject{currentmarker}{\pgfqpoint{-0.041667in}{0.000000in}}{\pgfqpoint{0.000000in}{0.000000in}}{%
\pgfpathmoveto{\pgfqpoint{0.000000in}{0.000000in}}%
\pgfpathlineto{\pgfqpoint{-0.041667in}{0.000000in}}%
\pgfusepath{stroke,fill}%
}%
\begin{pgfscope}%
\pgfsys@transformshift{1.125000in}{4.188414in}%
\pgfsys@useobject{currentmarker}{}%
\end{pgfscope}%
\end{pgfscope}%
\begin{pgfscope}%
\pgfsetbuttcap%
\pgfsetroundjoin%
\definecolor{currentfill}{rgb}{0.150000,0.150000,0.150000}%
\pgfsetfillcolor{currentfill}%
\pgfsetlinewidth{0.401500pt}%
\definecolor{currentstroke}{rgb}{0.150000,0.150000,0.150000}%
\pgfsetstrokecolor{currentstroke}%
\pgfsetdash{}{0pt}%
\pgfsys@defobject{currentmarker}{\pgfqpoint{-0.041667in}{0.000000in}}{\pgfqpoint{0.000000in}{0.000000in}}{%
\pgfpathmoveto{\pgfqpoint{0.000000in}{0.000000in}}%
\pgfpathlineto{\pgfqpoint{-0.041667in}{0.000000in}}%
\pgfusepath{stroke,fill}%
}%
\begin{pgfscope}%
\pgfsys@transformshift{1.125000in}{4.404102in}%
\pgfsys@useobject{currentmarker}{}%
\end{pgfscope}%
\end{pgfscope}%
\begin{pgfscope}%
\pgfsetbuttcap%
\pgfsetroundjoin%
\definecolor{currentfill}{rgb}{0.150000,0.150000,0.150000}%
\pgfsetfillcolor{currentfill}%
\pgfsetlinewidth{0.401500pt}%
\definecolor{currentstroke}{rgb}{0.150000,0.150000,0.150000}%
\pgfsetstrokecolor{currentstroke}%
\pgfsetdash{}{0pt}%
\pgfsys@defobject{currentmarker}{\pgfqpoint{-0.041667in}{0.000000in}}{\pgfqpoint{0.000000in}{0.000000in}}{%
\pgfpathmoveto{\pgfqpoint{0.000000in}{0.000000in}}%
\pgfpathlineto{\pgfqpoint{-0.041667in}{0.000000in}}%
\pgfusepath{stroke,fill}%
}%
\begin{pgfscope}%
\pgfsys@transformshift{1.125000in}{4.580332in}%
\pgfsys@useobject{currentmarker}{}%
\end{pgfscope}%
\end{pgfscope}%
\begin{pgfscope}%
\pgfsetbuttcap%
\pgfsetroundjoin%
\definecolor{currentfill}{rgb}{0.150000,0.150000,0.150000}%
\pgfsetfillcolor{currentfill}%
\pgfsetlinewidth{0.401500pt}%
\definecolor{currentstroke}{rgb}{0.150000,0.150000,0.150000}%
\pgfsetstrokecolor{currentstroke}%
\pgfsetdash{}{0pt}%
\pgfsys@defobject{currentmarker}{\pgfqpoint{-0.041667in}{0.000000in}}{\pgfqpoint{0.000000in}{0.000000in}}{%
\pgfpathmoveto{\pgfqpoint{0.000000in}{0.000000in}}%
\pgfpathlineto{\pgfqpoint{-0.041667in}{0.000000in}}%
\pgfusepath{stroke,fill}%
}%
\begin{pgfscope}%
\pgfsys@transformshift{1.125000in}{4.729333in}%
\pgfsys@useobject{currentmarker}{}%
\end{pgfscope}%
\end{pgfscope}%
\begin{pgfscope}%
\pgfsetbuttcap%
\pgfsetroundjoin%
\definecolor{currentfill}{rgb}{0.150000,0.150000,0.150000}%
\pgfsetfillcolor{currentfill}%
\pgfsetlinewidth{0.401500pt}%
\definecolor{currentstroke}{rgb}{0.150000,0.150000,0.150000}%
\pgfsetstrokecolor{currentstroke}%
\pgfsetdash{}{0pt}%
\pgfsys@defobject{currentmarker}{\pgfqpoint{-0.041667in}{0.000000in}}{\pgfqpoint{0.000000in}{0.000000in}}{%
\pgfpathmoveto{\pgfqpoint{0.000000in}{0.000000in}}%
\pgfpathlineto{\pgfqpoint{-0.041667in}{0.000000in}}%
\pgfusepath{stroke,fill}%
}%
\begin{pgfscope}%
\pgfsys@transformshift{1.125000in}{4.858403in}%
\pgfsys@useobject{currentmarker}{}%
\end{pgfscope}%
\end{pgfscope}%
\begin{pgfscope}%
\pgfsetbuttcap%
\pgfsetroundjoin%
\definecolor{currentfill}{rgb}{0.150000,0.150000,0.150000}%
\pgfsetfillcolor{currentfill}%
\pgfsetlinewidth{0.401500pt}%
\definecolor{currentstroke}{rgb}{0.150000,0.150000,0.150000}%
\pgfsetstrokecolor{currentstroke}%
\pgfsetdash{}{0pt}%
\pgfsys@defobject{currentmarker}{\pgfqpoint{-0.041667in}{0.000000in}}{\pgfqpoint{0.000000in}{0.000000in}}{%
\pgfpathmoveto{\pgfqpoint{0.000000in}{0.000000in}}%
\pgfpathlineto{\pgfqpoint{-0.041667in}{0.000000in}}%
\pgfusepath{stroke,fill}%
}%
\begin{pgfscope}%
\pgfsys@transformshift{1.125000in}{4.972251in}%
\pgfsys@useobject{currentmarker}{}%
\end{pgfscope}%
\end{pgfscope}%
\begin{pgfscope}%
\definecolor{textcolor}{rgb}{0.150000,0.150000,0.150000}%
\pgfsetstrokecolor{textcolor}%
\pgfsetfillcolor{textcolor}%
\pgftext[x=0.383703in,y=3.015000in,,bottom,rotate=90.000000]{\color{textcolor}\sffamily\fontsize{24.000000}{28.800000}\selectfont Spectral radius, \(\displaystyle \rho\)}%
\end{pgfscope}%
\begin{pgfscope}%
\pgfpathrectangle{\pgfqpoint{1.125000in}{0.750000in}}{\pgfqpoint{6.975000in}{4.530000in}} %
\pgfusepath{clip}%
\pgfsetroundcap%
\pgfsetroundjoin%
\pgfsetlinewidth{1.405250pt}%
\definecolor{currentstroke}{rgb}{0.133333,0.133333,0.133333}%
\pgfsetstrokecolor{currentstroke}%
\pgfsetdash{}{0pt}%
\pgfpathmoveto{\pgfqpoint{1.442045in}{3.252018in}}%
\pgfpathlineto{\pgfqpoint{3.047526in}{2.260870in}}%
\pgfpathlineto{\pgfqpoint{4.135329in}{2.306871in}}%
\pgfpathlineto{\pgfqpoint{4.959112in}{2.331252in}}%
\pgfpathlineto{\pgfqpoint{5.622327in}{2.344954in}}%
\pgfpathlineto{\pgfqpoint{6.177474in}{2.353311in}}%
\pgfpathlineto{\pgfqpoint{6.654885in}{2.358737in}}%
\pgfpathlineto{\pgfqpoint{7.073681in}{2.362436in}}%
\pgfpathlineto{\pgfqpoint{7.446694in}{2.365060in}}%
\pgfpathlineto{\pgfqpoint{7.782955in}{2.366984in}}%
\pgfusepath{stroke}%
\end{pgfscope}%
\begin{pgfscope}%
\pgfpathrectangle{\pgfqpoint{1.125000in}{0.750000in}}{\pgfqpoint{6.975000in}{4.530000in}} %
\pgfusepath{clip}%
\pgfsetbuttcap%
\pgfsetmiterjoin%
\definecolor{currentfill}{rgb}{0.133333,0.133333,0.133333}%
\pgfsetfillcolor{currentfill}%
\pgfsetlinewidth{0.000000pt}%
\definecolor{currentstroke}{rgb}{0.133333,0.133333,0.133333}%
\pgfsetstrokecolor{currentstroke}%
\pgfsetdash{}{0pt}%
\pgfsys@defobject{currentmarker}{\pgfqpoint{-0.062500in}{-0.062500in}}{\pgfqpoint{0.062500in}{0.062500in}}{%
\pgfpathmoveto{\pgfqpoint{-0.000000in}{-0.062500in}}%
\pgfpathlineto{\pgfqpoint{0.062500in}{0.062500in}}%
\pgfpathlineto{\pgfqpoint{-0.062500in}{0.062500in}}%
\pgfpathclose%
\pgfusepath{fill}%
}%
\begin{pgfscope}%
\pgfsys@transformshift{1.442045in}{3.252018in}%
\pgfsys@useobject{currentmarker}{}%
\end{pgfscope}%
\begin{pgfscope}%
\pgfsys@transformshift{3.047526in}{2.260870in}%
\pgfsys@useobject{currentmarker}{}%
\end{pgfscope}%
\begin{pgfscope}%
\pgfsys@transformshift{4.135329in}{2.306871in}%
\pgfsys@useobject{currentmarker}{}%
\end{pgfscope}%
\begin{pgfscope}%
\pgfsys@transformshift{4.959112in}{2.331252in}%
\pgfsys@useobject{currentmarker}{}%
\end{pgfscope}%
\begin{pgfscope}%
\pgfsys@transformshift{5.622327in}{2.344954in}%
\pgfsys@useobject{currentmarker}{}%
\end{pgfscope}%
\begin{pgfscope}%
\pgfsys@transformshift{6.177474in}{2.353311in}%
\pgfsys@useobject{currentmarker}{}%
\end{pgfscope}%
\begin{pgfscope}%
\pgfsys@transformshift{6.654885in}{2.358737in}%
\pgfsys@useobject{currentmarker}{}%
\end{pgfscope}%
\begin{pgfscope}%
\pgfsys@transformshift{7.073681in}{2.362436in}%
\pgfsys@useobject{currentmarker}{}%
\end{pgfscope}%
\begin{pgfscope}%
\pgfsys@transformshift{7.446694in}{2.365060in}%
\pgfsys@useobject{currentmarker}{}%
\end{pgfscope}%
\begin{pgfscope}%
\pgfsys@transformshift{7.782955in}{2.366984in}%
\pgfsys@useobject{currentmarker}{}%
\end{pgfscope}%
\end{pgfscope}%
\begin{pgfscope}%
\pgfpathrectangle{\pgfqpoint{1.125000in}{0.750000in}}{\pgfqpoint{6.975000in}{4.530000in}} %
\pgfusepath{clip}%
\pgfsetroundcap%
\pgfsetroundjoin%
\pgfsetlinewidth{1.405250pt}%
\definecolor{currentstroke}{rgb}{0.318370,0.318370,0.318370}%
\pgfsetstrokecolor{currentstroke}%
\pgfsetdash{}{0pt}%
\pgfpathmoveto{\pgfqpoint{1.442045in}{0.955909in}}%
\pgfpathlineto{\pgfqpoint{3.047526in}{1.855763in}}%
\pgfpathlineto{\pgfqpoint{4.135329in}{1.940423in}}%
\pgfpathlineto{\pgfqpoint{4.959112in}{1.981055in}}%
\pgfpathlineto{\pgfqpoint{5.622327in}{2.000131in}}%
\pgfpathlineto{\pgfqpoint{6.177474in}{2.027659in}}%
\pgfpathlineto{\pgfqpoint{6.654885in}{2.066823in}}%
\pgfpathlineto{\pgfqpoint{7.073681in}{2.081521in}}%
\pgfpathlineto{\pgfqpoint{7.446694in}{2.096743in}}%
\pgfpathlineto{\pgfqpoint{7.782955in}{2.120906in}}%
\pgfusepath{stroke}%
\end{pgfscope}%
\begin{pgfscope}%
\pgfpathrectangle{\pgfqpoint{1.125000in}{0.750000in}}{\pgfqpoint{6.975000in}{4.530000in}} %
\pgfusepath{clip}%
\pgfsetbuttcap%
\pgfsetroundjoin%
\definecolor{currentfill}{rgb}{0.318370,0.318370,0.318370}%
\pgfsetfillcolor{currentfill}%
\pgfsetlinewidth{0.000000pt}%
\definecolor{currentstroke}{rgb}{0.318370,0.318370,0.318370}%
\pgfsetstrokecolor{currentstroke}%
\pgfsetdash{}{0pt}%
\pgfsys@defobject{currentmarker}{\pgfqpoint{-0.048611in}{-0.048611in}}{\pgfqpoint{0.048611in}{0.048611in}}{%
\pgfpathmoveto{\pgfqpoint{0.000000in}{-0.048611in}}%
\pgfpathcurveto{\pgfqpoint{0.012892in}{-0.048611in}}{\pgfqpoint{0.025257in}{-0.043489in}}{\pgfqpoint{0.034373in}{-0.034373in}}%
\pgfpathcurveto{\pgfqpoint{0.043489in}{-0.025257in}}{\pgfqpoint{0.048611in}{-0.012892in}}{\pgfqpoint{0.048611in}{0.000000in}}%
\pgfpathcurveto{\pgfqpoint{0.048611in}{0.012892in}}{\pgfqpoint{0.043489in}{0.025257in}}{\pgfqpoint{0.034373in}{0.034373in}}%
\pgfpathcurveto{\pgfqpoint{0.025257in}{0.043489in}}{\pgfqpoint{0.012892in}{0.048611in}}{\pgfqpoint{0.000000in}{0.048611in}}%
\pgfpathcurveto{\pgfqpoint{-0.012892in}{0.048611in}}{\pgfqpoint{-0.025257in}{0.043489in}}{\pgfqpoint{-0.034373in}{0.034373in}}%
\pgfpathcurveto{\pgfqpoint{-0.043489in}{0.025257in}}{\pgfqpoint{-0.048611in}{0.012892in}}{\pgfqpoint{-0.048611in}{0.000000in}}%
\pgfpathcurveto{\pgfqpoint{-0.048611in}{-0.012892in}}{\pgfqpoint{-0.043489in}{-0.025257in}}{\pgfqpoint{-0.034373in}{-0.034373in}}%
\pgfpathcurveto{\pgfqpoint{-0.025257in}{-0.043489in}}{\pgfqpoint{-0.012892in}{-0.048611in}}{\pgfqpoint{0.000000in}{-0.048611in}}%
\pgfpathclose%
\pgfusepath{fill}%
}%
\begin{pgfscope}%
\pgfsys@transformshift{1.442045in}{0.955909in}%
\pgfsys@useobject{currentmarker}{}%
\end{pgfscope}%
\begin{pgfscope}%
\pgfsys@transformshift{3.047526in}{1.855763in}%
\pgfsys@useobject{currentmarker}{}%
\end{pgfscope}%
\begin{pgfscope}%
\pgfsys@transformshift{4.135329in}{1.940423in}%
\pgfsys@useobject{currentmarker}{}%
\end{pgfscope}%
\begin{pgfscope}%
\pgfsys@transformshift{4.959112in}{1.981055in}%
\pgfsys@useobject{currentmarker}{}%
\end{pgfscope}%
\begin{pgfscope}%
\pgfsys@transformshift{5.622327in}{2.000131in}%
\pgfsys@useobject{currentmarker}{}%
\end{pgfscope}%
\begin{pgfscope}%
\pgfsys@transformshift{6.177474in}{2.027659in}%
\pgfsys@useobject{currentmarker}{}%
\end{pgfscope}%
\begin{pgfscope}%
\pgfsys@transformshift{6.654885in}{2.066823in}%
\pgfsys@useobject{currentmarker}{}%
\end{pgfscope}%
\begin{pgfscope}%
\pgfsys@transformshift{7.073681in}{2.081521in}%
\pgfsys@useobject{currentmarker}{}%
\end{pgfscope}%
\begin{pgfscope}%
\pgfsys@transformshift{7.446694in}{2.096743in}%
\pgfsys@useobject{currentmarker}{}%
\end{pgfscope}%
\begin{pgfscope}%
\pgfsys@transformshift{7.782955in}{2.120906in}%
\pgfsys@useobject{currentmarker}{}%
\end{pgfscope}%
\end{pgfscope}%
\begin{pgfscope}%
\pgfpathrectangle{\pgfqpoint{1.125000in}{0.750000in}}{\pgfqpoint{6.975000in}{4.530000in}} %
\pgfusepath{clip}%
\pgfsetroundcap%
\pgfsetroundjoin%
\pgfsetlinewidth{1.405250pt}%
\definecolor{currentstroke}{rgb}{0.501961,0.501961,0.501961}%
\pgfsetstrokecolor{currentstroke}%
\pgfsetdash{}{0pt}%
\pgfpathmoveto{\pgfqpoint{1.442045in}{3.791013in}}%
\pgfpathlineto{\pgfqpoint{3.047526in}{3.484600in}}%
\pgfpathlineto{\pgfqpoint{4.135329in}{3.496593in}}%
\pgfpathlineto{\pgfqpoint{4.959112in}{3.493550in}}%
\pgfpathlineto{\pgfqpoint{5.622327in}{3.492872in}}%
\pgfpathlineto{\pgfqpoint{6.177474in}{3.498084in}}%
\pgfpathlineto{\pgfqpoint{6.654885in}{3.502575in}}%
\pgfpathlineto{\pgfqpoint{7.073681in}{3.499038in}}%
\pgfpathlineto{\pgfqpoint{7.446694in}{3.503625in}}%
\pgfpathlineto{\pgfqpoint{7.782955in}{3.503126in}}%
\pgfusepath{stroke}%
\end{pgfscope}%
\begin{pgfscope}%
\pgfpathrectangle{\pgfqpoint{1.125000in}{0.750000in}}{\pgfqpoint{6.975000in}{4.530000in}} %
\pgfusepath{clip}%
\pgfsetbuttcap%
\pgfsetroundjoin%
\definecolor{currentfill}{rgb}{0.501961,0.501961,0.501961}%
\pgfsetfillcolor{currentfill}%
\pgfsetlinewidth{0.000000pt}%
\definecolor{currentstroke}{rgb}{0.501961,0.501961,0.501961}%
\pgfsetstrokecolor{currentstroke}%
\pgfsetdash{}{0pt}%
\pgfsys@defobject{currentmarker}{\pgfqpoint{-0.034722in}{-0.034722in}}{\pgfqpoint{0.034722in}{0.034722in}}{%
\pgfpathmoveto{\pgfqpoint{0.000000in}{-0.034722in}}%
\pgfpathcurveto{\pgfqpoint{0.009208in}{-0.034722in}}{\pgfqpoint{0.018041in}{-0.031064in}}{\pgfqpoint{0.024552in}{-0.024552in}}%
\pgfpathcurveto{\pgfqpoint{0.031064in}{-0.018041in}}{\pgfqpoint{0.034722in}{-0.009208in}}{\pgfqpoint{0.034722in}{0.000000in}}%
\pgfpathcurveto{\pgfqpoint{0.034722in}{0.009208in}}{\pgfqpoint{0.031064in}{0.018041in}}{\pgfqpoint{0.024552in}{0.024552in}}%
\pgfpathcurveto{\pgfqpoint{0.018041in}{0.031064in}}{\pgfqpoint{0.009208in}{0.034722in}}{\pgfqpoint{0.000000in}{0.034722in}}%
\pgfpathcurveto{\pgfqpoint{-0.009208in}{0.034722in}}{\pgfqpoint{-0.018041in}{0.031064in}}{\pgfqpoint{-0.024552in}{0.024552in}}%
\pgfpathcurveto{\pgfqpoint{-0.031064in}{0.018041in}}{\pgfqpoint{-0.034722in}{0.009208in}}{\pgfqpoint{-0.034722in}{0.000000in}}%
\pgfpathcurveto{\pgfqpoint{-0.034722in}{-0.009208in}}{\pgfqpoint{-0.031064in}{-0.018041in}}{\pgfqpoint{-0.024552in}{-0.024552in}}%
\pgfpathcurveto{\pgfqpoint{-0.018041in}{-0.031064in}}{\pgfqpoint{-0.009208in}{-0.034722in}}{\pgfqpoint{0.000000in}{-0.034722in}}%
\pgfpathclose%
\pgfusepath{fill}%
}%
\begin{pgfscope}%
\pgfsys@transformshift{1.442045in}{3.791013in}%
\pgfsys@useobject{currentmarker}{}%
\end{pgfscope}%
\begin{pgfscope}%
\pgfsys@transformshift{3.047526in}{3.484600in}%
\pgfsys@useobject{currentmarker}{}%
\end{pgfscope}%
\begin{pgfscope}%
\pgfsys@transformshift{4.135329in}{3.496593in}%
\pgfsys@useobject{currentmarker}{}%
\end{pgfscope}%
\begin{pgfscope}%
\pgfsys@transformshift{4.959112in}{3.493550in}%
\pgfsys@useobject{currentmarker}{}%
\end{pgfscope}%
\begin{pgfscope}%
\pgfsys@transformshift{5.622327in}{3.492872in}%
\pgfsys@useobject{currentmarker}{}%
\end{pgfscope}%
\begin{pgfscope}%
\pgfsys@transformshift{6.177474in}{3.498084in}%
\pgfsys@useobject{currentmarker}{}%
\end{pgfscope}%
\begin{pgfscope}%
\pgfsys@transformshift{6.654885in}{3.502575in}%
\pgfsys@useobject{currentmarker}{}%
\end{pgfscope}%
\begin{pgfscope}%
\pgfsys@transformshift{7.073681in}{3.499038in}%
\pgfsys@useobject{currentmarker}{}%
\end{pgfscope}%
\begin{pgfscope}%
\pgfsys@transformshift{7.446694in}{3.503625in}%
\pgfsys@useobject{currentmarker}{}%
\end{pgfscope}%
\begin{pgfscope}%
\pgfsys@transformshift{7.782955in}{3.503126in}%
\pgfsys@useobject{currentmarker}{}%
\end{pgfscope}%
\end{pgfscope}%
\begin{pgfscope}%
\pgfpathrectangle{\pgfqpoint{1.125000in}{0.750000in}}{\pgfqpoint{6.975000in}{4.530000in}} %
\pgfusepath{clip}%
\pgfsetbuttcap%
\pgfsetroundjoin%
\pgfsetlinewidth{1.405250pt}%
\definecolor{currentstroke}{rgb}{0.000000,0.000000,0.000000}%
\pgfsetstrokecolor{currentstroke}%
\pgfsetdash{{5.180000pt}{2.240000pt}}{0.000000pt}%
\pgfpathmoveto{\pgfqpoint{1.442045in}{5.074091in}}%
\pgfpathlineto{\pgfqpoint{3.047526in}{5.074091in}}%
\pgfpathlineto{\pgfqpoint{4.135329in}{5.074091in}}%
\pgfpathlineto{\pgfqpoint{4.959112in}{5.074091in}}%
\pgfpathlineto{\pgfqpoint{5.622327in}{5.074091in}}%
\pgfpathlineto{\pgfqpoint{6.177474in}{5.074091in}}%
\pgfpathlineto{\pgfqpoint{6.654885in}{5.074091in}}%
\pgfpathlineto{\pgfqpoint{7.073681in}{5.074091in}}%
\pgfpathlineto{\pgfqpoint{7.446694in}{5.074091in}}%
\pgfpathlineto{\pgfqpoint{7.782955in}{5.074091in}}%
\pgfusepath{stroke}%
\end{pgfscope}%
\begin{pgfscope}%
\pgfsetrectcap%
\pgfsetmiterjoin%
\pgfsetlinewidth{1.254687pt}%
\definecolor{currentstroke}{rgb}{0.150000,0.150000,0.150000}%
\pgfsetstrokecolor{currentstroke}%
\pgfsetdash{}{0pt}%
\pgfpathmoveto{\pgfqpoint{1.125000in}{0.750000in}}%
\pgfpathlineto{\pgfqpoint{1.125000in}{5.280000in}}%
\pgfusepath{stroke}%
\end{pgfscope}%
\begin{pgfscope}%
\pgfsetrectcap%
\pgfsetmiterjoin%
\pgfsetlinewidth{1.254687pt}%
\definecolor{currentstroke}{rgb}{0.150000,0.150000,0.150000}%
\pgfsetstrokecolor{currentstroke}%
\pgfsetdash{}{0pt}%
\pgfpathmoveto{\pgfqpoint{8.100000in}{0.750000in}}%
\pgfpathlineto{\pgfqpoint{8.100000in}{5.280000in}}%
\pgfusepath{stroke}%
\end{pgfscope}%
\begin{pgfscope}%
\pgfsetrectcap%
\pgfsetmiterjoin%
\pgfsetlinewidth{1.254687pt}%
\definecolor{currentstroke}{rgb}{0.150000,0.150000,0.150000}%
\pgfsetstrokecolor{currentstroke}%
\pgfsetdash{}{0pt}%
\pgfpathmoveto{\pgfqpoint{1.125000in}{0.750000in}}%
\pgfpathlineto{\pgfqpoint{8.100000in}{0.750000in}}%
\pgfusepath{stroke}%
\end{pgfscope}%
\begin{pgfscope}%
\pgfsetrectcap%
\pgfsetmiterjoin%
\pgfsetlinewidth{1.254687pt}%
\definecolor{currentstroke}{rgb}{0.150000,0.150000,0.150000}%
\pgfsetstrokecolor{currentstroke}%
\pgfsetdash{}{0pt}%
\pgfpathmoveto{\pgfqpoint{1.125000in}{5.280000in}}%
\pgfpathlineto{\pgfqpoint{8.100000in}{5.280000in}}%
\pgfusepath{stroke}%
\end{pgfscope}%
\begin{pgfscope}%
\pgfsetroundcap%
\pgfsetroundjoin%
\pgfsetlinewidth{1.405250pt}%
\definecolor{currentstroke}{rgb}{0.133333,0.133333,0.133333}%
\pgfsetstrokecolor{currentstroke}%
\pgfsetdash{}{0pt}%
\pgfpathmoveto{\pgfqpoint{5.448533in}{4.780632in}}%
\pgfpathlineto{\pgfqpoint{5.865200in}{4.780632in}}%
\pgfusepath{stroke}%
\end{pgfscope}%
\begin{pgfscope}%
\pgfsetbuttcap%
\pgfsetmiterjoin%
\definecolor{currentfill}{rgb}{0.133333,0.133333,0.133333}%
\pgfsetfillcolor{currentfill}%
\pgfsetlinewidth{0.000000pt}%
\definecolor{currentstroke}{rgb}{0.133333,0.133333,0.133333}%
\pgfsetstrokecolor{currentstroke}%
\pgfsetdash{}{0pt}%
\pgfsys@defobject{currentmarker}{\pgfqpoint{-0.062500in}{-0.062500in}}{\pgfqpoint{0.062500in}{0.062500in}}{%
\pgfpathmoveto{\pgfqpoint{-0.000000in}{-0.062500in}}%
\pgfpathlineto{\pgfqpoint{0.062500in}{0.062500in}}%
\pgfpathlineto{\pgfqpoint{-0.062500in}{0.062500in}}%
\pgfpathclose%
\pgfusepath{fill}%
}%
\begin{pgfscope}%
\pgfsys@transformshift{5.656867in}{4.780632in}%
\pgfsys@useobject{currentmarker}{}%
\end{pgfscope}%
\end{pgfscope}%
\begin{pgfscope}%
\definecolor{textcolor}{rgb}{0.150000,0.150000,0.150000}%
\pgfsetstrokecolor{textcolor}%
\pgfsetfillcolor{textcolor}%
\pgftext[x=6.031867in,y=4.707716in,left,base]{\color{textcolor}\sffamily\fontsize{15.000000}{18.000000}\selectfont Linear}%
\end{pgfscope}%
\begin{pgfscope}%
\pgfsetroundcap%
\pgfsetroundjoin%
\pgfsetlinewidth{1.405250pt}%
\definecolor{currentstroke}{rgb}{0.318370,0.318370,0.318370}%
\pgfsetstrokecolor{currentstroke}%
\pgfsetdash{}{0pt}%
\pgfpathmoveto{\pgfqpoint{5.448533in}{4.474846in}}%
\pgfpathlineto{\pgfqpoint{5.865200in}{4.474846in}}%
\pgfusepath{stroke}%
\end{pgfscope}%
\begin{pgfscope}%
\pgfsetbuttcap%
\pgfsetroundjoin%
\definecolor{currentfill}{rgb}{0.318370,0.318370,0.318370}%
\pgfsetfillcolor{currentfill}%
\pgfsetlinewidth{0.000000pt}%
\definecolor{currentstroke}{rgb}{0.318370,0.318370,0.318370}%
\pgfsetstrokecolor{currentstroke}%
\pgfsetdash{}{0pt}%
\pgfsys@defobject{currentmarker}{\pgfqpoint{-0.048611in}{-0.048611in}}{\pgfqpoint{0.048611in}{0.048611in}}{%
\pgfpathmoveto{\pgfqpoint{0.000000in}{-0.048611in}}%
\pgfpathcurveto{\pgfqpoint{0.012892in}{-0.048611in}}{\pgfqpoint{0.025257in}{-0.043489in}}{\pgfqpoint{0.034373in}{-0.034373in}}%
\pgfpathcurveto{\pgfqpoint{0.043489in}{-0.025257in}}{\pgfqpoint{0.048611in}{-0.012892in}}{\pgfqpoint{0.048611in}{0.000000in}}%
\pgfpathcurveto{\pgfqpoint{0.048611in}{0.012892in}}{\pgfqpoint{0.043489in}{0.025257in}}{\pgfqpoint{0.034373in}{0.034373in}}%
\pgfpathcurveto{\pgfqpoint{0.025257in}{0.043489in}}{\pgfqpoint{0.012892in}{0.048611in}}{\pgfqpoint{0.000000in}{0.048611in}}%
\pgfpathcurveto{\pgfqpoint{-0.012892in}{0.048611in}}{\pgfqpoint{-0.025257in}{0.043489in}}{\pgfqpoint{-0.034373in}{0.034373in}}%
\pgfpathcurveto{\pgfqpoint{-0.043489in}{0.025257in}}{\pgfqpoint{-0.048611in}{0.012892in}}{\pgfqpoint{-0.048611in}{0.000000in}}%
\pgfpathcurveto{\pgfqpoint{-0.048611in}{-0.012892in}}{\pgfqpoint{-0.043489in}{-0.025257in}}{\pgfqpoint{-0.034373in}{-0.034373in}}%
\pgfpathcurveto{\pgfqpoint{-0.025257in}{-0.043489in}}{\pgfqpoint{-0.012892in}{-0.048611in}}{\pgfqpoint{0.000000in}{-0.048611in}}%
\pgfpathclose%
\pgfusepath{fill}%
}%
\begin{pgfscope}%
\pgfsys@transformshift{5.656867in}{4.474846in}%
\pgfsys@useobject{currentmarker}{}%
\end{pgfscope}%
\end{pgfscope}%
\begin{pgfscope}%
\definecolor{textcolor}{rgb}{0.150000,0.150000,0.150000}%
\pgfsetstrokecolor{textcolor}%
\pgfsetfillcolor{textcolor}%
\pgftext[x=6.031867in,y=4.401930in,left,base]{\color{textcolor}\sffamily\fontsize{15.000000}{18.000000}\selectfont DMG}%
\end{pgfscope}%
\begin{pgfscope}%
\pgfsetroundcap%
\pgfsetroundjoin%
\pgfsetlinewidth{1.405250pt}%
\definecolor{currentstroke}{rgb}{0.501961,0.501961,0.501961}%
\pgfsetstrokecolor{currentstroke}%
\pgfsetdash{}{0pt}%
\pgfpathmoveto{\pgfqpoint{5.448533in}{4.169061in}}%
\pgfpathlineto{\pgfqpoint{5.865200in}{4.169061in}}%
\pgfusepath{stroke}%
\end{pgfscope}%
\begin{pgfscope}%
\pgfsetbuttcap%
\pgfsetroundjoin%
\definecolor{currentfill}{rgb}{0.501961,0.501961,0.501961}%
\pgfsetfillcolor{currentfill}%
\pgfsetlinewidth{0.000000pt}%
\definecolor{currentstroke}{rgb}{0.501961,0.501961,0.501961}%
\pgfsetstrokecolor{currentstroke}%
\pgfsetdash{}{0pt}%
\pgfsys@defobject{currentmarker}{\pgfqpoint{-0.034722in}{-0.034722in}}{\pgfqpoint{0.034722in}{0.034722in}}{%
\pgfpathmoveto{\pgfqpoint{0.000000in}{-0.034722in}}%
\pgfpathcurveto{\pgfqpoint{0.009208in}{-0.034722in}}{\pgfqpoint{0.018041in}{-0.031064in}}{\pgfqpoint{0.024552in}{-0.024552in}}%
\pgfpathcurveto{\pgfqpoint{0.031064in}{-0.018041in}}{\pgfqpoint{0.034722in}{-0.009208in}}{\pgfqpoint{0.034722in}{0.000000in}}%
\pgfpathcurveto{\pgfqpoint{0.034722in}{0.009208in}}{\pgfqpoint{0.031064in}{0.018041in}}{\pgfqpoint{0.024552in}{0.024552in}}%
\pgfpathcurveto{\pgfqpoint{0.018041in}{0.031064in}}{\pgfqpoint{0.009208in}{0.034722in}}{\pgfqpoint{0.000000in}{0.034722in}}%
\pgfpathcurveto{\pgfqpoint{-0.009208in}{0.034722in}}{\pgfqpoint{-0.018041in}{0.031064in}}{\pgfqpoint{-0.024552in}{0.024552in}}%
\pgfpathcurveto{\pgfqpoint{-0.031064in}{0.018041in}}{\pgfqpoint{-0.034722in}{0.009208in}}{\pgfqpoint{-0.034722in}{0.000000in}}%
\pgfpathcurveto{\pgfqpoint{-0.034722in}{-0.009208in}}{\pgfqpoint{-0.031064in}{-0.018041in}}{\pgfqpoint{-0.024552in}{-0.024552in}}%
\pgfpathcurveto{\pgfqpoint{-0.018041in}{-0.031064in}}{\pgfqpoint{-0.009208in}{-0.034722in}}{\pgfqpoint{0.000000in}{-0.034722in}}%
\pgfpathclose%
\pgfusepath{fill}%
}%
\begin{pgfscope}%
\pgfsys@transformshift{5.656867in}{4.169061in}%
\pgfsys@useobject{currentmarker}{}%
\end{pgfscope}%
\end{pgfscope}%
\begin{pgfscope}%
\definecolor{textcolor}{rgb}{0.150000,0.150000,0.150000}%
\pgfsetstrokecolor{textcolor}%
\pgfsetfillcolor{textcolor}%
\pgftext[x=6.031867in,y=4.096144in,left,base]{\color{textcolor}\sffamily\fontsize{15.000000}{18.000000}\selectfont AMG}%
\end{pgfscope}%
\begin{pgfscope}%
\pgfsetbuttcap%
\pgfsetroundjoin%
\pgfsetlinewidth{1.405250pt}%
\definecolor{currentstroke}{rgb}{0.000000,0.000000,0.000000}%
\pgfsetstrokecolor{currentstroke}%
\pgfsetdash{{5.180000pt}{2.240000pt}}{0.000000pt}%
\pgfpathmoveto{\pgfqpoint{5.448533in}{3.863275in}}%
\pgfpathlineto{\pgfqpoint{5.865200in}{3.863275in}}%
\pgfusepath{stroke}%
\end{pgfscope}%
\begin{pgfscope}%
\definecolor{textcolor}{rgb}{0.150000,0.150000,0.150000}%
\pgfsetstrokecolor{textcolor}%
\pgfsetfillcolor{textcolor}%
\pgftext[x=6.031867in,y=3.790358in,left,base]{\color{textcolor}\sffamily\fontsize{15.000000}{18.000000}\selectfont Boundary \(\displaystyle \rho = 1\)}%
\end{pgfscope}%
\end{pgfpicture}%
\makeatother%
\endgroup%

%% file: diff_m3points_epsrange_n_127_bw.pgf
\begingroup%
\makeatletter%
\begin{pgfpicture}%
\pgfpathrectangle{\pgfpointorigin}{\pgfqpoint{9.000000in}{6.000000in}}%
\pgfusepath{use as bounding box, clip}%
\begin{pgfscope}%
\pgfsetbuttcap%
\pgfsetmiterjoin%
\definecolor{currentfill}{rgb}{1.000000,1.000000,1.000000}%
\pgfsetfillcolor{currentfill}%
\pgfsetlinewidth{0.000000pt}%
\definecolor{currentstroke}{rgb}{1.000000,1.000000,1.000000}%
\pgfsetstrokecolor{currentstroke}%
\pgfsetdash{}{0pt}%
\pgfpathmoveto{\pgfqpoint{0.000000in}{0.000000in}}%
\pgfpathlineto{\pgfqpoint{9.000000in}{0.000000in}}%
\pgfpathlineto{\pgfqpoint{9.000000in}{6.000000in}}%
\pgfpathlineto{\pgfqpoint{0.000000in}{6.000000in}}%
\pgfpathclose%
\pgfusepath{fill}%
\end{pgfscope}%
\begin{pgfscope}%
\pgfsetbuttcap%
\pgfsetmiterjoin%
\definecolor{currentfill}{rgb}{1.000000,1.000000,1.000000}%
\pgfsetfillcolor{currentfill}%
\pgfsetlinewidth{0.000000pt}%
\definecolor{currentstroke}{rgb}{0.000000,0.000000,0.000000}%
\pgfsetstrokecolor{currentstroke}%
\pgfsetstrokeopacity{0.000000}%
\pgfsetdash{}{0pt}%
\pgfpathmoveto{\pgfqpoint{1.125000in}{0.750000in}}%
\pgfpathlineto{\pgfqpoint{8.100000in}{0.750000in}}%
\pgfpathlineto{\pgfqpoint{8.100000in}{5.280000in}}%
\pgfpathlineto{\pgfqpoint{1.125000in}{5.280000in}}%
\pgfpathclose%
\pgfusepath{fill}%
\end{pgfscope}%
\begin{pgfscope}%
\pgfsetbuttcap%
\pgfsetroundjoin%
\definecolor{currentfill}{rgb}{0.150000,0.150000,0.150000}%
\pgfsetfillcolor{currentfill}%
\pgfsetlinewidth{0.803000pt}%
\definecolor{currentstroke}{rgb}{0.150000,0.150000,0.150000}%
\pgfsetstrokecolor{currentstroke}%
\pgfsetdash{}{0pt}%
\pgfsys@defobject{currentmarker}{\pgfqpoint{0.000000in}{-0.083333in}}{\pgfqpoint{0.000000in}{0.000000in}}{%
\pgfpathmoveto{\pgfqpoint{0.000000in}{0.000000in}}%
\pgfpathlineto{\pgfqpoint{0.000000in}{-0.083333in}}%
\pgfusepath{stroke,fill}%
}%
\begin{pgfscope}%
\pgfsys@transformshift{2.056029in}{0.750000in}%
\pgfsys@useobject{currentmarker}{}%
\end{pgfscope}%
\end{pgfscope}%
\begin{pgfscope}%
\definecolor{textcolor}{rgb}{0.150000,0.150000,0.150000}%
\pgfsetstrokecolor{textcolor}%
\pgfsetfillcolor{textcolor}%
\pgftext[x=2.056029in,y=0.588889in,,top]{\color{textcolor}\sffamily\fontsize{16.000000}{19.200000}\selectfont \(\displaystyle 10^{-2}\)}%
\end{pgfscope}%
\begin{pgfscope}%
\pgfsetbuttcap%
\pgfsetroundjoin%
\definecolor{currentfill}{rgb}{0.150000,0.150000,0.150000}%
\pgfsetfillcolor{currentfill}%
\pgfsetlinewidth{0.803000pt}%
\definecolor{currentstroke}{rgb}{0.150000,0.150000,0.150000}%
\pgfsetstrokecolor{currentstroke}%
\pgfsetdash{}{0pt}%
\pgfsys@defobject{currentmarker}{\pgfqpoint{0.000000in}{-0.083333in}}{\pgfqpoint{0.000000in}{0.000000in}}{%
\pgfpathmoveto{\pgfqpoint{0.000000in}{0.000000in}}%
\pgfpathlineto{\pgfqpoint{0.000000in}{-0.083333in}}%
\pgfusepath{stroke,fill}%
}%
\begin{pgfscope}%
\pgfsys@transformshift{7.782955in}{0.750000in}%
\pgfsys@useobject{currentmarker}{}%
\end{pgfscope}%
\end{pgfscope}%
\begin{pgfscope}%
\definecolor{textcolor}{rgb}{0.150000,0.150000,0.150000}%
\pgfsetstrokecolor{textcolor}%
\pgfsetfillcolor{textcolor}%
\pgftext[x=7.782955in,y=0.588889in,,top]{\color{textcolor}\sffamily\fontsize{16.000000}{19.200000}\selectfont \(\displaystyle 10^{-1}\)}%
\end{pgfscope}%
\begin{pgfscope}%
\pgfsetbuttcap%
\pgfsetroundjoin%
\definecolor{currentfill}{rgb}{0.150000,0.150000,0.150000}%
\pgfsetfillcolor{currentfill}%
\pgfsetlinewidth{0.401500pt}%
\definecolor{currentstroke}{rgb}{0.150000,0.150000,0.150000}%
\pgfsetstrokecolor{currentstroke}%
\pgfsetdash{}{0pt}%
\pgfsys@defobject{currentmarker}{\pgfqpoint{0.000000in}{-0.041667in}}{\pgfqpoint{0.000000in}{0.000000in}}{%
\pgfpathmoveto{\pgfqpoint{0.000000in}{0.000000in}}%
\pgfpathlineto{\pgfqpoint{0.000000in}{-0.041667in}}%
\pgfusepath{stroke,fill}%
}%
\begin{pgfscope}%
\pgfsys@transformshift{1.168917in}{0.750000in}%
\pgfsys@useobject{currentmarker}{}%
\end{pgfscope}%
\end{pgfscope}%
\begin{pgfscope}%
\pgfsetbuttcap%
\pgfsetroundjoin%
\definecolor{currentfill}{rgb}{0.150000,0.150000,0.150000}%
\pgfsetfillcolor{currentfill}%
\pgfsetlinewidth{0.401500pt}%
\definecolor{currentstroke}{rgb}{0.150000,0.150000,0.150000}%
\pgfsetstrokecolor{currentstroke}%
\pgfsetdash{}{0pt}%
\pgfsys@defobject{currentmarker}{\pgfqpoint{0.000000in}{-0.041667in}}{\pgfqpoint{0.000000in}{0.000000in}}{%
\pgfpathmoveto{\pgfqpoint{0.000000in}{0.000000in}}%
\pgfpathlineto{\pgfqpoint{0.000000in}{-0.041667in}}%
\pgfusepath{stroke,fill}%
}%
\begin{pgfscope}%
\pgfsys@transformshift{1.501033in}{0.750000in}%
\pgfsys@useobject{currentmarker}{}%
\end{pgfscope}%
\end{pgfscope}%
\begin{pgfscope}%
\pgfsetbuttcap%
\pgfsetroundjoin%
\definecolor{currentfill}{rgb}{0.150000,0.150000,0.150000}%
\pgfsetfillcolor{currentfill}%
\pgfsetlinewidth{0.401500pt}%
\definecolor{currentstroke}{rgb}{0.150000,0.150000,0.150000}%
\pgfsetstrokecolor{currentstroke}%
\pgfsetdash{}{0pt}%
\pgfsys@defobject{currentmarker}{\pgfqpoint{0.000000in}{-0.041667in}}{\pgfqpoint{0.000000in}{0.000000in}}{%
\pgfpathmoveto{\pgfqpoint{0.000000in}{0.000000in}}%
\pgfpathlineto{\pgfqpoint{0.000000in}{-0.041667in}}%
\pgfusepath{stroke,fill}%
}%
\begin{pgfscope}%
\pgfsys@transformshift{1.793979in}{0.750000in}%
\pgfsys@useobject{currentmarker}{}%
\end{pgfscope}%
\end{pgfscope}%
\begin{pgfscope}%
\pgfsetbuttcap%
\pgfsetroundjoin%
\definecolor{currentfill}{rgb}{0.150000,0.150000,0.150000}%
\pgfsetfillcolor{currentfill}%
\pgfsetlinewidth{0.401500pt}%
\definecolor{currentstroke}{rgb}{0.150000,0.150000,0.150000}%
\pgfsetstrokecolor{currentstroke}%
\pgfsetdash{}{0pt}%
\pgfsys@defobject{currentmarker}{\pgfqpoint{0.000000in}{-0.041667in}}{\pgfqpoint{0.000000in}{0.000000in}}{%
\pgfpathmoveto{\pgfqpoint{0.000000in}{0.000000in}}%
\pgfpathlineto{\pgfqpoint{0.000000in}{-0.041667in}}%
\pgfusepath{stroke,fill}%
}%
\begin{pgfscope}%
\pgfsys@transformshift{3.780005in}{0.750000in}%
\pgfsys@useobject{currentmarker}{}%
\end{pgfscope}%
\end{pgfscope}%
\begin{pgfscope}%
\pgfsetbuttcap%
\pgfsetroundjoin%
\definecolor{currentfill}{rgb}{0.150000,0.150000,0.150000}%
\pgfsetfillcolor{currentfill}%
\pgfsetlinewidth{0.401500pt}%
\definecolor{currentstroke}{rgb}{0.150000,0.150000,0.150000}%
\pgfsetstrokecolor{currentstroke}%
\pgfsetdash{}{0pt}%
\pgfsys@defobject{currentmarker}{\pgfqpoint{0.000000in}{-0.041667in}}{\pgfqpoint{0.000000in}{0.000000in}}{%
\pgfpathmoveto{\pgfqpoint{0.000000in}{0.000000in}}%
\pgfpathlineto{\pgfqpoint{0.000000in}{-0.041667in}}%
\pgfusepath{stroke,fill}%
}%
\begin{pgfscope}%
\pgfsys@transformshift{4.788467in}{0.750000in}%
\pgfsys@useobject{currentmarker}{}%
\end{pgfscope}%
\end{pgfscope}%
\begin{pgfscope}%
\pgfsetbuttcap%
\pgfsetroundjoin%
\definecolor{currentfill}{rgb}{0.150000,0.150000,0.150000}%
\pgfsetfillcolor{currentfill}%
\pgfsetlinewidth{0.401500pt}%
\definecolor{currentstroke}{rgb}{0.150000,0.150000,0.150000}%
\pgfsetstrokecolor{currentstroke}%
\pgfsetdash{}{0pt}%
\pgfsys@defobject{currentmarker}{\pgfqpoint{0.000000in}{-0.041667in}}{\pgfqpoint{0.000000in}{0.000000in}}{%
\pgfpathmoveto{\pgfqpoint{0.000000in}{0.000000in}}%
\pgfpathlineto{\pgfqpoint{0.000000in}{-0.041667in}}%
\pgfusepath{stroke,fill}%
}%
\begin{pgfscope}%
\pgfsys@transformshift{5.503982in}{0.750000in}%
\pgfsys@useobject{currentmarker}{}%
\end{pgfscope}%
\end{pgfscope}%
\begin{pgfscope}%
\pgfsetbuttcap%
\pgfsetroundjoin%
\definecolor{currentfill}{rgb}{0.150000,0.150000,0.150000}%
\pgfsetfillcolor{currentfill}%
\pgfsetlinewidth{0.401500pt}%
\definecolor{currentstroke}{rgb}{0.150000,0.150000,0.150000}%
\pgfsetstrokecolor{currentstroke}%
\pgfsetdash{}{0pt}%
\pgfsys@defobject{currentmarker}{\pgfqpoint{0.000000in}{-0.041667in}}{\pgfqpoint{0.000000in}{0.000000in}}{%
\pgfpathmoveto{\pgfqpoint{0.000000in}{0.000000in}}%
\pgfpathlineto{\pgfqpoint{0.000000in}{-0.041667in}}%
\pgfusepath{stroke,fill}%
}%
\begin{pgfscope}%
\pgfsys@transformshift{6.058978in}{0.750000in}%
\pgfsys@useobject{currentmarker}{}%
\end{pgfscope}%
\end{pgfscope}%
\begin{pgfscope}%
\pgfsetbuttcap%
\pgfsetroundjoin%
\definecolor{currentfill}{rgb}{0.150000,0.150000,0.150000}%
\pgfsetfillcolor{currentfill}%
\pgfsetlinewidth{0.401500pt}%
\definecolor{currentstroke}{rgb}{0.150000,0.150000,0.150000}%
\pgfsetstrokecolor{currentstroke}%
\pgfsetdash{}{0pt}%
\pgfsys@defobject{currentmarker}{\pgfqpoint{0.000000in}{-0.041667in}}{\pgfqpoint{0.000000in}{0.000000in}}{%
\pgfpathmoveto{\pgfqpoint{0.000000in}{0.000000in}}%
\pgfpathlineto{\pgfqpoint{0.000000in}{-0.041667in}}%
\pgfusepath{stroke,fill}%
}%
\begin{pgfscope}%
\pgfsys@transformshift{6.512443in}{0.750000in}%
\pgfsys@useobject{currentmarker}{}%
\end{pgfscope}%
\end{pgfscope}%
\begin{pgfscope}%
\pgfsetbuttcap%
\pgfsetroundjoin%
\definecolor{currentfill}{rgb}{0.150000,0.150000,0.150000}%
\pgfsetfillcolor{currentfill}%
\pgfsetlinewidth{0.401500pt}%
\definecolor{currentstroke}{rgb}{0.150000,0.150000,0.150000}%
\pgfsetstrokecolor{currentstroke}%
\pgfsetdash{}{0pt}%
\pgfsys@defobject{currentmarker}{\pgfqpoint{0.000000in}{-0.041667in}}{\pgfqpoint{0.000000in}{0.000000in}}{%
\pgfpathmoveto{\pgfqpoint{0.000000in}{0.000000in}}%
\pgfpathlineto{\pgfqpoint{0.000000in}{-0.041667in}}%
\pgfusepath{stroke,fill}%
}%
\begin{pgfscope}%
\pgfsys@transformshift{6.895843in}{0.750000in}%
\pgfsys@useobject{currentmarker}{}%
\end{pgfscope}%
\end{pgfscope}%
\begin{pgfscope}%
\pgfsetbuttcap%
\pgfsetroundjoin%
\definecolor{currentfill}{rgb}{0.150000,0.150000,0.150000}%
\pgfsetfillcolor{currentfill}%
\pgfsetlinewidth{0.401500pt}%
\definecolor{currentstroke}{rgb}{0.150000,0.150000,0.150000}%
\pgfsetstrokecolor{currentstroke}%
\pgfsetdash{}{0pt}%
\pgfsys@defobject{currentmarker}{\pgfqpoint{0.000000in}{-0.041667in}}{\pgfqpoint{0.000000in}{0.000000in}}{%
\pgfpathmoveto{\pgfqpoint{0.000000in}{0.000000in}}%
\pgfpathlineto{\pgfqpoint{0.000000in}{-0.041667in}}%
\pgfusepath{stroke,fill}%
}%
\begin{pgfscope}%
\pgfsys@transformshift{7.227958in}{0.750000in}%
\pgfsys@useobject{currentmarker}{}%
\end{pgfscope}%
\end{pgfscope}%
\begin{pgfscope}%
\pgfsetbuttcap%
\pgfsetroundjoin%
\definecolor{currentfill}{rgb}{0.150000,0.150000,0.150000}%
\pgfsetfillcolor{currentfill}%
\pgfsetlinewidth{0.401500pt}%
\definecolor{currentstroke}{rgb}{0.150000,0.150000,0.150000}%
\pgfsetstrokecolor{currentstroke}%
\pgfsetdash{}{0pt}%
\pgfsys@defobject{currentmarker}{\pgfqpoint{0.000000in}{-0.041667in}}{\pgfqpoint{0.000000in}{0.000000in}}{%
\pgfpathmoveto{\pgfqpoint{0.000000in}{0.000000in}}%
\pgfpathlineto{\pgfqpoint{0.000000in}{-0.041667in}}%
\pgfusepath{stroke,fill}%
}%
\begin{pgfscope}%
\pgfsys@transformshift{7.520905in}{0.750000in}%
\pgfsys@useobject{currentmarker}{}%
\end{pgfscope}%
\end{pgfscope}%
\begin{pgfscope}%
\definecolor{textcolor}{rgb}{0.150000,0.150000,0.150000}%
\pgfsetstrokecolor{textcolor}%
\pgfsetfillcolor{textcolor}%
\pgftext[x=4.612500in,y=0.318273in,,top]{\color{textcolor}\sffamily\fontsize{24.000000}{28.800000}\selectfont \(\displaystyle \varepsilon\)}%
\end{pgfscope}%
\begin{pgfscope}%
\pgfsetbuttcap%
\pgfsetroundjoin%
\definecolor{currentfill}{rgb}{0.150000,0.150000,0.150000}%
\pgfsetfillcolor{currentfill}%
\pgfsetlinewidth{0.803000pt}%
\definecolor{currentstroke}{rgb}{0.150000,0.150000,0.150000}%
\pgfsetstrokecolor{currentstroke}%
\pgfsetdash{}{0pt}%
\pgfsys@defobject{currentmarker}{\pgfqpoint{-0.083333in}{0.000000in}}{\pgfqpoint{0.000000in}{0.000000in}}{%
\pgfpathmoveto{\pgfqpoint{0.000000in}{0.000000in}}%
\pgfpathlineto{\pgfqpoint{-0.083333in}{0.000000in}}%
\pgfusepath{stroke,fill}%
}%
\begin{pgfscope}%
\pgfsys@transformshift{1.125000in}{2.805238in}%
\pgfsys@useobject{currentmarker}{}%
\end{pgfscope}%
\end{pgfscope}%
\begin{pgfscope}%
\definecolor{textcolor}{rgb}{0.150000,0.150000,0.150000}%
\pgfsetstrokecolor{textcolor}%
\pgfsetfillcolor{textcolor}%
\pgftext[x=0.439258in,y=2.689163in,left,base]{\color{textcolor}\sffamily\fontsize{22.000000}{26.400000}\selectfont \(\displaystyle 10^{-1}\)}%
\end{pgfscope}%
\begin{pgfscope}%
\pgfsetbuttcap%
\pgfsetroundjoin%
\definecolor{currentfill}{rgb}{0.150000,0.150000,0.150000}%
\pgfsetfillcolor{currentfill}%
\pgfsetlinewidth{0.803000pt}%
\definecolor{currentstroke}{rgb}{0.150000,0.150000,0.150000}%
\pgfsetstrokecolor{currentstroke}%
\pgfsetdash{}{0pt}%
\pgfsys@defobject{currentmarker}{\pgfqpoint{-0.083333in}{0.000000in}}{\pgfqpoint{0.000000in}{0.000000in}}{%
\pgfpathmoveto{\pgfqpoint{0.000000in}{0.000000in}}%
\pgfpathlineto{\pgfqpoint{-0.083333in}{0.000000in}}%
\pgfusepath{stroke,fill}%
}%
\begin{pgfscope}%
\pgfsys@transformshift{1.125000in}{5.074091in}%
\pgfsys@useobject{currentmarker}{}%
\end{pgfscope}%
\end{pgfscope}%
\begin{pgfscope}%
\definecolor{textcolor}{rgb}{0.150000,0.150000,0.150000}%
\pgfsetstrokecolor{textcolor}%
\pgfsetfillcolor{textcolor}%
\pgftext[x=0.594814in,y=4.958016in,left,base]{\color{textcolor}\sffamily\fontsize{22.000000}{26.400000}\selectfont \(\displaystyle 10^{0}\)}%
\end{pgfscope}%
\begin{pgfscope}%
\pgfsetbuttcap%
\pgfsetroundjoin%
\definecolor{currentfill}{rgb}{0.150000,0.150000,0.150000}%
\pgfsetfillcolor{currentfill}%
\pgfsetlinewidth{0.401500pt}%
\definecolor{currentstroke}{rgb}{0.150000,0.150000,0.150000}%
\pgfsetstrokecolor{currentstroke}%
\pgfsetdash{}{0pt}%
\pgfsys@defobject{currentmarker}{\pgfqpoint{-0.041667in}{0.000000in}}{\pgfqpoint{0.000000in}{0.000000in}}{%
\pgfpathmoveto{\pgfqpoint{0.000000in}{0.000000in}}%
\pgfpathlineto{\pgfqpoint{-0.041667in}{0.000000in}}%
\pgfusepath{stroke,fill}%
}%
\begin{pgfscope}%
\pgfsys@transformshift{1.125000in}{1.219378in}%
\pgfsys@useobject{currentmarker}{}%
\end{pgfscope}%
\end{pgfscope}%
\begin{pgfscope}%
\pgfsetbuttcap%
\pgfsetroundjoin%
\definecolor{currentfill}{rgb}{0.150000,0.150000,0.150000}%
\pgfsetfillcolor{currentfill}%
\pgfsetlinewidth{0.401500pt}%
\definecolor{currentstroke}{rgb}{0.150000,0.150000,0.150000}%
\pgfsetstrokecolor{currentstroke}%
\pgfsetdash{}{0pt}%
\pgfsys@defobject{currentmarker}{\pgfqpoint{-0.041667in}{0.000000in}}{\pgfqpoint{0.000000in}{0.000000in}}{%
\pgfpathmoveto{\pgfqpoint{0.000000in}{0.000000in}}%
\pgfpathlineto{\pgfqpoint{-0.041667in}{0.000000in}}%
\pgfusepath{stroke,fill}%
}%
\begin{pgfscope}%
\pgfsys@transformshift{1.125000in}{1.618903in}%
\pgfsys@useobject{currentmarker}{}%
\end{pgfscope}%
\end{pgfscope}%
\begin{pgfscope}%
\pgfsetbuttcap%
\pgfsetroundjoin%
\definecolor{currentfill}{rgb}{0.150000,0.150000,0.150000}%
\pgfsetfillcolor{currentfill}%
\pgfsetlinewidth{0.401500pt}%
\definecolor{currentstroke}{rgb}{0.150000,0.150000,0.150000}%
\pgfsetstrokecolor{currentstroke}%
\pgfsetdash{}{0pt}%
\pgfsys@defobject{currentmarker}{\pgfqpoint{-0.041667in}{0.000000in}}{\pgfqpoint{0.000000in}{0.000000in}}{%
\pgfpathmoveto{\pgfqpoint{0.000000in}{0.000000in}}%
\pgfpathlineto{\pgfqpoint{-0.041667in}{0.000000in}}%
\pgfusepath{stroke,fill}%
}%
\begin{pgfscope}%
\pgfsys@transformshift{1.125000in}{1.902371in}%
\pgfsys@useobject{currentmarker}{}%
\end{pgfscope}%
\end{pgfscope}%
\begin{pgfscope}%
\pgfsetbuttcap%
\pgfsetroundjoin%
\definecolor{currentfill}{rgb}{0.150000,0.150000,0.150000}%
\pgfsetfillcolor{currentfill}%
\pgfsetlinewidth{0.401500pt}%
\definecolor{currentstroke}{rgb}{0.150000,0.150000,0.150000}%
\pgfsetstrokecolor{currentstroke}%
\pgfsetdash{}{0pt}%
\pgfsys@defobject{currentmarker}{\pgfqpoint{-0.041667in}{0.000000in}}{\pgfqpoint{0.000000in}{0.000000in}}{%
\pgfpathmoveto{\pgfqpoint{0.000000in}{0.000000in}}%
\pgfpathlineto{\pgfqpoint{-0.041667in}{0.000000in}}%
\pgfusepath{stroke,fill}%
}%
\begin{pgfscope}%
\pgfsys@transformshift{1.125000in}{2.122245in}%
\pgfsys@useobject{currentmarker}{}%
\end{pgfscope}%
\end{pgfscope}%
\begin{pgfscope}%
\pgfsetbuttcap%
\pgfsetroundjoin%
\definecolor{currentfill}{rgb}{0.150000,0.150000,0.150000}%
\pgfsetfillcolor{currentfill}%
\pgfsetlinewidth{0.401500pt}%
\definecolor{currentstroke}{rgb}{0.150000,0.150000,0.150000}%
\pgfsetstrokecolor{currentstroke}%
\pgfsetdash{}{0pt}%
\pgfsys@defobject{currentmarker}{\pgfqpoint{-0.041667in}{0.000000in}}{\pgfqpoint{0.000000in}{0.000000in}}{%
\pgfpathmoveto{\pgfqpoint{0.000000in}{0.000000in}}%
\pgfpathlineto{\pgfqpoint{-0.041667in}{0.000000in}}%
\pgfusepath{stroke,fill}%
}%
\begin{pgfscope}%
\pgfsys@transformshift{1.125000in}{2.301896in}%
\pgfsys@useobject{currentmarker}{}%
\end{pgfscope}%
\end{pgfscope}%
\begin{pgfscope}%
\pgfsetbuttcap%
\pgfsetroundjoin%
\definecolor{currentfill}{rgb}{0.150000,0.150000,0.150000}%
\pgfsetfillcolor{currentfill}%
\pgfsetlinewidth{0.401500pt}%
\definecolor{currentstroke}{rgb}{0.150000,0.150000,0.150000}%
\pgfsetstrokecolor{currentstroke}%
\pgfsetdash{}{0pt}%
\pgfsys@defobject{currentmarker}{\pgfqpoint{-0.041667in}{0.000000in}}{\pgfqpoint{0.000000in}{0.000000in}}{%
\pgfpathmoveto{\pgfqpoint{0.000000in}{0.000000in}}%
\pgfpathlineto{\pgfqpoint{-0.041667in}{0.000000in}}%
\pgfusepath{stroke,fill}%
}%
\begin{pgfscope}%
\pgfsys@transformshift{1.125000in}{2.453788in}%
\pgfsys@useobject{currentmarker}{}%
\end{pgfscope}%
\end{pgfscope}%
\begin{pgfscope}%
\pgfsetbuttcap%
\pgfsetroundjoin%
\definecolor{currentfill}{rgb}{0.150000,0.150000,0.150000}%
\pgfsetfillcolor{currentfill}%
\pgfsetlinewidth{0.401500pt}%
\definecolor{currentstroke}{rgb}{0.150000,0.150000,0.150000}%
\pgfsetstrokecolor{currentstroke}%
\pgfsetdash{}{0pt}%
\pgfsys@defobject{currentmarker}{\pgfqpoint{-0.041667in}{0.000000in}}{\pgfqpoint{0.000000in}{0.000000in}}{%
\pgfpathmoveto{\pgfqpoint{0.000000in}{0.000000in}}%
\pgfpathlineto{\pgfqpoint{-0.041667in}{0.000000in}}%
\pgfusepath{stroke,fill}%
}%
\begin{pgfscope}%
\pgfsys@transformshift{1.125000in}{2.585364in}%
\pgfsys@useobject{currentmarker}{}%
\end{pgfscope}%
\end{pgfscope}%
\begin{pgfscope}%
\pgfsetbuttcap%
\pgfsetroundjoin%
\definecolor{currentfill}{rgb}{0.150000,0.150000,0.150000}%
\pgfsetfillcolor{currentfill}%
\pgfsetlinewidth{0.401500pt}%
\definecolor{currentstroke}{rgb}{0.150000,0.150000,0.150000}%
\pgfsetstrokecolor{currentstroke}%
\pgfsetdash{}{0pt}%
\pgfsys@defobject{currentmarker}{\pgfqpoint{-0.041667in}{0.000000in}}{\pgfqpoint{0.000000in}{0.000000in}}{%
\pgfpathmoveto{\pgfqpoint{0.000000in}{0.000000in}}%
\pgfpathlineto{\pgfqpoint{-0.041667in}{0.000000in}}%
\pgfusepath{stroke,fill}%
}%
\begin{pgfscope}%
\pgfsys@transformshift{1.125000in}{2.701421in}%
\pgfsys@useobject{currentmarker}{}%
\end{pgfscope}%
\end{pgfscope}%
\begin{pgfscope}%
\pgfsetbuttcap%
\pgfsetroundjoin%
\definecolor{currentfill}{rgb}{0.150000,0.150000,0.150000}%
\pgfsetfillcolor{currentfill}%
\pgfsetlinewidth{0.401500pt}%
\definecolor{currentstroke}{rgb}{0.150000,0.150000,0.150000}%
\pgfsetstrokecolor{currentstroke}%
\pgfsetdash{}{0pt}%
\pgfsys@defobject{currentmarker}{\pgfqpoint{-0.041667in}{0.000000in}}{\pgfqpoint{0.000000in}{0.000000in}}{%
\pgfpathmoveto{\pgfqpoint{0.000000in}{0.000000in}}%
\pgfpathlineto{\pgfqpoint{-0.041667in}{0.000000in}}%
\pgfusepath{stroke,fill}%
}%
\begin{pgfscope}%
\pgfsys@transformshift{1.125000in}{3.488231in}%
\pgfsys@useobject{currentmarker}{}%
\end{pgfscope}%
\end{pgfscope}%
\begin{pgfscope}%
\pgfsetbuttcap%
\pgfsetroundjoin%
\definecolor{currentfill}{rgb}{0.150000,0.150000,0.150000}%
\pgfsetfillcolor{currentfill}%
\pgfsetlinewidth{0.401500pt}%
\definecolor{currentstroke}{rgb}{0.150000,0.150000,0.150000}%
\pgfsetstrokecolor{currentstroke}%
\pgfsetdash{}{0pt}%
\pgfsys@defobject{currentmarker}{\pgfqpoint{-0.041667in}{0.000000in}}{\pgfqpoint{0.000000in}{0.000000in}}{%
\pgfpathmoveto{\pgfqpoint{0.000000in}{0.000000in}}%
\pgfpathlineto{\pgfqpoint{-0.041667in}{0.000000in}}%
\pgfusepath{stroke,fill}%
}%
\begin{pgfscope}%
\pgfsys@transformshift{1.125000in}{3.887756in}%
\pgfsys@useobject{currentmarker}{}%
\end{pgfscope}%
\end{pgfscope}%
\begin{pgfscope}%
\pgfsetbuttcap%
\pgfsetroundjoin%
\definecolor{currentfill}{rgb}{0.150000,0.150000,0.150000}%
\pgfsetfillcolor{currentfill}%
\pgfsetlinewidth{0.401500pt}%
\definecolor{currentstroke}{rgb}{0.150000,0.150000,0.150000}%
\pgfsetstrokecolor{currentstroke}%
\pgfsetdash{}{0pt}%
\pgfsys@defobject{currentmarker}{\pgfqpoint{-0.041667in}{0.000000in}}{\pgfqpoint{0.000000in}{0.000000in}}{%
\pgfpathmoveto{\pgfqpoint{0.000000in}{0.000000in}}%
\pgfpathlineto{\pgfqpoint{-0.041667in}{0.000000in}}%
\pgfusepath{stroke,fill}%
}%
\begin{pgfscope}%
\pgfsys@transformshift{1.125000in}{4.171224in}%
\pgfsys@useobject{currentmarker}{}%
\end{pgfscope}%
\end{pgfscope}%
\begin{pgfscope}%
\pgfsetbuttcap%
\pgfsetroundjoin%
\definecolor{currentfill}{rgb}{0.150000,0.150000,0.150000}%
\pgfsetfillcolor{currentfill}%
\pgfsetlinewidth{0.401500pt}%
\definecolor{currentstroke}{rgb}{0.150000,0.150000,0.150000}%
\pgfsetstrokecolor{currentstroke}%
\pgfsetdash{}{0pt}%
\pgfsys@defobject{currentmarker}{\pgfqpoint{-0.041667in}{0.000000in}}{\pgfqpoint{0.000000in}{0.000000in}}{%
\pgfpathmoveto{\pgfqpoint{0.000000in}{0.000000in}}%
\pgfpathlineto{\pgfqpoint{-0.041667in}{0.000000in}}%
\pgfusepath{stroke,fill}%
}%
\begin{pgfscope}%
\pgfsys@transformshift{1.125000in}{4.391098in}%
\pgfsys@useobject{currentmarker}{}%
\end{pgfscope}%
\end{pgfscope}%
\begin{pgfscope}%
\pgfsetbuttcap%
\pgfsetroundjoin%
\definecolor{currentfill}{rgb}{0.150000,0.150000,0.150000}%
\pgfsetfillcolor{currentfill}%
\pgfsetlinewidth{0.401500pt}%
\definecolor{currentstroke}{rgb}{0.150000,0.150000,0.150000}%
\pgfsetstrokecolor{currentstroke}%
\pgfsetdash{}{0pt}%
\pgfsys@defobject{currentmarker}{\pgfqpoint{-0.041667in}{0.000000in}}{\pgfqpoint{0.000000in}{0.000000in}}{%
\pgfpathmoveto{\pgfqpoint{0.000000in}{0.000000in}}%
\pgfpathlineto{\pgfqpoint{-0.041667in}{0.000000in}}%
\pgfusepath{stroke,fill}%
}%
\begin{pgfscope}%
\pgfsys@transformshift{1.125000in}{4.570749in}%
\pgfsys@useobject{currentmarker}{}%
\end{pgfscope}%
\end{pgfscope}%
\begin{pgfscope}%
\pgfsetbuttcap%
\pgfsetroundjoin%
\definecolor{currentfill}{rgb}{0.150000,0.150000,0.150000}%
\pgfsetfillcolor{currentfill}%
\pgfsetlinewidth{0.401500pt}%
\definecolor{currentstroke}{rgb}{0.150000,0.150000,0.150000}%
\pgfsetstrokecolor{currentstroke}%
\pgfsetdash{}{0pt}%
\pgfsys@defobject{currentmarker}{\pgfqpoint{-0.041667in}{0.000000in}}{\pgfqpoint{0.000000in}{0.000000in}}{%
\pgfpathmoveto{\pgfqpoint{0.000000in}{0.000000in}}%
\pgfpathlineto{\pgfqpoint{-0.041667in}{0.000000in}}%
\pgfusepath{stroke,fill}%
}%
\begin{pgfscope}%
\pgfsys@transformshift{1.125000in}{4.722641in}%
\pgfsys@useobject{currentmarker}{}%
\end{pgfscope}%
\end{pgfscope}%
\begin{pgfscope}%
\pgfsetbuttcap%
\pgfsetroundjoin%
\definecolor{currentfill}{rgb}{0.150000,0.150000,0.150000}%
\pgfsetfillcolor{currentfill}%
\pgfsetlinewidth{0.401500pt}%
\definecolor{currentstroke}{rgb}{0.150000,0.150000,0.150000}%
\pgfsetstrokecolor{currentstroke}%
\pgfsetdash{}{0pt}%
\pgfsys@defobject{currentmarker}{\pgfqpoint{-0.041667in}{0.000000in}}{\pgfqpoint{0.000000in}{0.000000in}}{%
\pgfpathmoveto{\pgfqpoint{0.000000in}{0.000000in}}%
\pgfpathlineto{\pgfqpoint{-0.041667in}{0.000000in}}%
\pgfusepath{stroke,fill}%
}%
\begin{pgfscope}%
\pgfsys@transformshift{1.125000in}{4.854216in}%
\pgfsys@useobject{currentmarker}{}%
\end{pgfscope}%
\end{pgfscope}%
\begin{pgfscope}%
\pgfsetbuttcap%
\pgfsetroundjoin%
\definecolor{currentfill}{rgb}{0.150000,0.150000,0.150000}%
\pgfsetfillcolor{currentfill}%
\pgfsetlinewidth{0.401500pt}%
\definecolor{currentstroke}{rgb}{0.150000,0.150000,0.150000}%
\pgfsetstrokecolor{currentstroke}%
\pgfsetdash{}{0pt}%
\pgfsys@defobject{currentmarker}{\pgfqpoint{-0.041667in}{0.000000in}}{\pgfqpoint{0.000000in}{0.000000in}}{%
\pgfpathmoveto{\pgfqpoint{0.000000in}{0.000000in}}%
\pgfpathlineto{\pgfqpoint{-0.041667in}{0.000000in}}%
\pgfusepath{stroke,fill}%
}%
\begin{pgfscope}%
\pgfsys@transformshift{1.125000in}{4.970274in}%
\pgfsys@useobject{currentmarker}{}%
\end{pgfscope}%
\end{pgfscope}%
\begin{pgfscope}%
\definecolor{textcolor}{rgb}{0.150000,0.150000,0.150000}%
\pgfsetstrokecolor{textcolor}%
\pgfsetfillcolor{textcolor}%
\pgftext[x=0.383703in,y=3.015000in,,bottom,rotate=90.000000]{\color{textcolor}\sffamily\fontsize{24.000000}{28.800000}\selectfont Spectral radius, \(\displaystyle \rho\)}%
\end{pgfscope}%
\begin{pgfscope}%
\pgfpathrectangle{\pgfqpoint{1.125000in}{0.750000in}}{\pgfqpoint{6.975000in}{4.530000in}} %
\pgfusepath{clip}%
\pgfsetroundcap%
\pgfsetroundjoin%
\pgfsetlinewidth{1.405250pt}%
\definecolor{currentstroke}{rgb}{0.133333,0.133333,0.133333}%
\pgfsetstrokecolor{currentstroke}%
\pgfsetdash{}{0pt}%
\pgfpathmoveto{\pgfqpoint{1.442045in}{3.228359in}}%
\pgfpathlineto{\pgfqpoint{3.525620in}{2.264004in}}%
\pgfpathlineto{\pgfqpoint{4.643254in}{2.300833in}}%
\pgfpathlineto{\pgfqpoint{5.411609in}{2.313436in}}%
\pgfpathlineto{\pgfqpoint{5.997779in}{2.319044in}}%
\pgfpathlineto{\pgfqpoint{6.471812in}{2.321984in}}%
\pgfpathlineto{\pgfqpoint{6.869799in}{2.323706in}}%
\pgfpathlineto{\pgfqpoint{7.212799in}{2.324798in}}%
\pgfpathlineto{\pgfqpoint{7.514179in}{2.325533in}}%
\pgfpathlineto{\pgfqpoint{7.782955in}{2.326052in}}%
\pgfusepath{stroke}%
\end{pgfscope}%
\begin{pgfscope}%
\pgfpathrectangle{\pgfqpoint{1.125000in}{0.750000in}}{\pgfqpoint{6.975000in}{4.530000in}} %
\pgfusepath{clip}%
\pgfsetbuttcap%
\pgfsetmiterjoin%
\definecolor{currentfill}{rgb}{0.133333,0.133333,0.133333}%
\pgfsetfillcolor{currentfill}%
\pgfsetlinewidth{0.000000pt}%
\definecolor{currentstroke}{rgb}{0.133333,0.133333,0.133333}%
\pgfsetstrokecolor{currentstroke}%
\pgfsetdash{}{0pt}%
\pgfsys@defobject{currentmarker}{\pgfqpoint{-0.062500in}{-0.062500in}}{\pgfqpoint{0.062500in}{0.062500in}}{%
\pgfpathmoveto{\pgfqpoint{-0.000000in}{-0.062500in}}%
\pgfpathlineto{\pgfqpoint{0.062500in}{0.062500in}}%
\pgfpathlineto{\pgfqpoint{-0.062500in}{0.062500in}}%
\pgfpathclose%
\pgfusepath{fill}%
}%
\begin{pgfscope}%
\pgfsys@transformshift{1.442045in}{3.228359in}%
\pgfsys@useobject{currentmarker}{}%
\end{pgfscope}%
\begin{pgfscope}%
\pgfsys@transformshift{3.525620in}{2.264004in}%
\pgfsys@useobject{currentmarker}{}%
\end{pgfscope}%
\begin{pgfscope}%
\pgfsys@transformshift{4.643254in}{2.300833in}%
\pgfsys@useobject{currentmarker}{}%
\end{pgfscope}%
\begin{pgfscope}%
\pgfsys@transformshift{5.411609in}{2.313436in}%
\pgfsys@useobject{currentmarker}{}%
\end{pgfscope}%
\begin{pgfscope}%
\pgfsys@transformshift{5.997779in}{2.319044in}%
\pgfsys@useobject{currentmarker}{}%
\end{pgfscope}%
\begin{pgfscope}%
\pgfsys@transformshift{6.471812in}{2.321984in}%
\pgfsys@useobject{currentmarker}{}%
\end{pgfscope}%
\begin{pgfscope}%
\pgfsys@transformshift{6.869799in}{2.323706in}%
\pgfsys@useobject{currentmarker}{}%
\end{pgfscope}%
\begin{pgfscope}%
\pgfsys@transformshift{7.212799in}{2.324798in}%
\pgfsys@useobject{currentmarker}{}%
\end{pgfscope}%
\begin{pgfscope}%
\pgfsys@transformshift{7.514179in}{2.325533in}%
\pgfsys@useobject{currentmarker}{}%
\end{pgfscope}%
\begin{pgfscope}%
\pgfsys@transformshift{7.782955in}{2.326052in}%
\pgfsys@useobject{currentmarker}{}%
\end{pgfscope}%
\end{pgfscope}%
\begin{pgfscope}%
\pgfpathrectangle{\pgfqpoint{1.125000in}{0.750000in}}{\pgfqpoint{6.975000in}{4.530000in}} %
\pgfusepath{clip}%
\pgfsetroundcap%
\pgfsetroundjoin%
\pgfsetlinewidth{1.405250pt}%
\definecolor{currentstroke}{rgb}{0.318370,0.318370,0.318370}%
\pgfsetstrokecolor{currentstroke}%
\pgfsetdash{}{0pt}%
\pgfpathmoveto{\pgfqpoint{1.442045in}{0.955909in}}%
\pgfpathlineto{\pgfqpoint{3.525620in}{2.010261in}}%
\pgfpathlineto{\pgfqpoint{4.643254in}{2.025801in}}%
\pgfpathlineto{\pgfqpoint{5.411609in}{2.031099in}}%
\pgfpathlineto{\pgfqpoint{5.997779in}{1.998478in}}%
\pgfpathlineto{\pgfqpoint{6.471812in}{2.001144in}}%
\pgfpathlineto{\pgfqpoint{6.869799in}{1.970135in}}%
\pgfpathlineto{\pgfqpoint{7.212799in}{1.989631in}}%
\pgfpathlineto{\pgfqpoint{7.514179in}{1.981984in}}%
\pgfpathlineto{\pgfqpoint{7.782955in}{1.981743in}}%
\pgfusepath{stroke}%
\end{pgfscope}%
\begin{pgfscope}%
\pgfpathrectangle{\pgfqpoint{1.125000in}{0.750000in}}{\pgfqpoint{6.975000in}{4.530000in}} %
\pgfusepath{clip}%
\pgfsetbuttcap%
\pgfsetroundjoin%
\definecolor{currentfill}{rgb}{0.318370,0.318370,0.318370}%
\pgfsetfillcolor{currentfill}%
\pgfsetlinewidth{0.000000pt}%
\definecolor{currentstroke}{rgb}{0.318370,0.318370,0.318370}%
\pgfsetstrokecolor{currentstroke}%
\pgfsetdash{}{0pt}%
\pgfsys@defobject{currentmarker}{\pgfqpoint{-0.048611in}{-0.048611in}}{\pgfqpoint{0.048611in}{0.048611in}}{%
\pgfpathmoveto{\pgfqpoint{0.000000in}{-0.048611in}}%
\pgfpathcurveto{\pgfqpoint{0.012892in}{-0.048611in}}{\pgfqpoint{0.025257in}{-0.043489in}}{\pgfqpoint{0.034373in}{-0.034373in}}%
\pgfpathcurveto{\pgfqpoint{0.043489in}{-0.025257in}}{\pgfqpoint{0.048611in}{-0.012892in}}{\pgfqpoint{0.048611in}{0.000000in}}%
\pgfpathcurveto{\pgfqpoint{0.048611in}{0.012892in}}{\pgfqpoint{0.043489in}{0.025257in}}{\pgfqpoint{0.034373in}{0.034373in}}%
\pgfpathcurveto{\pgfqpoint{0.025257in}{0.043489in}}{\pgfqpoint{0.012892in}{0.048611in}}{\pgfqpoint{0.000000in}{0.048611in}}%
\pgfpathcurveto{\pgfqpoint{-0.012892in}{0.048611in}}{\pgfqpoint{-0.025257in}{0.043489in}}{\pgfqpoint{-0.034373in}{0.034373in}}%
\pgfpathcurveto{\pgfqpoint{-0.043489in}{0.025257in}}{\pgfqpoint{-0.048611in}{0.012892in}}{\pgfqpoint{-0.048611in}{0.000000in}}%
\pgfpathcurveto{\pgfqpoint{-0.048611in}{-0.012892in}}{\pgfqpoint{-0.043489in}{-0.025257in}}{\pgfqpoint{-0.034373in}{-0.034373in}}%
\pgfpathcurveto{\pgfqpoint{-0.025257in}{-0.043489in}}{\pgfqpoint{-0.012892in}{-0.048611in}}{\pgfqpoint{0.000000in}{-0.048611in}}%
\pgfpathclose%
\pgfusepath{fill}%
}%
\begin{pgfscope}%
\pgfsys@transformshift{1.442045in}{0.955909in}%
\pgfsys@useobject{currentmarker}{}%
\end{pgfscope}%
\begin{pgfscope}%
\pgfsys@transformshift{3.525620in}{2.010261in}%
\pgfsys@useobject{currentmarker}{}%
\end{pgfscope}%
\begin{pgfscope}%
\pgfsys@transformshift{4.643254in}{2.025801in}%
\pgfsys@useobject{currentmarker}{}%
\end{pgfscope}%
\begin{pgfscope}%
\pgfsys@transformshift{5.411609in}{2.031099in}%
\pgfsys@useobject{currentmarker}{}%
\end{pgfscope}%
\begin{pgfscope}%
\pgfsys@transformshift{5.997779in}{1.998478in}%
\pgfsys@useobject{currentmarker}{}%
\end{pgfscope}%
\begin{pgfscope}%
\pgfsys@transformshift{6.471812in}{2.001144in}%
\pgfsys@useobject{currentmarker}{}%
\end{pgfscope}%
\begin{pgfscope}%
\pgfsys@transformshift{6.869799in}{1.970135in}%
\pgfsys@useobject{currentmarker}{}%
\end{pgfscope}%
\begin{pgfscope}%
\pgfsys@transformshift{7.212799in}{1.989631in}%
\pgfsys@useobject{currentmarker}{}%
\end{pgfscope}%
\begin{pgfscope}%
\pgfsys@transformshift{7.514179in}{1.981984in}%
\pgfsys@useobject{currentmarker}{}%
\end{pgfscope}%
\begin{pgfscope}%
\pgfsys@transformshift{7.782955in}{1.981743in}%
\pgfsys@useobject{currentmarker}{}%
\end{pgfscope}%
\end{pgfscope}%
\begin{pgfscope}%
\pgfpathrectangle{\pgfqpoint{1.125000in}{0.750000in}}{\pgfqpoint{6.975000in}{4.530000in}} %
\pgfusepath{clip}%
\pgfsetroundcap%
\pgfsetroundjoin%
\pgfsetlinewidth{1.405250pt}%
\definecolor{currentstroke}{rgb}{0.501961,0.501961,0.501961}%
\pgfsetstrokecolor{currentstroke}%
\pgfsetdash{}{0pt}%
\pgfpathmoveto{\pgfqpoint{1.442045in}{3.915180in}}%
\pgfpathlineto{\pgfqpoint{3.525620in}{3.459474in}}%
\pgfpathlineto{\pgfqpoint{4.643254in}{3.471959in}}%
\pgfpathlineto{\pgfqpoint{5.411609in}{3.472300in}}%
\pgfpathlineto{\pgfqpoint{5.997779in}{3.471420in}}%
\pgfpathlineto{\pgfqpoint{6.471812in}{3.471775in}}%
\pgfpathlineto{\pgfqpoint{6.869799in}{3.471892in}}%
\pgfpathlineto{\pgfqpoint{7.212799in}{3.473092in}}%
\pgfpathlineto{\pgfqpoint{7.514179in}{3.473717in}}%
\pgfpathlineto{\pgfqpoint{7.782955in}{3.474202in}}%
\pgfusepath{stroke}%
\end{pgfscope}%
\begin{pgfscope}%
\pgfpathrectangle{\pgfqpoint{1.125000in}{0.750000in}}{\pgfqpoint{6.975000in}{4.530000in}} %
\pgfusepath{clip}%
\pgfsetbuttcap%
\pgfsetroundjoin%
\definecolor{currentfill}{rgb}{0.501961,0.501961,0.501961}%
\pgfsetfillcolor{currentfill}%
\pgfsetlinewidth{0.000000pt}%
\definecolor{currentstroke}{rgb}{0.501961,0.501961,0.501961}%
\pgfsetstrokecolor{currentstroke}%
\pgfsetdash{}{0pt}%
\pgfsys@defobject{currentmarker}{\pgfqpoint{-0.034722in}{-0.034722in}}{\pgfqpoint{0.034722in}{0.034722in}}{%
\pgfpathmoveto{\pgfqpoint{0.000000in}{-0.034722in}}%
\pgfpathcurveto{\pgfqpoint{0.009208in}{-0.034722in}}{\pgfqpoint{0.018041in}{-0.031064in}}{\pgfqpoint{0.024552in}{-0.024552in}}%
\pgfpathcurveto{\pgfqpoint{0.031064in}{-0.018041in}}{\pgfqpoint{0.034722in}{-0.009208in}}{\pgfqpoint{0.034722in}{0.000000in}}%
\pgfpathcurveto{\pgfqpoint{0.034722in}{0.009208in}}{\pgfqpoint{0.031064in}{0.018041in}}{\pgfqpoint{0.024552in}{0.024552in}}%
\pgfpathcurveto{\pgfqpoint{0.018041in}{0.031064in}}{\pgfqpoint{0.009208in}{0.034722in}}{\pgfqpoint{0.000000in}{0.034722in}}%
\pgfpathcurveto{\pgfqpoint{-0.009208in}{0.034722in}}{\pgfqpoint{-0.018041in}{0.031064in}}{\pgfqpoint{-0.024552in}{0.024552in}}%
\pgfpathcurveto{\pgfqpoint{-0.031064in}{0.018041in}}{\pgfqpoint{-0.034722in}{0.009208in}}{\pgfqpoint{-0.034722in}{0.000000in}}%
\pgfpathcurveto{\pgfqpoint{-0.034722in}{-0.009208in}}{\pgfqpoint{-0.031064in}{-0.018041in}}{\pgfqpoint{-0.024552in}{-0.024552in}}%
\pgfpathcurveto{\pgfqpoint{-0.018041in}{-0.031064in}}{\pgfqpoint{-0.009208in}{-0.034722in}}{\pgfqpoint{0.000000in}{-0.034722in}}%
\pgfpathclose%
\pgfusepath{fill}%
}%
\begin{pgfscope}%
\pgfsys@transformshift{1.442045in}{3.915180in}%
\pgfsys@useobject{currentmarker}{}%
\end{pgfscope}%
\begin{pgfscope}%
\pgfsys@transformshift{3.525620in}{3.459474in}%
\pgfsys@useobject{currentmarker}{}%
\end{pgfscope}%
\begin{pgfscope}%
\pgfsys@transformshift{4.643254in}{3.471959in}%
\pgfsys@useobject{currentmarker}{}%
\end{pgfscope}%
\begin{pgfscope}%
\pgfsys@transformshift{5.411609in}{3.472300in}%
\pgfsys@useobject{currentmarker}{}%
\end{pgfscope}%
\begin{pgfscope}%
\pgfsys@transformshift{5.997779in}{3.471420in}%
\pgfsys@useobject{currentmarker}{}%
\end{pgfscope}%
\begin{pgfscope}%
\pgfsys@transformshift{6.471812in}{3.471775in}%
\pgfsys@useobject{currentmarker}{}%
\end{pgfscope}%
\begin{pgfscope}%
\pgfsys@transformshift{6.869799in}{3.471892in}%
\pgfsys@useobject{currentmarker}{}%
\end{pgfscope}%
\begin{pgfscope}%
\pgfsys@transformshift{7.212799in}{3.473092in}%
\pgfsys@useobject{currentmarker}{}%
\end{pgfscope}%
\begin{pgfscope}%
\pgfsys@transformshift{7.514179in}{3.473717in}%
\pgfsys@useobject{currentmarker}{}%
\end{pgfscope}%
\begin{pgfscope}%
\pgfsys@transformshift{7.782955in}{3.474202in}%
\pgfsys@useobject{currentmarker}{}%
\end{pgfscope}%
\end{pgfscope}%
\begin{pgfscope}%
\pgfpathrectangle{\pgfqpoint{1.125000in}{0.750000in}}{\pgfqpoint{6.975000in}{4.530000in}} %
\pgfusepath{clip}%
\pgfsetbuttcap%
\pgfsetroundjoin%
\pgfsetlinewidth{1.405250pt}%
\definecolor{currentstroke}{rgb}{0.000000,0.000000,0.000000}%
\pgfsetstrokecolor{currentstroke}%
\pgfsetdash{{5.180000pt}{2.240000pt}}{0.000000pt}%
\pgfpathmoveto{\pgfqpoint{1.442045in}{5.074091in}}%
\pgfpathlineto{\pgfqpoint{3.525620in}{5.074091in}}%
\pgfpathlineto{\pgfqpoint{4.643254in}{5.074091in}}%
\pgfpathlineto{\pgfqpoint{5.411609in}{5.074091in}}%
\pgfpathlineto{\pgfqpoint{5.997779in}{5.074091in}}%
\pgfpathlineto{\pgfqpoint{6.471812in}{5.074091in}}%
\pgfpathlineto{\pgfqpoint{6.869799in}{5.074091in}}%
\pgfpathlineto{\pgfqpoint{7.212799in}{5.074091in}}%
\pgfpathlineto{\pgfqpoint{7.514179in}{5.074091in}}%
\pgfpathlineto{\pgfqpoint{7.782955in}{5.074091in}}%
\pgfusepath{stroke}%
\end{pgfscope}%
\begin{pgfscope}%
\pgfsetrectcap%
\pgfsetmiterjoin%
\pgfsetlinewidth{1.254687pt}%
\definecolor{currentstroke}{rgb}{0.150000,0.150000,0.150000}%
\pgfsetstrokecolor{currentstroke}%
\pgfsetdash{}{0pt}%
\pgfpathmoveto{\pgfqpoint{1.125000in}{0.750000in}}%
\pgfpathlineto{\pgfqpoint{1.125000in}{5.280000in}}%
\pgfusepath{stroke}%
\end{pgfscope}%
\begin{pgfscope}%
\pgfsetrectcap%
\pgfsetmiterjoin%
\pgfsetlinewidth{1.254687pt}%
\definecolor{currentstroke}{rgb}{0.150000,0.150000,0.150000}%
\pgfsetstrokecolor{currentstroke}%
\pgfsetdash{}{0pt}%
\pgfpathmoveto{\pgfqpoint{8.100000in}{0.750000in}}%
\pgfpathlineto{\pgfqpoint{8.100000in}{5.280000in}}%
\pgfusepath{stroke}%
\end{pgfscope}%
\begin{pgfscope}%
\pgfsetrectcap%
\pgfsetmiterjoin%
\pgfsetlinewidth{1.254687pt}%
\definecolor{currentstroke}{rgb}{0.150000,0.150000,0.150000}%
\pgfsetstrokecolor{currentstroke}%
\pgfsetdash{}{0pt}%
\pgfpathmoveto{\pgfqpoint{1.125000in}{0.750000in}}%
\pgfpathlineto{\pgfqpoint{8.100000in}{0.750000in}}%
\pgfusepath{stroke}%
\end{pgfscope}%
\begin{pgfscope}%
\pgfsetrectcap%
\pgfsetmiterjoin%
\pgfsetlinewidth{1.254687pt}%
\definecolor{currentstroke}{rgb}{0.150000,0.150000,0.150000}%
\pgfsetstrokecolor{currentstroke}%
\pgfsetdash{}{0pt}%
\pgfpathmoveto{\pgfqpoint{1.125000in}{5.280000in}}%
\pgfpathlineto{\pgfqpoint{8.100000in}{5.280000in}}%
\pgfusepath{stroke}%
\end{pgfscope}%
\begin{pgfscope}%
\pgfsetroundcap%
\pgfsetroundjoin%
\pgfsetlinewidth{1.405250pt}%
\definecolor{currentstroke}{rgb}{0.133333,0.133333,0.133333}%
\pgfsetstrokecolor{currentstroke}%
\pgfsetdash{}{0pt}%
\pgfpathmoveto{\pgfqpoint{5.448533in}{4.780632in}}%
\pgfpathlineto{\pgfqpoint{5.865200in}{4.780632in}}%
\pgfusepath{stroke}%
\end{pgfscope}%
\begin{pgfscope}%
\pgfsetbuttcap%
\pgfsetmiterjoin%
\definecolor{currentfill}{rgb}{0.133333,0.133333,0.133333}%
\pgfsetfillcolor{currentfill}%
\pgfsetlinewidth{0.000000pt}%
\definecolor{currentstroke}{rgb}{0.133333,0.133333,0.133333}%
\pgfsetstrokecolor{currentstroke}%
\pgfsetdash{}{0pt}%
\pgfsys@defobject{currentmarker}{\pgfqpoint{-0.062500in}{-0.062500in}}{\pgfqpoint{0.062500in}{0.062500in}}{%
\pgfpathmoveto{\pgfqpoint{-0.000000in}{-0.062500in}}%
\pgfpathlineto{\pgfqpoint{0.062500in}{0.062500in}}%
\pgfpathlineto{\pgfqpoint{-0.062500in}{0.062500in}}%
\pgfpathclose%
\pgfusepath{fill}%
}%
\begin{pgfscope}%
\pgfsys@transformshift{5.656867in}{4.780632in}%
\pgfsys@useobject{currentmarker}{}%
\end{pgfscope}%
\end{pgfscope}%
\begin{pgfscope}%
\definecolor{textcolor}{rgb}{0.150000,0.150000,0.150000}%
\pgfsetstrokecolor{textcolor}%
\pgfsetfillcolor{textcolor}%
\pgftext[x=6.031867in,y=4.707716in,left,base]{\color{textcolor}\sffamily\fontsize{15.000000}{18.000000}\selectfont Linear}%
\end{pgfscope}%
\begin{pgfscope}%
\pgfsetroundcap%
\pgfsetroundjoin%
\pgfsetlinewidth{1.405250pt}%
\definecolor{currentstroke}{rgb}{0.318370,0.318370,0.318370}%
\pgfsetstrokecolor{currentstroke}%
\pgfsetdash{}{0pt}%
\pgfpathmoveto{\pgfqpoint{5.448533in}{4.474846in}}%
\pgfpathlineto{\pgfqpoint{5.865200in}{4.474846in}}%
\pgfusepath{stroke}%
\end{pgfscope}%
\begin{pgfscope}%
\pgfsetbuttcap%
\pgfsetroundjoin%
\definecolor{currentfill}{rgb}{0.318370,0.318370,0.318370}%
\pgfsetfillcolor{currentfill}%
\pgfsetlinewidth{0.000000pt}%
\definecolor{currentstroke}{rgb}{0.318370,0.318370,0.318370}%
\pgfsetstrokecolor{currentstroke}%
\pgfsetdash{}{0pt}%
\pgfsys@defobject{currentmarker}{\pgfqpoint{-0.048611in}{-0.048611in}}{\pgfqpoint{0.048611in}{0.048611in}}{%
\pgfpathmoveto{\pgfqpoint{0.000000in}{-0.048611in}}%
\pgfpathcurveto{\pgfqpoint{0.012892in}{-0.048611in}}{\pgfqpoint{0.025257in}{-0.043489in}}{\pgfqpoint{0.034373in}{-0.034373in}}%
\pgfpathcurveto{\pgfqpoint{0.043489in}{-0.025257in}}{\pgfqpoint{0.048611in}{-0.012892in}}{\pgfqpoint{0.048611in}{0.000000in}}%
\pgfpathcurveto{\pgfqpoint{0.048611in}{0.012892in}}{\pgfqpoint{0.043489in}{0.025257in}}{\pgfqpoint{0.034373in}{0.034373in}}%
\pgfpathcurveto{\pgfqpoint{0.025257in}{0.043489in}}{\pgfqpoint{0.012892in}{0.048611in}}{\pgfqpoint{0.000000in}{0.048611in}}%
\pgfpathcurveto{\pgfqpoint{-0.012892in}{0.048611in}}{\pgfqpoint{-0.025257in}{0.043489in}}{\pgfqpoint{-0.034373in}{0.034373in}}%
\pgfpathcurveto{\pgfqpoint{-0.043489in}{0.025257in}}{\pgfqpoint{-0.048611in}{0.012892in}}{\pgfqpoint{-0.048611in}{0.000000in}}%
\pgfpathcurveto{\pgfqpoint{-0.048611in}{-0.012892in}}{\pgfqpoint{-0.043489in}{-0.025257in}}{\pgfqpoint{-0.034373in}{-0.034373in}}%
\pgfpathcurveto{\pgfqpoint{-0.025257in}{-0.043489in}}{\pgfqpoint{-0.012892in}{-0.048611in}}{\pgfqpoint{0.000000in}{-0.048611in}}%
\pgfpathclose%
\pgfusepath{fill}%
}%
\begin{pgfscope}%
\pgfsys@transformshift{5.656867in}{4.474846in}%
\pgfsys@useobject{currentmarker}{}%
\end{pgfscope}%
\end{pgfscope}%
\begin{pgfscope}%
\definecolor{textcolor}{rgb}{0.150000,0.150000,0.150000}%
\pgfsetstrokecolor{textcolor}%
\pgfsetfillcolor{textcolor}%
\pgftext[x=6.031867in,y=4.401930in,left,base]{\color{textcolor}\sffamily\fontsize{15.000000}{18.000000}\selectfont DMG}%
\end{pgfscope}%
\begin{pgfscope}%
\pgfsetroundcap%
\pgfsetroundjoin%
\pgfsetlinewidth{1.405250pt}%
\definecolor{currentstroke}{rgb}{0.501961,0.501961,0.501961}%
\pgfsetstrokecolor{currentstroke}%
\pgfsetdash{}{0pt}%
\pgfpathmoveto{\pgfqpoint{5.448533in}{4.169061in}}%
\pgfpathlineto{\pgfqpoint{5.865200in}{4.169061in}}%
\pgfusepath{stroke}%
\end{pgfscope}%
\begin{pgfscope}%
\pgfsetbuttcap%
\pgfsetroundjoin%
\definecolor{currentfill}{rgb}{0.501961,0.501961,0.501961}%
\pgfsetfillcolor{currentfill}%
\pgfsetlinewidth{0.000000pt}%
\definecolor{currentstroke}{rgb}{0.501961,0.501961,0.501961}%
\pgfsetstrokecolor{currentstroke}%
\pgfsetdash{}{0pt}%
\pgfsys@defobject{currentmarker}{\pgfqpoint{-0.034722in}{-0.034722in}}{\pgfqpoint{0.034722in}{0.034722in}}{%
\pgfpathmoveto{\pgfqpoint{0.000000in}{-0.034722in}}%
\pgfpathcurveto{\pgfqpoint{0.009208in}{-0.034722in}}{\pgfqpoint{0.018041in}{-0.031064in}}{\pgfqpoint{0.024552in}{-0.024552in}}%
\pgfpathcurveto{\pgfqpoint{0.031064in}{-0.018041in}}{\pgfqpoint{0.034722in}{-0.009208in}}{\pgfqpoint{0.034722in}{0.000000in}}%
\pgfpathcurveto{\pgfqpoint{0.034722in}{0.009208in}}{\pgfqpoint{0.031064in}{0.018041in}}{\pgfqpoint{0.024552in}{0.024552in}}%
\pgfpathcurveto{\pgfqpoint{0.018041in}{0.031064in}}{\pgfqpoint{0.009208in}{0.034722in}}{\pgfqpoint{0.000000in}{0.034722in}}%
\pgfpathcurveto{\pgfqpoint{-0.009208in}{0.034722in}}{\pgfqpoint{-0.018041in}{0.031064in}}{\pgfqpoint{-0.024552in}{0.024552in}}%
\pgfpathcurveto{\pgfqpoint{-0.031064in}{0.018041in}}{\pgfqpoint{-0.034722in}{0.009208in}}{\pgfqpoint{-0.034722in}{0.000000in}}%
\pgfpathcurveto{\pgfqpoint{-0.034722in}{-0.009208in}}{\pgfqpoint{-0.031064in}{-0.018041in}}{\pgfqpoint{-0.024552in}{-0.024552in}}%
\pgfpathcurveto{\pgfqpoint{-0.018041in}{-0.031064in}}{\pgfqpoint{-0.009208in}{-0.034722in}}{\pgfqpoint{0.000000in}{-0.034722in}}%
\pgfpathclose%
\pgfusepath{fill}%
}%
\begin{pgfscope}%
\pgfsys@transformshift{5.656867in}{4.169061in}%
\pgfsys@useobject{currentmarker}{}%
\end{pgfscope}%
\end{pgfscope}%
\begin{pgfscope}%
\definecolor{textcolor}{rgb}{0.150000,0.150000,0.150000}%
\pgfsetstrokecolor{textcolor}%
\pgfsetfillcolor{textcolor}%
\pgftext[x=6.031867in,y=4.096144in,left,base]{\color{textcolor}\sffamily\fontsize{15.000000}{18.000000}\selectfont AMG}%
\end{pgfscope}%
\begin{pgfscope}%
\pgfsetbuttcap%
\pgfsetroundjoin%
\pgfsetlinewidth{1.405250pt}%
\definecolor{currentstroke}{rgb}{0.000000,0.000000,0.000000}%
\pgfsetstrokecolor{currentstroke}%
\pgfsetdash{{5.180000pt}{2.240000pt}}{0.000000pt}%
\pgfpathmoveto{\pgfqpoint{5.448533in}{3.863275in}}%
\pgfpathlineto{\pgfqpoint{5.865200in}{3.863275in}}%
\pgfusepath{stroke}%
\end{pgfscope}%
\begin{pgfscope}%
\definecolor{textcolor}{rgb}{0.150000,0.150000,0.150000}%
\pgfsetstrokecolor{textcolor}%
\pgfsetfillcolor{textcolor}%
\pgftext[x=6.031867in,y=3.790358in,left,base]{\color{textcolor}\sffamily\fontsize{15.000000}{18.000000}\selectfont Boundary \(\displaystyle \rho = 1\)}%
\end{pgfscope}%
\end{pgfpicture}%
\makeatother%
\endgroup%

%% file: diff_m3points_epsrange_n_255_bw.pgf
\begingroup%
\makeatletter%
\begin{pgfpicture}%
\pgfpathrectangle{\pgfpointorigin}{\pgfqpoint{9.000000in}{6.000000in}}%
\pgfusepath{use as bounding box, clip}%
\begin{pgfscope}%
\pgfsetbuttcap%
\pgfsetmiterjoin%
\definecolor{currentfill}{rgb}{1.000000,1.000000,1.000000}%
\pgfsetfillcolor{currentfill}%
\pgfsetlinewidth{0.000000pt}%
\definecolor{currentstroke}{rgb}{1.000000,1.000000,1.000000}%
\pgfsetstrokecolor{currentstroke}%
\pgfsetdash{}{0pt}%
\pgfpathmoveto{\pgfqpoint{0.000000in}{0.000000in}}%
\pgfpathlineto{\pgfqpoint{9.000000in}{0.000000in}}%
\pgfpathlineto{\pgfqpoint{9.000000in}{6.000000in}}%
\pgfpathlineto{\pgfqpoint{0.000000in}{6.000000in}}%
\pgfpathclose%
\pgfusepath{fill}%
\end{pgfscope}%
\begin{pgfscope}%
\pgfsetbuttcap%
\pgfsetmiterjoin%
\definecolor{currentfill}{rgb}{1.000000,1.000000,1.000000}%
\pgfsetfillcolor{currentfill}%
\pgfsetlinewidth{0.000000pt}%
\definecolor{currentstroke}{rgb}{0.000000,0.000000,0.000000}%
\pgfsetstrokecolor{currentstroke}%
\pgfsetstrokeopacity{0.000000}%
\pgfsetdash{}{0pt}%
\pgfpathmoveto{\pgfqpoint{1.125000in}{0.750000in}}%
\pgfpathlineto{\pgfqpoint{8.100000in}{0.750000in}}%
\pgfpathlineto{\pgfqpoint{8.100000in}{5.280000in}}%
\pgfpathlineto{\pgfqpoint{1.125000in}{5.280000in}}%
\pgfpathclose%
\pgfusepath{fill}%
\end{pgfscope}%
\begin{pgfscope}%
\pgfsetbuttcap%
\pgfsetroundjoin%
\definecolor{currentfill}{rgb}{0.150000,0.150000,0.150000}%
\pgfsetfillcolor{currentfill}%
\pgfsetlinewidth{0.803000pt}%
\definecolor{currentstroke}{rgb}{0.150000,0.150000,0.150000}%
\pgfsetstrokecolor{currentstroke}%
\pgfsetdash{}{0pt}%
\pgfsys@defobject{currentmarker}{\pgfqpoint{0.000000in}{-0.083333in}}{\pgfqpoint{0.000000in}{0.000000in}}{%
\pgfpathmoveto{\pgfqpoint{0.000000in}{0.000000in}}%
\pgfpathlineto{\pgfqpoint{0.000000in}{-0.083333in}}%
\pgfusepath{stroke,fill}%
}%
\begin{pgfscope}%
\pgfsys@transformshift{3.280235in}{0.750000in}%
\pgfsys@useobject{currentmarker}{}%
\end{pgfscope}%
\end{pgfscope}%
\begin{pgfscope}%
\definecolor{textcolor}{rgb}{0.150000,0.150000,0.150000}%
\pgfsetstrokecolor{textcolor}%
\pgfsetfillcolor{textcolor}%
\pgftext[x=3.280235in,y=0.588889in,,top]{\color{textcolor}\sffamily\fontsize{16.000000}{19.200000}\selectfont \(\displaystyle 10^{-2}\)}%
\end{pgfscope}%
\begin{pgfscope}%
\pgfsetbuttcap%
\pgfsetroundjoin%
\definecolor{currentfill}{rgb}{0.150000,0.150000,0.150000}%
\pgfsetfillcolor{currentfill}%
\pgfsetlinewidth{0.803000pt}%
\definecolor{currentstroke}{rgb}{0.150000,0.150000,0.150000}%
\pgfsetstrokecolor{currentstroke}%
\pgfsetdash{}{0pt}%
\pgfsys@defobject{currentmarker}{\pgfqpoint{0.000000in}{-0.083333in}}{\pgfqpoint{0.000000in}{0.000000in}}{%
\pgfpathmoveto{\pgfqpoint{0.000000in}{0.000000in}}%
\pgfpathlineto{\pgfqpoint{0.000000in}{-0.083333in}}%
\pgfusepath{stroke,fill}%
}%
\begin{pgfscope}%
\pgfsys@transformshift{7.782955in}{0.750000in}%
\pgfsys@useobject{currentmarker}{}%
\end{pgfscope}%
\end{pgfscope}%
\begin{pgfscope}%
\definecolor{textcolor}{rgb}{0.150000,0.150000,0.150000}%
\pgfsetstrokecolor{textcolor}%
\pgfsetfillcolor{textcolor}%
\pgftext[x=7.782955in,y=0.588889in,,top]{\color{textcolor}\sffamily\fontsize{16.000000}{19.200000}\selectfont \(\displaystyle 10^{-1}\)}%
\end{pgfscope}%
\begin{pgfscope}%
\pgfsetbuttcap%
\pgfsetroundjoin%
\definecolor{currentfill}{rgb}{0.150000,0.150000,0.150000}%
\pgfsetfillcolor{currentfill}%
\pgfsetlinewidth{0.401500pt}%
\definecolor{currentstroke}{rgb}{0.150000,0.150000,0.150000}%
\pgfsetstrokecolor{currentstroke}%
\pgfsetdash{}{0pt}%
\pgfsys@defobject{currentmarker}{\pgfqpoint{0.000000in}{-0.041667in}}{\pgfqpoint{0.000000in}{0.000000in}}{%
\pgfpathmoveto{\pgfqpoint{0.000000in}{0.000000in}}%
\pgfpathlineto{\pgfqpoint{0.000000in}{-0.041667in}}%
\pgfusepath{stroke,fill}%
}%
\begin{pgfscope}%
\pgfsys@transformshift{1.488423in}{0.750000in}%
\pgfsys@useobject{currentmarker}{}%
\end{pgfscope}%
\end{pgfscope}%
\begin{pgfscope}%
\pgfsetbuttcap%
\pgfsetroundjoin%
\definecolor{currentfill}{rgb}{0.150000,0.150000,0.150000}%
\pgfsetfillcolor{currentfill}%
\pgfsetlinewidth{0.401500pt}%
\definecolor{currentstroke}{rgb}{0.150000,0.150000,0.150000}%
\pgfsetstrokecolor{currentstroke}%
\pgfsetdash{}{0pt}%
\pgfsys@defobject{currentmarker}{\pgfqpoint{0.000000in}{-0.041667in}}{\pgfqpoint{0.000000in}{0.000000in}}{%
\pgfpathmoveto{\pgfqpoint{0.000000in}{0.000000in}}%
\pgfpathlineto{\pgfqpoint{0.000000in}{-0.041667in}}%
\pgfusepath{stroke,fill}%
}%
\begin{pgfscope}%
\pgfsys@transformshift{1.924782in}{0.750000in}%
\pgfsys@useobject{currentmarker}{}%
\end{pgfscope}%
\end{pgfscope}%
\begin{pgfscope}%
\pgfsetbuttcap%
\pgfsetroundjoin%
\definecolor{currentfill}{rgb}{0.150000,0.150000,0.150000}%
\pgfsetfillcolor{currentfill}%
\pgfsetlinewidth{0.401500pt}%
\definecolor{currentstroke}{rgb}{0.150000,0.150000,0.150000}%
\pgfsetstrokecolor{currentstroke}%
\pgfsetdash{}{0pt}%
\pgfsys@defobject{currentmarker}{\pgfqpoint{0.000000in}{-0.041667in}}{\pgfqpoint{0.000000in}{0.000000in}}{%
\pgfpathmoveto{\pgfqpoint{0.000000in}{0.000000in}}%
\pgfpathlineto{\pgfqpoint{0.000000in}{-0.041667in}}%
\pgfusepath{stroke,fill}%
}%
\begin{pgfscope}%
\pgfsys@transformshift{2.281313in}{0.750000in}%
\pgfsys@useobject{currentmarker}{}%
\end{pgfscope}%
\end{pgfscope}%
\begin{pgfscope}%
\pgfsetbuttcap%
\pgfsetroundjoin%
\definecolor{currentfill}{rgb}{0.150000,0.150000,0.150000}%
\pgfsetfillcolor{currentfill}%
\pgfsetlinewidth{0.401500pt}%
\definecolor{currentstroke}{rgb}{0.150000,0.150000,0.150000}%
\pgfsetstrokecolor{currentstroke}%
\pgfsetdash{}{0pt}%
\pgfsys@defobject{currentmarker}{\pgfqpoint{0.000000in}{-0.041667in}}{\pgfqpoint{0.000000in}{0.000000in}}{%
\pgfpathmoveto{\pgfqpoint{0.000000in}{0.000000in}}%
\pgfpathlineto{\pgfqpoint{0.000000in}{-0.041667in}}%
\pgfusepath{stroke,fill}%
}%
\begin{pgfscope}%
\pgfsys@transformshift{2.582755in}{0.750000in}%
\pgfsys@useobject{currentmarker}{}%
\end{pgfscope}%
\end{pgfscope}%
\begin{pgfscope}%
\pgfsetbuttcap%
\pgfsetroundjoin%
\definecolor{currentfill}{rgb}{0.150000,0.150000,0.150000}%
\pgfsetfillcolor{currentfill}%
\pgfsetlinewidth{0.401500pt}%
\definecolor{currentstroke}{rgb}{0.150000,0.150000,0.150000}%
\pgfsetstrokecolor{currentstroke}%
\pgfsetdash{}{0pt}%
\pgfsys@defobject{currentmarker}{\pgfqpoint{0.000000in}{-0.041667in}}{\pgfqpoint{0.000000in}{0.000000in}}{%
\pgfpathmoveto{\pgfqpoint{0.000000in}{0.000000in}}%
\pgfpathlineto{\pgfqpoint{0.000000in}{-0.041667in}}%
\pgfusepath{stroke,fill}%
}%
\begin{pgfscope}%
\pgfsys@transformshift{2.843877in}{0.750000in}%
\pgfsys@useobject{currentmarker}{}%
\end{pgfscope}%
\end{pgfscope}%
\begin{pgfscope}%
\pgfsetbuttcap%
\pgfsetroundjoin%
\definecolor{currentfill}{rgb}{0.150000,0.150000,0.150000}%
\pgfsetfillcolor{currentfill}%
\pgfsetlinewidth{0.401500pt}%
\definecolor{currentstroke}{rgb}{0.150000,0.150000,0.150000}%
\pgfsetstrokecolor{currentstroke}%
\pgfsetdash{}{0pt}%
\pgfsys@defobject{currentmarker}{\pgfqpoint{0.000000in}{-0.041667in}}{\pgfqpoint{0.000000in}{0.000000in}}{%
\pgfpathmoveto{\pgfqpoint{0.000000in}{0.000000in}}%
\pgfpathlineto{\pgfqpoint{0.000000in}{-0.041667in}}%
\pgfusepath{stroke,fill}%
}%
\begin{pgfscope}%
\pgfsys@transformshift{3.074202in}{0.750000in}%
\pgfsys@useobject{currentmarker}{}%
\end{pgfscope}%
\end{pgfscope}%
\begin{pgfscope}%
\pgfsetbuttcap%
\pgfsetroundjoin%
\definecolor{currentfill}{rgb}{0.150000,0.150000,0.150000}%
\pgfsetfillcolor{currentfill}%
\pgfsetlinewidth{0.401500pt}%
\definecolor{currentstroke}{rgb}{0.150000,0.150000,0.150000}%
\pgfsetstrokecolor{currentstroke}%
\pgfsetdash{}{0pt}%
\pgfsys@defobject{currentmarker}{\pgfqpoint{0.000000in}{-0.041667in}}{\pgfqpoint{0.000000in}{0.000000in}}{%
\pgfpathmoveto{\pgfqpoint{0.000000in}{0.000000in}}%
\pgfpathlineto{\pgfqpoint{0.000000in}{-0.041667in}}%
\pgfusepath{stroke,fill}%
}%
\begin{pgfscope}%
\pgfsys@transformshift{4.635689in}{0.750000in}%
\pgfsys@useobject{currentmarker}{}%
\end{pgfscope}%
\end{pgfscope}%
\begin{pgfscope}%
\pgfsetbuttcap%
\pgfsetroundjoin%
\definecolor{currentfill}{rgb}{0.150000,0.150000,0.150000}%
\pgfsetfillcolor{currentfill}%
\pgfsetlinewidth{0.401500pt}%
\definecolor{currentstroke}{rgb}{0.150000,0.150000,0.150000}%
\pgfsetstrokecolor{currentstroke}%
\pgfsetdash{}{0pt}%
\pgfsys@defobject{currentmarker}{\pgfqpoint{0.000000in}{-0.041667in}}{\pgfqpoint{0.000000in}{0.000000in}}{%
\pgfpathmoveto{\pgfqpoint{0.000000in}{0.000000in}}%
\pgfpathlineto{\pgfqpoint{0.000000in}{-0.041667in}}%
\pgfusepath{stroke,fill}%
}%
\begin{pgfscope}%
\pgfsys@transformshift{5.428578in}{0.750000in}%
\pgfsys@useobject{currentmarker}{}%
\end{pgfscope}%
\end{pgfscope}%
\begin{pgfscope}%
\pgfsetbuttcap%
\pgfsetroundjoin%
\definecolor{currentfill}{rgb}{0.150000,0.150000,0.150000}%
\pgfsetfillcolor{currentfill}%
\pgfsetlinewidth{0.401500pt}%
\definecolor{currentstroke}{rgb}{0.150000,0.150000,0.150000}%
\pgfsetstrokecolor{currentstroke}%
\pgfsetdash{}{0pt}%
\pgfsys@defobject{currentmarker}{\pgfqpoint{0.000000in}{-0.041667in}}{\pgfqpoint{0.000000in}{0.000000in}}{%
\pgfpathmoveto{\pgfqpoint{0.000000in}{0.000000in}}%
\pgfpathlineto{\pgfqpoint{0.000000in}{-0.041667in}}%
\pgfusepath{stroke,fill}%
}%
\begin{pgfscope}%
\pgfsys@transformshift{5.991142in}{0.750000in}%
\pgfsys@useobject{currentmarker}{}%
\end{pgfscope}%
\end{pgfscope}%
\begin{pgfscope}%
\pgfsetbuttcap%
\pgfsetroundjoin%
\definecolor{currentfill}{rgb}{0.150000,0.150000,0.150000}%
\pgfsetfillcolor{currentfill}%
\pgfsetlinewidth{0.401500pt}%
\definecolor{currentstroke}{rgb}{0.150000,0.150000,0.150000}%
\pgfsetstrokecolor{currentstroke}%
\pgfsetdash{}{0pt}%
\pgfsys@defobject{currentmarker}{\pgfqpoint{0.000000in}{-0.041667in}}{\pgfqpoint{0.000000in}{0.000000in}}{%
\pgfpathmoveto{\pgfqpoint{0.000000in}{0.000000in}}%
\pgfpathlineto{\pgfqpoint{0.000000in}{-0.041667in}}%
\pgfusepath{stroke,fill}%
}%
\begin{pgfscope}%
\pgfsys@transformshift{6.427501in}{0.750000in}%
\pgfsys@useobject{currentmarker}{}%
\end{pgfscope}%
\end{pgfscope}%
\begin{pgfscope}%
\pgfsetbuttcap%
\pgfsetroundjoin%
\definecolor{currentfill}{rgb}{0.150000,0.150000,0.150000}%
\pgfsetfillcolor{currentfill}%
\pgfsetlinewidth{0.401500pt}%
\definecolor{currentstroke}{rgb}{0.150000,0.150000,0.150000}%
\pgfsetstrokecolor{currentstroke}%
\pgfsetdash{}{0pt}%
\pgfsys@defobject{currentmarker}{\pgfqpoint{0.000000in}{-0.041667in}}{\pgfqpoint{0.000000in}{0.000000in}}{%
\pgfpathmoveto{\pgfqpoint{0.000000in}{0.000000in}}%
\pgfpathlineto{\pgfqpoint{0.000000in}{-0.041667in}}%
\pgfusepath{stroke,fill}%
}%
\begin{pgfscope}%
\pgfsys@transformshift{6.784032in}{0.750000in}%
\pgfsys@useobject{currentmarker}{}%
\end{pgfscope}%
\end{pgfscope}%
\begin{pgfscope}%
\pgfsetbuttcap%
\pgfsetroundjoin%
\definecolor{currentfill}{rgb}{0.150000,0.150000,0.150000}%
\pgfsetfillcolor{currentfill}%
\pgfsetlinewidth{0.401500pt}%
\definecolor{currentstroke}{rgb}{0.150000,0.150000,0.150000}%
\pgfsetstrokecolor{currentstroke}%
\pgfsetdash{}{0pt}%
\pgfsys@defobject{currentmarker}{\pgfqpoint{0.000000in}{-0.041667in}}{\pgfqpoint{0.000000in}{0.000000in}}{%
\pgfpathmoveto{\pgfqpoint{0.000000in}{0.000000in}}%
\pgfpathlineto{\pgfqpoint{0.000000in}{-0.041667in}}%
\pgfusepath{stroke,fill}%
}%
\begin{pgfscope}%
\pgfsys@transformshift{7.085475in}{0.750000in}%
\pgfsys@useobject{currentmarker}{}%
\end{pgfscope}%
\end{pgfscope}%
\begin{pgfscope}%
\pgfsetbuttcap%
\pgfsetroundjoin%
\definecolor{currentfill}{rgb}{0.150000,0.150000,0.150000}%
\pgfsetfillcolor{currentfill}%
\pgfsetlinewidth{0.401500pt}%
\definecolor{currentstroke}{rgb}{0.150000,0.150000,0.150000}%
\pgfsetstrokecolor{currentstroke}%
\pgfsetdash{}{0pt}%
\pgfsys@defobject{currentmarker}{\pgfqpoint{0.000000in}{-0.041667in}}{\pgfqpoint{0.000000in}{0.000000in}}{%
\pgfpathmoveto{\pgfqpoint{0.000000in}{0.000000in}}%
\pgfpathlineto{\pgfqpoint{0.000000in}{-0.041667in}}%
\pgfusepath{stroke,fill}%
}%
\begin{pgfscope}%
\pgfsys@transformshift{7.346596in}{0.750000in}%
\pgfsys@useobject{currentmarker}{}%
\end{pgfscope}%
\end{pgfscope}%
\begin{pgfscope}%
\pgfsetbuttcap%
\pgfsetroundjoin%
\definecolor{currentfill}{rgb}{0.150000,0.150000,0.150000}%
\pgfsetfillcolor{currentfill}%
\pgfsetlinewidth{0.401500pt}%
\definecolor{currentstroke}{rgb}{0.150000,0.150000,0.150000}%
\pgfsetstrokecolor{currentstroke}%
\pgfsetdash{}{0pt}%
\pgfsys@defobject{currentmarker}{\pgfqpoint{0.000000in}{-0.041667in}}{\pgfqpoint{0.000000in}{0.000000in}}{%
\pgfpathmoveto{\pgfqpoint{0.000000in}{0.000000in}}%
\pgfpathlineto{\pgfqpoint{0.000000in}{-0.041667in}}%
\pgfusepath{stroke,fill}%
}%
\begin{pgfscope}%
\pgfsys@transformshift{7.576921in}{0.750000in}%
\pgfsys@useobject{currentmarker}{}%
\end{pgfscope}%
\end{pgfscope}%
\begin{pgfscope}%
\definecolor{textcolor}{rgb}{0.150000,0.150000,0.150000}%
\pgfsetstrokecolor{textcolor}%
\pgfsetfillcolor{textcolor}%
\pgftext[x=4.612500in,y=0.318273in,,top]{\color{textcolor}\sffamily\fontsize{24.000000}{28.800000}\selectfont \(\displaystyle \varepsilon\)}%
\end{pgfscope}%
\begin{pgfscope}%
\pgfsetbuttcap%
\pgfsetroundjoin%
\definecolor{currentfill}{rgb}{0.150000,0.150000,0.150000}%
\pgfsetfillcolor{currentfill}%
\pgfsetlinewidth{0.803000pt}%
\definecolor{currentstroke}{rgb}{0.150000,0.150000,0.150000}%
\pgfsetstrokecolor{currentstroke}%
\pgfsetdash{}{0pt}%
\pgfsys@defobject{currentmarker}{\pgfqpoint{-0.083333in}{0.000000in}}{\pgfqpoint{0.000000in}{0.000000in}}{%
\pgfpathmoveto{\pgfqpoint{0.000000in}{0.000000in}}%
\pgfpathlineto{\pgfqpoint{-0.083333in}{0.000000in}}%
\pgfusepath{stroke,fill}%
}%
\begin{pgfscope}%
\pgfsys@transformshift{1.125000in}{2.725974in}%
\pgfsys@useobject{currentmarker}{}%
\end{pgfscope}%
\end{pgfscope}%
\begin{pgfscope}%
\definecolor{textcolor}{rgb}{0.150000,0.150000,0.150000}%
\pgfsetstrokecolor{textcolor}%
\pgfsetfillcolor{textcolor}%
\pgftext[x=0.439258in,y=2.609899in,left,base]{\color{textcolor}\sffamily\fontsize{22.000000}{26.400000}\selectfont \(\displaystyle 10^{-1}\)}%
\end{pgfscope}%
\begin{pgfscope}%
\pgfsetbuttcap%
\pgfsetroundjoin%
\definecolor{currentfill}{rgb}{0.150000,0.150000,0.150000}%
\pgfsetfillcolor{currentfill}%
\pgfsetlinewidth{0.803000pt}%
\definecolor{currentstroke}{rgb}{0.150000,0.150000,0.150000}%
\pgfsetstrokecolor{currentstroke}%
\pgfsetdash{}{0pt}%
\pgfsys@defobject{currentmarker}{\pgfqpoint{-0.083333in}{0.000000in}}{\pgfqpoint{0.000000in}{0.000000in}}{%
\pgfpathmoveto{\pgfqpoint{0.000000in}{0.000000in}}%
\pgfpathlineto{\pgfqpoint{-0.083333in}{0.000000in}}%
\pgfusepath{stroke,fill}%
}%
\begin{pgfscope}%
\pgfsys@transformshift{1.125000in}{5.074091in}%
\pgfsys@useobject{currentmarker}{}%
\end{pgfscope}%
\end{pgfscope}%
\begin{pgfscope}%
\definecolor{textcolor}{rgb}{0.150000,0.150000,0.150000}%
\pgfsetstrokecolor{textcolor}%
\pgfsetfillcolor{textcolor}%
\pgftext[x=0.594814in,y=4.958016in,left,base]{\color{textcolor}\sffamily\fontsize{22.000000}{26.400000}\selectfont \(\displaystyle 10^{0}\)}%
\end{pgfscope}%
\begin{pgfscope}%
\pgfsetbuttcap%
\pgfsetroundjoin%
\definecolor{currentfill}{rgb}{0.150000,0.150000,0.150000}%
\pgfsetfillcolor{currentfill}%
\pgfsetlinewidth{0.401500pt}%
\definecolor{currentstroke}{rgb}{0.150000,0.150000,0.150000}%
\pgfsetstrokecolor{currentstroke}%
\pgfsetdash{}{0pt}%
\pgfsys@defobject{currentmarker}{\pgfqpoint{-0.041667in}{0.000000in}}{\pgfqpoint{0.000000in}{0.000000in}}{%
\pgfpathmoveto{\pgfqpoint{0.000000in}{0.000000in}}%
\pgfpathlineto{\pgfqpoint{-0.041667in}{0.000000in}}%
\pgfusepath{stroke,fill}%
}%
\begin{pgfscope}%
\pgfsys@transformshift{1.125000in}{1.084711in}%
\pgfsys@useobject{currentmarker}{}%
\end{pgfscope}%
\end{pgfscope}%
\begin{pgfscope}%
\pgfsetbuttcap%
\pgfsetroundjoin%
\definecolor{currentfill}{rgb}{0.150000,0.150000,0.150000}%
\pgfsetfillcolor{currentfill}%
\pgfsetlinewidth{0.401500pt}%
\definecolor{currentstroke}{rgb}{0.150000,0.150000,0.150000}%
\pgfsetstrokecolor{currentstroke}%
\pgfsetdash{}{0pt}%
\pgfsys@defobject{currentmarker}{\pgfqpoint{-0.041667in}{0.000000in}}{\pgfqpoint{0.000000in}{0.000000in}}{%
\pgfpathmoveto{\pgfqpoint{0.000000in}{0.000000in}}%
\pgfpathlineto{\pgfqpoint{-0.041667in}{0.000000in}}%
\pgfusepath{stroke,fill}%
}%
\begin{pgfscope}%
\pgfsys@transformshift{1.125000in}{1.498194in}%
\pgfsys@useobject{currentmarker}{}%
\end{pgfscope}%
\end{pgfscope}%
\begin{pgfscope}%
\pgfsetbuttcap%
\pgfsetroundjoin%
\definecolor{currentfill}{rgb}{0.150000,0.150000,0.150000}%
\pgfsetfillcolor{currentfill}%
\pgfsetlinewidth{0.401500pt}%
\definecolor{currentstroke}{rgb}{0.150000,0.150000,0.150000}%
\pgfsetstrokecolor{currentstroke}%
\pgfsetdash{}{0pt}%
\pgfsys@defobject{currentmarker}{\pgfqpoint{-0.041667in}{0.000000in}}{\pgfqpoint{0.000000in}{0.000000in}}{%
\pgfpathmoveto{\pgfqpoint{0.000000in}{0.000000in}}%
\pgfpathlineto{\pgfqpoint{-0.041667in}{0.000000in}}%
\pgfusepath{stroke,fill}%
}%
\begin{pgfscope}%
\pgfsys@transformshift{1.125000in}{1.791564in}%
\pgfsys@useobject{currentmarker}{}%
\end{pgfscope}%
\end{pgfscope}%
\begin{pgfscope}%
\pgfsetbuttcap%
\pgfsetroundjoin%
\definecolor{currentfill}{rgb}{0.150000,0.150000,0.150000}%
\pgfsetfillcolor{currentfill}%
\pgfsetlinewidth{0.401500pt}%
\definecolor{currentstroke}{rgb}{0.150000,0.150000,0.150000}%
\pgfsetstrokecolor{currentstroke}%
\pgfsetdash{}{0pt}%
\pgfsys@defobject{currentmarker}{\pgfqpoint{-0.041667in}{0.000000in}}{\pgfqpoint{0.000000in}{0.000000in}}{%
\pgfpathmoveto{\pgfqpoint{0.000000in}{0.000000in}}%
\pgfpathlineto{\pgfqpoint{-0.041667in}{0.000000in}}%
\pgfusepath{stroke,fill}%
}%
\begin{pgfscope}%
\pgfsys@transformshift{1.125000in}{2.019121in}%
\pgfsys@useobject{currentmarker}{}%
\end{pgfscope}%
\end{pgfscope}%
\begin{pgfscope}%
\pgfsetbuttcap%
\pgfsetroundjoin%
\definecolor{currentfill}{rgb}{0.150000,0.150000,0.150000}%
\pgfsetfillcolor{currentfill}%
\pgfsetlinewidth{0.401500pt}%
\definecolor{currentstroke}{rgb}{0.150000,0.150000,0.150000}%
\pgfsetstrokecolor{currentstroke}%
\pgfsetdash{}{0pt}%
\pgfsys@defobject{currentmarker}{\pgfqpoint{-0.041667in}{0.000000in}}{\pgfqpoint{0.000000in}{0.000000in}}{%
\pgfpathmoveto{\pgfqpoint{0.000000in}{0.000000in}}%
\pgfpathlineto{\pgfqpoint{-0.041667in}{0.000000in}}%
\pgfusepath{stroke,fill}%
}%
\begin{pgfscope}%
\pgfsys@transformshift{1.125000in}{2.205047in}%
\pgfsys@useobject{currentmarker}{}%
\end{pgfscope}%
\end{pgfscope}%
\begin{pgfscope}%
\pgfsetbuttcap%
\pgfsetroundjoin%
\definecolor{currentfill}{rgb}{0.150000,0.150000,0.150000}%
\pgfsetfillcolor{currentfill}%
\pgfsetlinewidth{0.401500pt}%
\definecolor{currentstroke}{rgb}{0.150000,0.150000,0.150000}%
\pgfsetstrokecolor{currentstroke}%
\pgfsetdash{}{0pt}%
\pgfsys@defobject{currentmarker}{\pgfqpoint{-0.041667in}{0.000000in}}{\pgfqpoint{0.000000in}{0.000000in}}{%
\pgfpathmoveto{\pgfqpoint{0.000000in}{0.000000in}}%
\pgfpathlineto{\pgfqpoint{-0.041667in}{0.000000in}}%
\pgfusepath{stroke,fill}%
}%
\begin{pgfscope}%
\pgfsys@transformshift{1.125000in}{2.362246in}%
\pgfsys@useobject{currentmarker}{}%
\end{pgfscope}%
\end{pgfscope}%
\begin{pgfscope}%
\pgfsetbuttcap%
\pgfsetroundjoin%
\definecolor{currentfill}{rgb}{0.150000,0.150000,0.150000}%
\pgfsetfillcolor{currentfill}%
\pgfsetlinewidth{0.401500pt}%
\definecolor{currentstroke}{rgb}{0.150000,0.150000,0.150000}%
\pgfsetstrokecolor{currentstroke}%
\pgfsetdash{}{0pt}%
\pgfsys@defobject{currentmarker}{\pgfqpoint{-0.041667in}{0.000000in}}{\pgfqpoint{0.000000in}{0.000000in}}{%
\pgfpathmoveto{\pgfqpoint{0.000000in}{0.000000in}}%
\pgfpathlineto{\pgfqpoint{-0.041667in}{0.000000in}}%
\pgfusepath{stroke,fill}%
}%
\begin{pgfscope}%
\pgfsys@transformshift{1.125000in}{2.498418in}%
\pgfsys@useobject{currentmarker}{}%
\end{pgfscope}%
\end{pgfscope}%
\begin{pgfscope}%
\pgfsetbuttcap%
\pgfsetroundjoin%
\definecolor{currentfill}{rgb}{0.150000,0.150000,0.150000}%
\pgfsetfillcolor{currentfill}%
\pgfsetlinewidth{0.401500pt}%
\definecolor{currentstroke}{rgb}{0.150000,0.150000,0.150000}%
\pgfsetstrokecolor{currentstroke}%
\pgfsetdash{}{0pt}%
\pgfsys@defobject{currentmarker}{\pgfqpoint{-0.041667in}{0.000000in}}{\pgfqpoint{0.000000in}{0.000000in}}{%
\pgfpathmoveto{\pgfqpoint{0.000000in}{0.000000in}}%
\pgfpathlineto{\pgfqpoint{-0.041667in}{0.000000in}}%
\pgfusepath{stroke,fill}%
}%
\begin{pgfscope}%
\pgfsys@transformshift{1.125000in}{2.618530in}%
\pgfsys@useobject{currentmarker}{}%
\end{pgfscope}%
\end{pgfscope}%
\begin{pgfscope}%
\pgfsetbuttcap%
\pgfsetroundjoin%
\definecolor{currentfill}{rgb}{0.150000,0.150000,0.150000}%
\pgfsetfillcolor{currentfill}%
\pgfsetlinewidth{0.401500pt}%
\definecolor{currentstroke}{rgb}{0.150000,0.150000,0.150000}%
\pgfsetstrokecolor{currentstroke}%
\pgfsetdash{}{0pt}%
\pgfsys@defobject{currentmarker}{\pgfqpoint{-0.041667in}{0.000000in}}{\pgfqpoint{0.000000in}{0.000000in}}{%
\pgfpathmoveto{\pgfqpoint{0.000000in}{0.000000in}}%
\pgfpathlineto{\pgfqpoint{-0.041667in}{0.000000in}}%
\pgfusepath{stroke,fill}%
}%
\begin{pgfscope}%
\pgfsys@transformshift{1.125000in}{3.432828in}%
\pgfsys@useobject{currentmarker}{}%
\end{pgfscope}%
\end{pgfscope}%
\begin{pgfscope}%
\pgfsetbuttcap%
\pgfsetroundjoin%
\definecolor{currentfill}{rgb}{0.150000,0.150000,0.150000}%
\pgfsetfillcolor{currentfill}%
\pgfsetlinewidth{0.401500pt}%
\definecolor{currentstroke}{rgb}{0.150000,0.150000,0.150000}%
\pgfsetstrokecolor{currentstroke}%
\pgfsetdash{}{0pt}%
\pgfsys@defobject{currentmarker}{\pgfqpoint{-0.041667in}{0.000000in}}{\pgfqpoint{0.000000in}{0.000000in}}{%
\pgfpathmoveto{\pgfqpoint{0.000000in}{0.000000in}}%
\pgfpathlineto{\pgfqpoint{-0.041667in}{0.000000in}}%
\pgfusepath{stroke,fill}%
}%
\begin{pgfscope}%
\pgfsys@transformshift{1.125000in}{3.846311in}%
\pgfsys@useobject{currentmarker}{}%
\end{pgfscope}%
\end{pgfscope}%
\begin{pgfscope}%
\pgfsetbuttcap%
\pgfsetroundjoin%
\definecolor{currentfill}{rgb}{0.150000,0.150000,0.150000}%
\pgfsetfillcolor{currentfill}%
\pgfsetlinewidth{0.401500pt}%
\definecolor{currentstroke}{rgb}{0.150000,0.150000,0.150000}%
\pgfsetstrokecolor{currentstroke}%
\pgfsetdash{}{0pt}%
\pgfsys@defobject{currentmarker}{\pgfqpoint{-0.041667in}{0.000000in}}{\pgfqpoint{0.000000in}{0.000000in}}{%
\pgfpathmoveto{\pgfqpoint{0.000000in}{0.000000in}}%
\pgfpathlineto{\pgfqpoint{-0.041667in}{0.000000in}}%
\pgfusepath{stroke,fill}%
}%
\begin{pgfscope}%
\pgfsys@transformshift{1.125000in}{4.139681in}%
\pgfsys@useobject{currentmarker}{}%
\end{pgfscope}%
\end{pgfscope}%
\begin{pgfscope}%
\pgfsetbuttcap%
\pgfsetroundjoin%
\definecolor{currentfill}{rgb}{0.150000,0.150000,0.150000}%
\pgfsetfillcolor{currentfill}%
\pgfsetlinewidth{0.401500pt}%
\definecolor{currentstroke}{rgb}{0.150000,0.150000,0.150000}%
\pgfsetstrokecolor{currentstroke}%
\pgfsetdash{}{0pt}%
\pgfsys@defobject{currentmarker}{\pgfqpoint{-0.041667in}{0.000000in}}{\pgfqpoint{0.000000in}{0.000000in}}{%
\pgfpathmoveto{\pgfqpoint{0.000000in}{0.000000in}}%
\pgfpathlineto{\pgfqpoint{-0.041667in}{0.000000in}}%
\pgfusepath{stroke,fill}%
}%
\begin{pgfscope}%
\pgfsys@transformshift{1.125000in}{4.367237in}%
\pgfsys@useobject{currentmarker}{}%
\end{pgfscope}%
\end{pgfscope}%
\begin{pgfscope}%
\pgfsetbuttcap%
\pgfsetroundjoin%
\definecolor{currentfill}{rgb}{0.150000,0.150000,0.150000}%
\pgfsetfillcolor{currentfill}%
\pgfsetlinewidth{0.401500pt}%
\definecolor{currentstroke}{rgb}{0.150000,0.150000,0.150000}%
\pgfsetstrokecolor{currentstroke}%
\pgfsetdash{}{0pt}%
\pgfsys@defobject{currentmarker}{\pgfqpoint{-0.041667in}{0.000000in}}{\pgfqpoint{0.000000in}{0.000000in}}{%
\pgfpathmoveto{\pgfqpoint{0.000000in}{0.000000in}}%
\pgfpathlineto{\pgfqpoint{-0.041667in}{0.000000in}}%
\pgfusepath{stroke,fill}%
}%
\begin{pgfscope}%
\pgfsys@transformshift{1.125000in}{4.553164in}%
\pgfsys@useobject{currentmarker}{}%
\end{pgfscope}%
\end{pgfscope}%
\begin{pgfscope}%
\pgfsetbuttcap%
\pgfsetroundjoin%
\definecolor{currentfill}{rgb}{0.150000,0.150000,0.150000}%
\pgfsetfillcolor{currentfill}%
\pgfsetlinewidth{0.401500pt}%
\definecolor{currentstroke}{rgb}{0.150000,0.150000,0.150000}%
\pgfsetstrokecolor{currentstroke}%
\pgfsetdash{}{0pt}%
\pgfsys@defobject{currentmarker}{\pgfqpoint{-0.041667in}{0.000000in}}{\pgfqpoint{0.000000in}{0.000000in}}{%
\pgfpathmoveto{\pgfqpoint{0.000000in}{0.000000in}}%
\pgfpathlineto{\pgfqpoint{-0.041667in}{0.000000in}}%
\pgfusepath{stroke,fill}%
}%
\begin{pgfscope}%
\pgfsys@transformshift{1.125000in}{4.710363in}%
\pgfsys@useobject{currentmarker}{}%
\end{pgfscope}%
\end{pgfscope}%
\begin{pgfscope}%
\pgfsetbuttcap%
\pgfsetroundjoin%
\definecolor{currentfill}{rgb}{0.150000,0.150000,0.150000}%
\pgfsetfillcolor{currentfill}%
\pgfsetlinewidth{0.401500pt}%
\definecolor{currentstroke}{rgb}{0.150000,0.150000,0.150000}%
\pgfsetstrokecolor{currentstroke}%
\pgfsetdash{}{0pt}%
\pgfsys@defobject{currentmarker}{\pgfqpoint{-0.041667in}{0.000000in}}{\pgfqpoint{0.000000in}{0.000000in}}{%
\pgfpathmoveto{\pgfqpoint{0.000000in}{0.000000in}}%
\pgfpathlineto{\pgfqpoint{-0.041667in}{0.000000in}}%
\pgfusepath{stroke,fill}%
}%
\begin{pgfscope}%
\pgfsys@transformshift{1.125000in}{4.846535in}%
\pgfsys@useobject{currentmarker}{}%
\end{pgfscope}%
\end{pgfscope}%
\begin{pgfscope}%
\pgfsetbuttcap%
\pgfsetroundjoin%
\definecolor{currentfill}{rgb}{0.150000,0.150000,0.150000}%
\pgfsetfillcolor{currentfill}%
\pgfsetlinewidth{0.401500pt}%
\definecolor{currentstroke}{rgb}{0.150000,0.150000,0.150000}%
\pgfsetstrokecolor{currentstroke}%
\pgfsetdash{}{0pt}%
\pgfsys@defobject{currentmarker}{\pgfqpoint{-0.041667in}{0.000000in}}{\pgfqpoint{0.000000in}{0.000000in}}{%
\pgfpathmoveto{\pgfqpoint{0.000000in}{0.000000in}}%
\pgfpathlineto{\pgfqpoint{-0.041667in}{0.000000in}}%
\pgfusepath{stroke,fill}%
}%
\begin{pgfscope}%
\pgfsys@transformshift{1.125000in}{4.966647in}%
\pgfsys@useobject{currentmarker}{}%
\end{pgfscope}%
\end{pgfscope}%
\begin{pgfscope}%
\definecolor{textcolor}{rgb}{0.150000,0.150000,0.150000}%
\pgfsetstrokecolor{textcolor}%
\pgfsetfillcolor{textcolor}%
\pgftext[x=0.383703in,y=3.015000in,,bottom,rotate=90.000000]{\color{textcolor}\sffamily\fontsize{24.000000}{28.800000}\selectfont Spectral radius, \(\displaystyle \rho\)}%
\end{pgfscope}%
\begin{pgfscope}%
\pgfpathrectangle{\pgfqpoint{1.125000in}{0.750000in}}{\pgfqpoint{6.975000in}{4.530000in}} %
\pgfusepath{clip}%
\pgfsetroundcap%
\pgfsetroundjoin%
\pgfsetlinewidth{1.405250pt}%
\definecolor{currentstroke}{rgb}{0.133333,0.133333,0.133333}%
\pgfsetstrokecolor{currentstroke}%
\pgfsetdash{}{0pt}%
\pgfpathmoveto{\pgfqpoint{1.442045in}{3.253660in}}%
\pgfpathlineto{\pgfqpoint{4.018037in}{2.207008in}}%
\pgfpathlineto{\pgfqpoint{5.092312in}{2.224842in}}%
\pgfpathlineto{\pgfqpoint{5.781711in}{2.229330in}}%
\pgfpathlineto{\pgfqpoint{6.290401in}{2.231074in}}%
\pgfpathlineto{\pgfqpoint{6.693708in}{2.231925in}}%
\pgfpathlineto{\pgfqpoint{7.027890in}{2.232402in}}%
\pgfpathlineto{\pgfqpoint{7.313212in}{2.232695in}}%
\pgfpathlineto{\pgfqpoint{7.562154in}{2.232889in}}%
\pgfpathlineto{\pgfqpoint{7.782955in}{2.233024in}}%
\pgfusepath{stroke}%
\end{pgfscope}%
\begin{pgfscope}%
\pgfpathrectangle{\pgfqpoint{1.125000in}{0.750000in}}{\pgfqpoint{6.975000in}{4.530000in}} %
\pgfusepath{clip}%
\pgfsetbuttcap%
\pgfsetmiterjoin%
\definecolor{currentfill}{rgb}{0.133333,0.133333,0.133333}%
\pgfsetfillcolor{currentfill}%
\pgfsetlinewidth{0.000000pt}%
\definecolor{currentstroke}{rgb}{0.133333,0.133333,0.133333}%
\pgfsetstrokecolor{currentstroke}%
\pgfsetdash{}{0pt}%
\pgfsys@defobject{currentmarker}{\pgfqpoint{-0.062500in}{-0.062500in}}{\pgfqpoint{0.062500in}{0.062500in}}{%
\pgfpathmoveto{\pgfqpoint{-0.000000in}{-0.062500in}}%
\pgfpathlineto{\pgfqpoint{0.062500in}{0.062500in}}%
\pgfpathlineto{\pgfqpoint{-0.062500in}{0.062500in}}%
\pgfpathclose%
\pgfusepath{fill}%
}%
\begin{pgfscope}%
\pgfsys@transformshift{1.442045in}{3.253660in}%
\pgfsys@useobject{currentmarker}{}%
\end{pgfscope}%
\begin{pgfscope}%
\pgfsys@transformshift{4.018037in}{2.207008in}%
\pgfsys@useobject{currentmarker}{}%
\end{pgfscope}%
\begin{pgfscope}%
\pgfsys@transformshift{5.092312in}{2.224842in}%
\pgfsys@useobject{currentmarker}{}%
\end{pgfscope}%
\begin{pgfscope}%
\pgfsys@transformshift{5.781711in}{2.229330in}%
\pgfsys@useobject{currentmarker}{}%
\end{pgfscope}%
\begin{pgfscope}%
\pgfsys@transformshift{6.290401in}{2.231074in}%
\pgfsys@useobject{currentmarker}{}%
\end{pgfscope}%
\begin{pgfscope}%
\pgfsys@transformshift{6.693708in}{2.231925in}%
\pgfsys@useobject{currentmarker}{}%
\end{pgfscope}%
\begin{pgfscope}%
\pgfsys@transformshift{7.027890in}{2.232402in}%
\pgfsys@useobject{currentmarker}{}%
\end{pgfscope}%
\begin{pgfscope}%
\pgfsys@transformshift{7.313212in}{2.232695in}%
\pgfsys@useobject{currentmarker}{}%
\end{pgfscope}%
\begin{pgfscope}%
\pgfsys@transformshift{7.562154in}{2.232889in}%
\pgfsys@useobject{currentmarker}{}%
\end{pgfscope}%
\begin{pgfscope}%
\pgfsys@transformshift{7.782955in}{2.233024in}%
\pgfsys@useobject{currentmarker}{}%
\end{pgfscope}%
\end{pgfscope}%
\begin{pgfscope}%
\pgfpathrectangle{\pgfqpoint{1.125000in}{0.750000in}}{\pgfqpoint{6.975000in}{4.530000in}} %
\pgfusepath{clip}%
\pgfsetroundcap%
\pgfsetroundjoin%
\pgfsetlinewidth{1.405250pt}%
\definecolor{currentstroke}{rgb}{0.318370,0.318370,0.318370}%
\pgfsetstrokecolor{currentstroke}%
\pgfsetdash{}{0pt}%
\pgfpathmoveto{\pgfqpoint{1.442045in}{0.955909in}}%
\pgfpathlineto{\pgfqpoint{4.018037in}{1.872075in}}%
\pgfpathlineto{\pgfqpoint{5.092312in}{1.863960in}}%
\pgfpathlineto{\pgfqpoint{5.781711in}{1.872908in}}%
\pgfpathlineto{\pgfqpoint{6.290401in}{1.896328in}}%
\pgfpathlineto{\pgfqpoint{6.693708in}{1.917193in}}%
\pgfpathlineto{\pgfqpoint{7.027890in}{1.969642in}}%
\pgfpathlineto{\pgfqpoint{7.313212in}{2.000273in}}%
\pgfpathlineto{\pgfqpoint{7.562154in}{2.050771in}}%
\pgfpathlineto{\pgfqpoint{7.782955in}{2.061944in}}%
\pgfusepath{stroke}%
\end{pgfscope}%
\begin{pgfscope}%
\pgfpathrectangle{\pgfqpoint{1.125000in}{0.750000in}}{\pgfqpoint{6.975000in}{4.530000in}} %
\pgfusepath{clip}%
\pgfsetbuttcap%
\pgfsetroundjoin%
\definecolor{currentfill}{rgb}{0.318370,0.318370,0.318370}%
\pgfsetfillcolor{currentfill}%
\pgfsetlinewidth{0.000000pt}%
\definecolor{currentstroke}{rgb}{0.318370,0.318370,0.318370}%
\pgfsetstrokecolor{currentstroke}%
\pgfsetdash{}{0pt}%
\pgfsys@defobject{currentmarker}{\pgfqpoint{-0.048611in}{-0.048611in}}{\pgfqpoint{0.048611in}{0.048611in}}{%
\pgfpathmoveto{\pgfqpoint{0.000000in}{-0.048611in}}%
\pgfpathcurveto{\pgfqpoint{0.012892in}{-0.048611in}}{\pgfqpoint{0.025257in}{-0.043489in}}{\pgfqpoint{0.034373in}{-0.034373in}}%
\pgfpathcurveto{\pgfqpoint{0.043489in}{-0.025257in}}{\pgfqpoint{0.048611in}{-0.012892in}}{\pgfqpoint{0.048611in}{0.000000in}}%
\pgfpathcurveto{\pgfqpoint{0.048611in}{0.012892in}}{\pgfqpoint{0.043489in}{0.025257in}}{\pgfqpoint{0.034373in}{0.034373in}}%
\pgfpathcurveto{\pgfqpoint{0.025257in}{0.043489in}}{\pgfqpoint{0.012892in}{0.048611in}}{\pgfqpoint{0.000000in}{0.048611in}}%
\pgfpathcurveto{\pgfqpoint{-0.012892in}{0.048611in}}{\pgfqpoint{-0.025257in}{0.043489in}}{\pgfqpoint{-0.034373in}{0.034373in}}%
\pgfpathcurveto{\pgfqpoint{-0.043489in}{0.025257in}}{\pgfqpoint{-0.048611in}{0.012892in}}{\pgfqpoint{-0.048611in}{0.000000in}}%
\pgfpathcurveto{\pgfqpoint{-0.048611in}{-0.012892in}}{\pgfqpoint{-0.043489in}{-0.025257in}}{\pgfqpoint{-0.034373in}{-0.034373in}}%
\pgfpathcurveto{\pgfqpoint{-0.025257in}{-0.043489in}}{\pgfqpoint{-0.012892in}{-0.048611in}}{\pgfqpoint{0.000000in}{-0.048611in}}%
\pgfpathclose%
\pgfusepath{fill}%
}%
\begin{pgfscope}%
\pgfsys@transformshift{1.442045in}{0.955909in}%
\pgfsys@useobject{currentmarker}{}%
\end{pgfscope}%
\begin{pgfscope}%
\pgfsys@transformshift{4.018037in}{1.872075in}%
\pgfsys@useobject{currentmarker}{}%
\end{pgfscope}%
\begin{pgfscope}%
\pgfsys@transformshift{5.092312in}{1.863960in}%
\pgfsys@useobject{currentmarker}{}%
\end{pgfscope}%
\begin{pgfscope}%
\pgfsys@transformshift{5.781711in}{1.872908in}%
\pgfsys@useobject{currentmarker}{}%
\end{pgfscope}%
\begin{pgfscope}%
\pgfsys@transformshift{6.290401in}{1.896328in}%
\pgfsys@useobject{currentmarker}{}%
\end{pgfscope}%
\begin{pgfscope}%
\pgfsys@transformshift{6.693708in}{1.917193in}%
\pgfsys@useobject{currentmarker}{}%
\end{pgfscope}%
\begin{pgfscope}%
\pgfsys@transformshift{7.027890in}{1.969642in}%
\pgfsys@useobject{currentmarker}{}%
\end{pgfscope}%
\begin{pgfscope}%
\pgfsys@transformshift{7.313212in}{2.000273in}%
\pgfsys@useobject{currentmarker}{}%
\end{pgfscope}%
\begin{pgfscope}%
\pgfsys@transformshift{7.562154in}{2.050771in}%
\pgfsys@useobject{currentmarker}{}%
\end{pgfscope}%
\begin{pgfscope}%
\pgfsys@transformshift{7.782955in}{2.061944in}%
\pgfsys@useobject{currentmarker}{}%
\end{pgfscope}%
\end{pgfscope}%
\begin{pgfscope}%
\pgfpathrectangle{\pgfqpoint{1.125000in}{0.750000in}}{\pgfqpoint{6.975000in}{4.530000in}} %
\pgfusepath{clip}%
\pgfsetroundcap%
\pgfsetroundjoin%
\pgfsetlinewidth{1.405250pt}%
\definecolor{currentstroke}{rgb}{0.501961,0.501961,0.501961}%
\pgfsetstrokecolor{currentstroke}%
\pgfsetdash{}{0pt}%
\pgfpathmoveto{\pgfqpoint{1.442045in}{3.992211in}}%
\pgfpathlineto{\pgfqpoint{4.018037in}{3.407403in}}%
\pgfpathlineto{\pgfqpoint{5.092312in}{3.417938in}}%
\pgfpathlineto{\pgfqpoint{5.781711in}{3.418386in}}%
\pgfpathlineto{\pgfqpoint{6.290401in}{3.417125in}}%
\pgfpathlineto{\pgfqpoint{6.693708in}{3.418395in}}%
\pgfpathlineto{\pgfqpoint{7.027890in}{3.414490in}}%
\pgfpathlineto{\pgfqpoint{7.313212in}{3.418785in}}%
\pgfpathlineto{\pgfqpoint{7.562154in}{3.416974in}}%
\pgfpathlineto{\pgfqpoint{7.782955in}{3.417466in}}%
\pgfusepath{stroke}%
\end{pgfscope}%
\begin{pgfscope}%
\pgfpathrectangle{\pgfqpoint{1.125000in}{0.750000in}}{\pgfqpoint{6.975000in}{4.530000in}} %
\pgfusepath{clip}%
\pgfsetbuttcap%
\pgfsetroundjoin%
\definecolor{currentfill}{rgb}{0.501961,0.501961,0.501961}%
\pgfsetfillcolor{currentfill}%
\pgfsetlinewidth{0.000000pt}%
\definecolor{currentstroke}{rgb}{0.501961,0.501961,0.501961}%
\pgfsetstrokecolor{currentstroke}%
\pgfsetdash{}{0pt}%
\pgfsys@defobject{currentmarker}{\pgfqpoint{-0.034722in}{-0.034722in}}{\pgfqpoint{0.034722in}{0.034722in}}{%
\pgfpathmoveto{\pgfqpoint{0.000000in}{-0.034722in}}%
\pgfpathcurveto{\pgfqpoint{0.009208in}{-0.034722in}}{\pgfqpoint{0.018041in}{-0.031064in}}{\pgfqpoint{0.024552in}{-0.024552in}}%
\pgfpathcurveto{\pgfqpoint{0.031064in}{-0.018041in}}{\pgfqpoint{0.034722in}{-0.009208in}}{\pgfqpoint{0.034722in}{0.000000in}}%
\pgfpathcurveto{\pgfqpoint{0.034722in}{0.009208in}}{\pgfqpoint{0.031064in}{0.018041in}}{\pgfqpoint{0.024552in}{0.024552in}}%
\pgfpathcurveto{\pgfqpoint{0.018041in}{0.031064in}}{\pgfqpoint{0.009208in}{0.034722in}}{\pgfqpoint{0.000000in}{0.034722in}}%
\pgfpathcurveto{\pgfqpoint{-0.009208in}{0.034722in}}{\pgfqpoint{-0.018041in}{0.031064in}}{\pgfqpoint{-0.024552in}{0.024552in}}%
\pgfpathcurveto{\pgfqpoint{-0.031064in}{0.018041in}}{\pgfqpoint{-0.034722in}{0.009208in}}{\pgfqpoint{-0.034722in}{0.000000in}}%
\pgfpathcurveto{\pgfqpoint{-0.034722in}{-0.009208in}}{\pgfqpoint{-0.031064in}{-0.018041in}}{\pgfqpoint{-0.024552in}{-0.024552in}}%
\pgfpathcurveto{\pgfqpoint{-0.018041in}{-0.031064in}}{\pgfqpoint{-0.009208in}{-0.034722in}}{\pgfqpoint{0.000000in}{-0.034722in}}%
\pgfpathclose%
\pgfusepath{fill}%
}%
\begin{pgfscope}%
\pgfsys@transformshift{1.442045in}{3.992211in}%
\pgfsys@useobject{currentmarker}{}%
\end{pgfscope}%
\begin{pgfscope}%
\pgfsys@transformshift{4.018037in}{3.407403in}%
\pgfsys@useobject{currentmarker}{}%
\end{pgfscope}%
\begin{pgfscope}%
\pgfsys@transformshift{5.092312in}{3.417938in}%
\pgfsys@useobject{currentmarker}{}%
\end{pgfscope}%
\begin{pgfscope}%
\pgfsys@transformshift{5.781711in}{3.418386in}%
\pgfsys@useobject{currentmarker}{}%
\end{pgfscope}%
\begin{pgfscope}%
\pgfsys@transformshift{6.290401in}{3.417125in}%
\pgfsys@useobject{currentmarker}{}%
\end{pgfscope}%
\begin{pgfscope}%
\pgfsys@transformshift{6.693708in}{3.418395in}%
\pgfsys@useobject{currentmarker}{}%
\end{pgfscope}%
\begin{pgfscope}%
\pgfsys@transformshift{7.027890in}{3.414490in}%
\pgfsys@useobject{currentmarker}{}%
\end{pgfscope}%
\begin{pgfscope}%
\pgfsys@transformshift{7.313212in}{3.418785in}%
\pgfsys@useobject{currentmarker}{}%
\end{pgfscope}%
\begin{pgfscope}%
\pgfsys@transformshift{7.562154in}{3.416974in}%
\pgfsys@useobject{currentmarker}{}%
\end{pgfscope}%
\begin{pgfscope}%
\pgfsys@transformshift{7.782955in}{3.417466in}%
\pgfsys@useobject{currentmarker}{}%
\end{pgfscope}%
\end{pgfscope}%
\begin{pgfscope}%
\pgfpathrectangle{\pgfqpoint{1.125000in}{0.750000in}}{\pgfqpoint{6.975000in}{4.530000in}} %
\pgfusepath{clip}%
\pgfsetbuttcap%
\pgfsetroundjoin%
\pgfsetlinewidth{1.405250pt}%
\definecolor{currentstroke}{rgb}{0.000000,0.000000,0.000000}%
\pgfsetstrokecolor{currentstroke}%
\pgfsetdash{{5.180000pt}{2.240000pt}}{0.000000pt}%
\pgfpathmoveto{\pgfqpoint{1.442045in}{5.074091in}}%
\pgfpathlineto{\pgfqpoint{4.018037in}{5.074091in}}%
\pgfpathlineto{\pgfqpoint{5.092312in}{5.074091in}}%
\pgfpathlineto{\pgfqpoint{5.781711in}{5.074091in}}%
\pgfpathlineto{\pgfqpoint{6.290401in}{5.074091in}}%
\pgfpathlineto{\pgfqpoint{6.693708in}{5.074091in}}%
\pgfpathlineto{\pgfqpoint{7.027890in}{5.074091in}}%
\pgfpathlineto{\pgfqpoint{7.313212in}{5.074091in}}%
\pgfpathlineto{\pgfqpoint{7.562154in}{5.074091in}}%
\pgfpathlineto{\pgfqpoint{7.782955in}{5.074091in}}%
\pgfusepath{stroke}%
\end{pgfscope}%
\begin{pgfscope}%
\pgfsetrectcap%
\pgfsetmiterjoin%
\pgfsetlinewidth{1.254687pt}%
\definecolor{currentstroke}{rgb}{0.150000,0.150000,0.150000}%
\pgfsetstrokecolor{currentstroke}%
\pgfsetdash{}{0pt}%
\pgfpathmoveto{\pgfqpoint{1.125000in}{0.750000in}}%
\pgfpathlineto{\pgfqpoint{1.125000in}{5.280000in}}%
\pgfusepath{stroke}%
\end{pgfscope}%
\begin{pgfscope}%
\pgfsetrectcap%
\pgfsetmiterjoin%
\pgfsetlinewidth{1.254687pt}%
\definecolor{currentstroke}{rgb}{0.150000,0.150000,0.150000}%
\pgfsetstrokecolor{currentstroke}%
\pgfsetdash{}{0pt}%
\pgfpathmoveto{\pgfqpoint{8.100000in}{0.750000in}}%
\pgfpathlineto{\pgfqpoint{8.100000in}{5.280000in}}%
\pgfusepath{stroke}%
\end{pgfscope}%
\begin{pgfscope}%
\pgfsetrectcap%
\pgfsetmiterjoin%
\pgfsetlinewidth{1.254687pt}%
\definecolor{currentstroke}{rgb}{0.150000,0.150000,0.150000}%
\pgfsetstrokecolor{currentstroke}%
\pgfsetdash{}{0pt}%
\pgfpathmoveto{\pgfqpoint{1.125000in}{0.750000in}}%
\pgfpathlineto{\pgfqpoint{8.100000in}{0.750000in}}%
\pgfusepath{stroke}%
\end{pgfscope}%
\begin{pgfscope}%
\pgfsetrectcap%
\pgfsetmiterjoin%
\pgfsetlinewidth{1.254687pt}%
\definecolor{currentstroke}{rgb}{0.150000,0.150000,0.150000}%
\pgfsetstrokecolor{currentstroke}%
\pgfsetdash{}{0pt}%
\pgfpathmoveto{\pgfqpoint{1.125000in}{5.280000in}}%
\pgfpathlineto{\pgfqpoint{8.100000in}{5.280000in}}%
\pgfusepath{stroke}%
\end{pgfscope}%
\begin{pgfscope}%
\pgfsetroundcap%
\pgfsetroundjoin%
\pgfsetlinewidth{1.405250pt}%
\definecolor{currentstroke}{rgb}{0.133333,0.133333,0.133333}%
\pgfsetstrokecolor{currentstroke}%
\pgfsetdash{}{0pt}%
\pgfpathmoveto{\pgfqpoint{1.312500in}{2.809621in}}%
\pgfpathlineto{\pgfqpoint{1.729167in}{2.809621in}}%
\pgfusepath{stroke}%
\end{pgfscope}%
\begin{pgfscope}%
\pgfsetbuttcap%
\pgfsetmiterjoin%
\definecolor{currentfill}{rgb}{0.133333,0.133333,0.133333}%
\pgfsetfillcolor{currentfill}%
\pgfsetlinewidth{0.000000pt}%
\definecolor{currentstroke}{rgb}{0.133333,0.133333,0.133333}%
\pgfsetstrokecolor{currentstroke}%
\pgfsetdash{}{0pt}%
\pgfsys@defobject{currentmarker}{\pgfqpoint{-0.062500in}{-0.062500in}}{\pgfqpoint{0.062500in}{0.062500in}}{%
\pgfpathmoveto{\pgfqpoint{-0.000000in}{-0.062500in}}%
\pgfpathlineto{\pgfqpoint{0.062500in}{0.062500in}}%
\pgfpathlineto{\pgfqpoint{-0.062500in}{0.062500in}}%
\pgfpathclose%
\pgfusepath{fill}%
}%
\begin{pgfscope}%
\pgfsys@transformshift{1.520833in}{2.809621in}%
\pgfsys@useobject{currentmarker}{}%
\end{pgfscope}%
\end{pgfscope}%
\begin{pgfscope}%
\definecolor{textcolor}{rgb}{0.150000,0.150000,0.150000}%
\pgfsetstrokecolor{textcolor}%
\pgfsetfillcolor{textcolor}%
\pgftext[x=1.895833in,y=2.736704in,left,base]{\color{textcolor}\sffamily\fontsize{15.000000}{18.000000}\selectfont Linear}%
\end{pgfscope}%
\begin{pgfscope}%
\pgfsetroundcap%
\pgfsetroundjoin%
\pgfsetlinewidth{1.405250pt}%
\definecolor{currentstroke}{rgb}{0.318370,0.318370,0.318370}%
\pgfsetstrokecolor{currentstroke}%
\pgfsetdash{}{0pt}%
\pgfpathmoveto{\pgfqpoint{1.312500in}{2.503835in}}%
\pgfpathlineto{\pgfqpoint{1.729167in}{2.503835in}}%
\pgfusepath{stroke}%
\end{pgfscope}%
\begin{pgfscope}%
\pgfsetbuttcap%
\pgfsetroundjoin%
\definecolor{currentfill}{rgb}{0.318370,0.318370,0.318370}%
\pgfsetfillcolor{currentfill}%
\pgfsetlinewidth{0.000000pt}%
\definecolor{currentstroke}{rgb}{0.318370,0.318370,0.318370}%
\pgfsetstrokecolor{currentstroke}%
\pgfsetdash{}{0pt}%
\pgfsys@defobject{currentmarker}{\pgfqpoint{-0.048611in}{-0.048611in}}{\pgfqpoint{0.048611in}{0.048611in}}{%
\pgfpathmoveto{\pgfqpoint{0.000000in}{-0.048611in}}%
\pgfpathcurveto{\pgfqpoint{0.012892in}{-0.048611in}}{\pgfqpoint{0.025257in}{-0.043489in}}{\pgfqpoint{0.034373in}{-0.034373in}}%
\pgfpathcurveto{\pgfqpoint{0.043489in}{-0.025257in}}{\pgfqpoint{0.048611in}{-0.012892in}}{\pgfqpoint{0.048611in}{0.000000in}}%
\pgfpathcurveto{\pgfqpoint{0.048611in}{0.012892in}}{\pgfqpoint{0.043489in}{0.025257in}}{\pgfqpoint{0.034373in}{0.034373in}}%
\pgfpathcurveto{\pgfqpoint{0.025257in}{0.043489in}}{\pgfqpoint{0.012892in}{0.048611in}}{\pgfqpoint{0.000000in}{0.048611in}}%
\pgfpathcurveto{\pgfqpoint{-0.012892in}{0.048611in}}{\pgfqpoint{-0.025257in}{0.043489in}}{\pgfqpoint{-0.034373in}{0.034373in}}%
\pgfpathcurveto{\pgfqpoint{-0.043489in}{0.025257in}}{\pgfqpoint{-0.048611in}{0.012892in}}{\pgfqpoint{-0.048611in}{0.000000in}}%
\pgfpathcurveto{\pgfqpoint{-0.048611in}{-0.012892in}}{\pgfqpoint{-0.043489in}{-0.025257in}}{\pgfqpoint{-0.034373in}{-0.034373in}}%
\pgfpathcurveto{\pgfqpoint{-0.025257in}{-0.043489in}}{\pgfqpoint{-0.012892in}{-0.048611in}}{\pgfqpoint{0.000000in}{-0.048611in}}%
\pgfpathclose%
\pgfusepath{fill}%
}%
\begin{pgfscope}%
\pgfsys@transformshift{1.520833in}{2.503835in}%
\pgfsys@useobject{currentmarker}{}%
\end{pgfscope}%
\end{pgfscope}%
\begin{pgfscope}%
\definecolor{textcolor}{rgb}{0.150000,0.150000,0.150000}%
\pgfsetstrokecolor{textcolor}%
\pgfsetfillcolor{textcolor}%
\pgftext[x=1.895833in,y=2.430918in,left,base]{\color{textcolor}\sffamily\fontsize{15.000000}{18.000000}\selectfont DMG}%
\end{pgfscope}%
\begin{pgfscope}%
\pgfsetroundcap%
\pgfsetroundjoin%
\pgfsetlinewidth{1.405250pt}%
\definecolor{currentstroke}{rgb}{0.501961,0.501961,0.501961}%
\pgfsetstrokecolor{currentstroke}%
\pgfsetdash{}{0pt}%
\pgfpathmoveto{\pgfqpoint{1.312500in}{2.198049in}}%
\pgfpathlineto{\pgfqpoint{1.729167in}{2.198049in}}%
\pgfusepath{stroke}%
\end{pgfscope}%
\begin{pgfscope}%
\pgfsetbuttcap%
\pgfsetroundjoin%
\definecolor{currentfill}{rgb}{0.501961,0.501961,0.501961}%
\pgfsetfillcolor{currentfill}%
\pgfsetlinewidth{0.000000pt}%
\definecolor{currentstroke}{rgb}{0.501961,0.501961,0.501961}%
\pgfsetstrokecolor{currentstroke}%
\pgfsetdash{}{0pt}%
\pgfsys@defobject{currentmarker}{\pgfqpoint{-0.034722in}{-0.034722in}}{\pgfqpoint{0.034722in}{0.034722in}}{%
\pgfpathmoveto{\pgfqpoint{0.000000in}{-0.034722in}}%
\pgfpathcurveto{\pgfqpoint{0.009208in}{-0.034722in}}{\pgfqpoint{0.018041in}{-0.031064in}}{\pgfqpoint{0.024552in}{-0.024552in}}%
\pgfpathcurveto{\pgfqpoint{0.031064in}{-0.018041in}}{\pgfqpoint{0.034722in}{-0.009208in}}{\pgfqpoint{0.034722in}{0.000000in}}%
\pgfpathcurveto{\pgfqpoint{0.034722in}{0.009208in}}{\pgfqpoint{0.031064in}{0.018041in}}{\pgfqpoint{0.024552in}{0.024552in}}%
\pgfpathcurveto{\pgfqpoint{0.018041in}{0.031064in}}{\pgfqpoint{0.009208in}{0.034722in}}{\pgfqpoint{0.000000in}{0.034722in}}%
\pgfpathcurveto{\pgfqpoint{-0.009208in}{0.034722in}}{\pgfqpoint{-0.018041in}{0.031064in}}{\pgfqpoint{-0.024552in}{0.024552in}}%
\pgfpathcurveto{\pgfqpoint{-0.031064in}{0.018041in}}{\pgfqpoint{-0.034722in}{0.009208in}}{\pgfqpoint{-0.034722in}{0.000000in}}%
\pgfpathcurveto{\pgfqpoint{-0.034722in}{-0.009208in}}{\pgfqpoint{-0.031064in}{-0.018041in}}{\pgfqpoint{-0.024552in}{-0.024552in}}%
\pgfpathcurveto{\pgfqpoint{-0.018041in}{-0.031064in}}{\pgfqpoint{-0.009208in}{-0.034722in}}{\pgfqpoint{0.000000in}{-0.034722in}}%
\pgfpathclose%
\pgfusepath{fill}%
}%
\begin{pgfscope}%
\pgfsys@transformshift{1.520833in}{2.198049in}%
\pgfsys@useobject{currentmarker}{}%
\end{pgfscope}%
\end{pgfscope}%
\begin{pgfscope}%
\definecolor{textcolor}{rgb}{0.150000,0.150000,0.150000}%
\pgfsetstrokecolor{textcolor}%
\pgfsetfillcolor{textcolor}%
\pgftext[x=1.895833in,y=2.125132in,left,base]{\color{textcolor}\sffamily\fontsize{15.000000}{18.000000}\selectfont AMG}%
\end{pgfscope}%
\begin{pgfscope}%
\pgfsetbuttcap%
\pgfsetroundjoin%
\pgfsetlinewidth{1.405250pt}%
\definecolor{currentstroke}{rgb}{0.000000,0.000000,0.000000}%
\pgfsetstrokecolor{currentstroke}%
\pgfsetdash{{5.180000pt}{2.240000pt}}{0.000000pt}%
\pgfpathmoveto{\pgfqpoint{1.312500in}{1.892263in}}%
\pgfpathlineto{\pgfqpoint{1.729167in}{1.892263in}}%
\pgfusepath{stroke}%
\end{pgfscope}%
\begin{pgfscope}%
\definecolor{textcolor}{rgb}{0.150000,0.150000,0.150000}%
\pgfsetstrokecolor{textcolor}%
\pgfsetfillcolor{textcolor}%
\pgftext[x=1.895833in,y=1.819346in,left,base]{\color{textcolor}\sffamily\fontsize{15.000000}{18.000000}\selectfont Boundary \(\displaystyle \rho = 1\)}%
\end{pgfscope}%
\end{pgfpicture}%
\makeatother%
\endgroup%

%% file: diff_m3points_epsrange_n_511_bw.pgf
\begingroup%
\makeatletter%
\begin{pgfpicture}%
\pgfpathrectangle{\pgfpointorigin}{\pgfqpoint{9.000000in}{6.000000in}}%
\pgfusepath{use as bounding box, clip}%
\begin{pgfscope}%
\pgfsetbuttcap%
\pgfsetmiterjoin%
\definecolor{currentfill}{rgb}{1.000000,1.000000,1.000000}%
\pgfsetfillcolor{currentfill}%
\pgfsetlinewidth{0.000000pt}%
\definecolor{currentstroke}{rgb}{1.000000,1.000000,1.000000}%
\pgfsetstrokecolor{currentstroke}%
\pgfsetdash{}{0pt}%
\pgfpathmoveto{\pgfqpoint{0.000000in}{0.000000in}}%
\pgfpathlineto{\pgfqpoint{9.000000in}{0.000000in}}%
\pgfpathlineto{\pgfqpoint{9.000000in}{6.000000in}}%
\pgfpathlineto{\pgfqpoint{0.000000in}{6.000000in}}%
\pgfpathclose%
\pgfusepath{fill}%
\end{pgfscope}%
\begin{pgfscope}%
\pgfsetbuttcap%
\pgfsetmiterjoin%
\definecolor{currentfill}{rgb}{1.000000,1.000000,1.000000}%
\pgfsetfillcolor{currentfill}%
\pgfsetlinewidth{0.000000pt}%
\definecolor{currentstroke}{rgb}{0.000000,0.000000,0.000000}%
\pgfsetstrokecolor{currentstroke}%
\pgfsetstrokeopacity{0.000000}%
\pgfsetdash{}{0pt}%
\pgfpathmoveto{\pgfqpoint{1.125000in}{0.750000in}}%
\pgfpathlineto{\pgfqpoint{8.100000in}{0.750000in}}%
\pgfpathlineto{\pgfqpoint{8.100000in}{5.280000in}}%
\pgfpathlineto{\pgfqpoint{1.125000in}{5.280000in}}%
\pgfpathclose%
\pgfusepath{fill}%
\end{pgfscope}%
\begin{pgfscope}%
\pgfsetbuttcap%
\pgfsetroundjoin%
\definecolor{currentfill}{rgb}{0.150000,0.150000,0.150000}%
\pgfsetfillcolor{currentfill}%
\pgfsetlinewidth{0.803000pt}%
\definecolor{currentstroke}{rgb}{0.150000,0.150000,0.150000}%
\pgfsetstrokecolor{currentstroke}%
\pgfsetdash{}{0pt}%
\pgfsys@defobject{currentmarker}{\pgfqpoint{0.000000in}{-0.083333in}}{\pgfqpoint{0.000000in}{0.000000in}}{%
\pgfpathmoveto{\pgfqpoint{0.000000in}{0.000000in}}%
\pgfpathlineto{\pgfqpoint{0.000000in}{-0.083333in}}%
\pgfusepath{stroke,fill}%
}%
\begin{pgfscope}%
\pgfsys@transformshift{4.073237in}{0.750000in}%
\pgfsys@useobject{currentmarker}{}%
\end{pgfscope}%
\end{pgfscope}%
\begin{pgfscope}%
\definecolor{textcolor}{rgb}{0.150000,0.150000,0.150000}%
\pgfsetstrokecolor{textcolor}%
\pgfsetfillcolor{textcolor}%
\pgftext[x=4.073237in,y=0.588889in,,top]{\color{textcolor}\sffamily\fontsize{16.000000}{19.200000}\selectfont \(\displaystyle 10^{-2}\)}%
\end{pgfscope}%
\begin{pgfscope}%
\pgfsetbuttcap%
\pgfsetroundjoin%
\definecolor{currentfill}{rgb}{0.150000,0.150000,0.150000}%
\pgfsetfillcolor{currentfill}%
\pgfsetlinewidth{0.803000pt}%
\definecolor{currentstroke}{rgb}{0.150000,0.150000,0.150000}%
\pgfsetstrokecolor{currentstroke}%
\pgfsetdash{}{0pt}%
\pgfsys@defobject{currentmarker}{\pgfqpoint{0.000000in}{-0.083333in}}{\pgfqpoint{0.000000in}{0.000000in}}{%
\pgfpathmoveto{\pgfqpoint{0.000000in}{0.000000in}}%
\pgfpathlineto{\pgfqpoint{0.000000in}{-0.083333in}}%
\pgfusepath{stroke,fill}%
}%
\begin{pgfscope}%
\pgfsys@transformshift{7.782955in}{0.750000in}%
\pgfsys@useobject{currentmarker}{}%
\end{pgfscope}%
\end{pgfscope}%
\begin{pgfscope}%
\definecolor{textcolor}{rgb}{0.150000,0.150000,0.150000}%
\pgfsetstrokecolor{textcolor}%
\pgfsetfillcolor{textcolor}%
\pgftext[x=7.782955in,y=0.588889in,,top]{\color{textcolor}\sffamily\fontsize{16.000000}{19.200000}\selectfont \(\displaystyle 10^{-1}\)}%
\end{pgfscope}%
\begin{pgfscope}%
\pgfsetbuttcap%
\pgfsetroundjoin%
\definecolor{currentfill}{rgb}{0.150000,0.150000,0.150000}%
\pgfsetfillcolor{currentfill}%
\pgfsetlinewidth{0.401500pt}%
\definecolor{currentstroke}{rgb}{0.150000,0.150000,0.150000}%
\pgfsetstrokecolor{currentstroke}%
\pgfsetdash{}{0pt}%
\pgfsys@defobject{currentmarker}{\pgfqpoint{0.000000in}{-0.041667in}}{\pgfqpoint{0.000000in}{0.000000in}}{%
\pgfpathmoveto{\pgfqpoint{0.000000in}{0.000000in}}%
\pgfpathlineto{\pgfqpoint{0.000000in}{-0.041667in}}%
\pgfusepath{stroke,fill}%
}%
\begin{pgfscope}%
\pgfsys@transformshift{1.480255in}{0.750000in}%
\pgfsys@useobject{currentmarker}{}%
\end{pgfscope}%
\end{pgfscope}%
\begin{pgfscope}%
\pgfsetbuttcap%
\pgfsetroundjoin%
\definecolor{currentfill}{rgb}{0.150000,0.150000,0.150000}%
\pgfsetfillcolor{currentfill}%
\pgfsetlinewidth{0.401500pt}%
\definecolor{currentstroke}{rgb}{0.150000,0.150000,0.150000}%
\pgfsetstrokecolor{currentstroke}%
\pgfsetdash{}{0pt}%
\pgfsys@defobject{currentmarker}{\pgfqpoint{0.000000in}{-0.041667in}}{\pgfqpoint{0.000000in}{0.000000in}}{%
\pgfpathmoveto{\pgfqpoint{0.000000in}{0.000000in}}%
\pgfpathlineto{\pgfqpoint{0.000000in}{-0.041667in}}%
\pgfusepath{stroke,fill}%
}%
\begin{pgfscope}%
\pgfsys@transformshift{2.133504in}{0.750000in}%
\pgfsys@useobject{currentmarker}{}%
\end{pgfscope}%
\end{pgfscope}%
\begin{pgfscope}%
\pgfsetbuttcap%
\pgfsetroundjoin%
\definecolor{currentfill}{rgb}{0.150000,0.150000,0.150000}%
\pgfsetfillcolor{currentfill}%
\pgfsetlinewidth{0.401500pt}%
\definecolor{currentstroke}{rgb}{0.150000,0.150000,0.150000}%
\pgfsetstrokecolor{currentstroke}%
\pgfsetdash{}{0pt}%
\pgfsys@defobject{currentmarker}{\pgfqpoint{0.000000in}{-0.041667in}}{\pgfqpoint{0.000000in}{0.000000in}}{%
\pgfpathmoveto{\pgfqpoint{0.000000in}{0.000000in}}%
\pgfpathlineto{\pgfqpoint{0.000000in}{-0.041667in}}%
\pgfusepath{stroke,fill}%
}%
\begin{pgfscope}%
\pgfsys@transformshift{2.596992in}{0.750000in}%
\pgfsys@useobject{currentmarker}{}%
\end{pgfscope}%
\end{pgfscope}%
\begin{pgfscope}%
\pgfsetbuttcap%
\pgfsetroundjoin%
\definecolor{currentfill}{rgb}{0.150000,0.150000,0.150000}%
\pgfsetfillcolor{currentfill}%
\pgfsetlinewidth{0.401500pt}%
\definecolor{currentstroke}{rgb}{0.150000,0.150000,0.150000}%
\pgfsetstrokecolor{currentstroke}%
\pgfsetdash{}{0pt}%
\pgfsys@defobject{currentmarker}{\pgfqpoint{0.000000in}{-0.041667in}}{\pgfqpoint{0.000000in}{0.000000in}}{%
\pgfpathmoveto{\pgfqpoint{0.000000in}{0.000000in}}%
\pgfpathlineto{\pgfqpoint{0.000000in}{-0.041667in}}%
\pgfusepath{stroke,fill}%
}%
\begin{pgfscope}%
\pgfsys@transformshift{2.956500in}{0.750000in}%
\pgfsys@useobject{currentmarker}{}%
\end{pgfscope}%
\end{pgfscope}%
\begin{pgfscope}%
\pgfsetbuttcap%
\pgfsetroundjoin%
\definecolor{currentfill}{rgb}{0.150000,0.150000,0.150000}%
\pgfsetfillcolor{currentfill}%
\pgfsetlinewidth{0.401500pt}%
\definecolor{currentstroke}{rgb}{0.150000,0.150000,0.150000}%
\pgfsetstrokecolor{currentstroke}%
\pgfsetdash{}{0pt}%
\pgfsys@defobject{currentmarker}{\pgfqpoint{0.000000in}{-0.041667in}}{\pgfqpoint{0.000000in}{0.000000in}}{%
\pgfpathmoveto{\pgfqpoint{0.000000in}{0.000000in}}%
\pgfpathlineto{\pgfqpoint{0.000000in}{-0.041667in}}%
\pgfusepath{stroke,fill}%
}%
\begin{pgfscope}%
\pgfsys@transformshift{3.250241in}{0.750000in}%
\pgfsys@useobject{currentmarker}{}%
\end{pgfscope}%
\end{pgfscope}%
\begin{pgfscope}%
\pgfsetbuttcap%
\pgfsetroundjoin%
\definecolor{currentfill}{rgb}{0.150000,0.150000,0.150000}%
\pgfsetfillcolor{currentfill}%
\pgfsetlinewidth{0.401500pt}%
\definecolor{currentstroke}{rgb}{0.150000,0.150000,0.150000}%
\pgfsetstrokecolor{currentstroke}%
\pgfsetdash{}{0pt}%
\pgfsys@defobject{currentmarker}{\pgfqpoint{0.000000in}{-0.041667in}}{\pgfqpoint{0.000000in}{0.000000in}}{%
\pgfpathmoveto{\pgfqpoint{0.000000in}{0.000000in}}%
\pgfpathlineto{\pgfqpoint{0.000000in}{-0.041667in}}%
\pgfusepath{stroke,fill}%
}%
\begin{pgfscope}%
\pgfsys@transformshift{3.498594in}{0.750000in}%
\pgfsys@useobject{currentmarker}{}%
\end{pgfscope}%
\end{pgfscope}%
\begin{pgfscope}%
\pgfsetbuttcap%
\pgfsetroundjoin%
\definecolor{currentfill}{rgb}{0.150000,0.150000,0.150000}%
\pgfsetfillcolor{currentfill}%
\pgfsetlinewidth{0.401500pt}%
\definecolor{currentstroke}{rgb}{0.150000,0.150000,0.150000}%
\pgfsetstrokecolor{currentstroke}%
\pgfsetdash{}{0pt}%
\pgfsys@defobject{currentmarker}{\pgfqpoint{0.000000in}{-0.041667in}}{\pgfqpoint{0.000000in}{0.000000in}}{%
\pgfpathmoveto{\pgfqpoint{0.000000in}{0.000000in}}%
\pgfpathlineto{\pgfqpoint{0.000000in}{-0.041667in}}%
\pgfusepath{stroke,fill}%
}%
\begin{pgfscope}%
\pgfsys@transformshift{3.713728in}{0.750000in}%
\pgfsys@useobject{currentmarker}{}%
\end{pgfscope}%
\end{pgfscope}%
\begin{pgfscope}%
\pgfsetbuttcap%
\pgfsetroundjoin%
\definecolor{currentfill}{rgb}{0.150000,0.150000,0.150000}%
\pgfsetfillcolor{currentfill}%
\pgfsetlinewidth{0.401500pt}%
\definecolor{currentstroke}{rgb}{0.150000,0.150000,0.150000}%
\pgfsetstrokecolor{currentstroke}%
\pgfsetdash{}{0pt}%
\pgfsys@defobject{currentmarker}{\pgfqpoint{0.000000in}{-0.041667in}}{\pgfqpoint{0.000000in}{0.000000in}}{%
\pgfpathmoveto{\pgfqpoint{0.000000in}{0.000000in}}%
\pgfpathlineto{\pgfqpoint{0.000000in}{-0.041667in}}%
\pgfusepath{stroke,fill}%
}%
\begin{pgfscope}%
\pgfsys@transformshift{3.903489in}{0.750000in}%
\pgfsys@useobject{currentmarker}{}%
\end{pgfscope}%
\end{pgfscope}%
\begin{pgfscope}%
\pgfsetbuttcap%
\pgfsetroundjoin%
\definecolor{currentfill}{rgb}{0.150000,0.150000,0.150000}%
\pgfsetfillcolor{currentfill}%
\pgfsetlinewidth{0.401500pt}%
\definecolor{currentstroke}{rgb}{0.150000,0.150000,0.150000}%
\pgfsetstrokecolor{currentstroke}%
\pgfsetdash{}{0pt}%
\pgfsys@defobject{currentmarker}{\pgfqpoint{0.000000in}{-0.041667in}}{\pgfqpoint{0.000000in}{0.000000in}}{%
\pgfpathmoveto{\pgfqpoint{0.000000in}{0.000000in}}%
\pgfpathlineto{\pgfqpoint{0.000000in}{-0.041667in}}%
\pgfusepath{stroke,fill}%
}%
\begin{pgfscope}%
\pgfsys@transformshift{5.189973in}{0.750000in}%
\pgfsys@useobject{currentmarker}{}%
\end{pgfscope}%
\end{pgfscope}%
\begin{pgfscope}%
\pgfsetbuttcap%
\pgfsetroundjoin%
\definecolor{currentfill}{rgb}{0.150000,0.150000,0.150000}%
\pgfsetfillcolor{currentfill}%
\pgfsetlinewidth{0.401500pt}%
\definecolor{currentstroke}{rgb}{0.150000,0.150000,0.150000}%
\pgfsetstrokecolor{currentstroke}%
\pgfsetdash{}{0pt}%
\pgfsys@defobject{currentmarker}{\pgfqpoint{0.000000in}{-0.041667in}}{\pgfqpoint{0.000000in}{0.000000in}}{%
\pgfpathmoveto{\pgfqpoint{0.000000in}{0.000000in}}%
\pgfpathlineto{\pgfqpoint{0.000000in}{-0.041667in}}%
\pgfusepath{stroke,fill}%
}%
\begin{pgfscope}%
\pgfsys@transformshift{5.843222in}{0.750000in}%
\pgfsys@useobject{currentmarker}{}%
\end{pgfscope}%
\end{pgfscope}%
\begin{pgfscope}%
\pgfsetbuttcap%
\pgfsetroundjoin%
\definecolor{currentfill}{rgb}{0.150000,0.150000,0.150000}%
\pgfsetfillcolor{currentfill}%
\pgfsetlinewidth{0.401500pt}%
\definecolor{currentstroke}{rgb}{0.150000,0.150000,0.150000}%
\pgfsetstrokecolor{currentstroke}%
\pgfsetdash{}{0pt}%
\pgfsys@defobject{currentmarker}{\pgfqpoint{0.000000in}{-0.041667in}}{\pgfqpoint{0.000000in}{0.000000in}}{%
\pgfpathmoveto{\pgfqpoint{0.000000in}{0.000000in}}%
\pgfpathlineto{\pgfqpoint{0.000000in}{-0.041667in}}%
\pgfusepath{stroke,fill}%
}%
\begin{pgfscope}%
\pgfsys@transformshift{6.306709in}{0.750000in}%
\pgfsys@useobject{currentmarker}{}%
\end{pgfscope}%
\end{pgfscope}%
\begin{pgfscope}%
\pgfsetbuttcap%
\pgfsetroundjoin%
\definecolor{currentfill}{rgb}{0.150000,0.150000,0.150000}%
\pgfsetfillcolor{currentfill}%
\pgfsetlinewidth{0.401500pt}%
\definecolor{currentstroke}{rgb}{0.150000,0.150000,0.150000}%
\pgfsetstrokecolor{currentstroke}%
\pgfsetdash{}{0pt}%
\pgfsys@defobject{currentmarker}{\pgfqpoint{0.000000in}{-0.041667in}}{\pgfqpoint{0.000000in}{0.000000in}}{%
\pgfpathmoveto{\pgfqpoint{0.000000in}{0.000000in}}%
\pgfpathlineto{\pgfqpoint{0.000000in}{-0.041667in}}%
\pgfusepath{stroke,fill}%
}%
\begin{pgfscope}%
\pgfsys@transformshift{6.666218in}{0.750000in}%
\pgfsys@useobject{currentmarker}{}%
\end{pgfscope}%
\end{pgfscope}%
\begin{pgfscope}%
\pgfsetbuttcap%
\pgfsetroundjoin%
\definecolor{currentfill}{rgb}{0.150000,0.150000,0.150000}%
\pgfsetfillcolor{currentfill}%
\pgfsetlinewidth{0.401500pt}%
\definecolor{currentstroke}{rgb}{0.150000,0.150000,0.150000}%
\pgfsetstrokecolor{currentstroke}%
\pgfsetdash{}{0pt}%
\pgfsys@defobject{currentmarker}{\pgfqpoint{0.000000in}{-0.041667in}}{\pgfqpoint{0.000000in}{0.000000in}}{%
\pgfpathmoveto{\pgfqpoint{0.000000in}{0.000000in}}%
\pgfpathlineto{\pgfqpoint{0.000000in}{-0.041667in}}%
\pgfusepath{stroke,fill}%
}%
\begin{pgfscope}%
\pgfsys@transformshift{6.959958in}{0.750000in}%
\pgfsys@useobject{currentmarker}{}%
\end{pgfscope}%
\end{pgfscope}%
\begin{pgfscope}%
\pgfsetbuttcap%
\pgfsetroundjoin%
\definecolor{currentfill}{rgb}{0.150000,0.150000,0.150000}%
\pgfsetfillcolor{currentfill}%
\pgfsetlinewidth{0.401500pt}%
\definecolor{currentstroke}{rgb}{0.150000,0.150000,0.150000}%
\pgfsetstrokecolor{currentstroke}%
\pgfsetdash{}{0pt}%
\pgfsys@defobject{currentmarker}{\pgfqpoint{0.000000in}{-0.041667in}}{\pgfqpoint{0.000000in}{0.000000in}}{%
\pgfpathmoveto{\pgfqpoint{0.000000in}{0.000000in}}%
\pgfpathlineto{\pgfqpoint{0.000000in}{-0.041667in}}%
\pgfusepath{stroke,fill}%
}%
\begin{pgfscope}%
\pgfsys@transformshift{7.208312in}{0.750000in}%
\pgfsys@useobject{currentmarker}{}%
\end{pgfscope}%
\end{pgfscope}%
\begin{pgfscope}%
\pgfsetbuttcap%
\pgfsetroundjoin%
\definecolor{currentfill}{rgb}{0.150000,0.150000,0.150000}%
\pgfsetfillcolor{currentfill}%
\pgfsetlinewidth{0.401500pt}%
\definecolor{currentstroke}{rgb}{0.150000,0.150000,0.150000}%
\pgfsetstrokecolor{currentstroke}%
\pgfsetdash{}{0pt}%
\pgfsys@defobject{currentmarker}{\pgfqpoint{0.000000in}{-0.041667in}}{\pgfqpoint{0.000000in}{0.000000in}}{%
\pgfpathmoveto{\pgfqpoint{0.000000in}{0.000000in}}%
\pgfpathlineto{\pgfqpoint{0.000000in}{-0.041667in}}%
\pgfusepath{stroke,fill}%
}%
\begin{pgfscope}%
\pgfsys@transformshift{7.423446in}{0.750000in}%
\pgfsys@useobject{currentmarker}{}%
\end{pgfscope}%
\end{pgfscope}%
\begin{pgfscope}%
\pgfsetbuttcap%
\pgfsetroundjoin%
\definecolor{currentfill}{rgb}{0.150000,0.150000,0.150000}%
\pgfsetfillcolor{currentfill}%
\pgfsetlinewidth{0.401500pt}%
\definecolor{currentstroke}{rgb}{0.150000,0.150000,0.150000}%
\pgfsetstrokecolor{currentstroke}%
\pgfsetdash{}{0pt}%
\pgfsys@defobject{currentmarker}{\pgfqpoint{0.000000in}{-0.041667in}}{\pgfqpoint{0.000000in}{0.000000in}}{%
\pgfpathmoveto{\pgfqpoint{0.000000in}{0.000000in}}%
\pgfpathlineto{\pgfqpoint{0.000000in}{-0.041667in}}%
\pgfusepath{stroke,fill}%
}%
\begin{pgfscope}%
\pgfsys@transformshift{7.613207in}{0.750000in}%
\pgfsys@useobject{currentmarker}{}%
\end{pgfscope}%
\end{pgfscope}%
\begin{pgfscope}%
\definecolor{textcolor}{rgb}{0.150000,0.150000,0.150000}%
\pgfsetstrokecolor{textcolor}%
\pgfsetfillcolor{textcolor}%
\pgftext[x=4.612500in,y=0.318273in,,top]{\color{textcolor}\sffamily\fontsize{24.000000}{28.800000}\selectfont \(\displaystyle \varepsilon\)}%
\end{pgfscope}%
\begin{pgfscope}%
\pgfsetbuttcap%
\pgfsetroundjoin%
\definecolor{currentfill}{rgb}{0.150000,0.150000,0.150000}%
\pgfsetfillcolor{currentfill}%
\pgfsetlinewidth{0.803000pt}%
\definecolor{currentstroke}{rgb}{0.150000,0.150000,0.150000}%
\pgfsetstrokecolor{currentstroke}%
\pgfsetdash{}{0pt}%
\pgfsys@defobject{currentmarker}{\pgfqpoint{-0.083333in}{0.000000in}}{\pgfqpoint{0.000000in}{0.000000in}}{%
\pgfpathmoveto{\pgfqpoint{0.000000in}{0.000000in}}%
\pgfpathlineto{\pgfqpoint{-0.083333in}{0.000000in}}%
\pgfusepath{stroke,fill}%
}%
\begin{pgfscope}%
\pgfsys@transformshift{1.125000in}{2.718314in}%
\pgfsys@useobject{currentmarker}{}%
\end{pgfscope}%
\end{pgfscope}%
\begin{pgfscope}%
\definecolor{textcolor}{rgb}{0.150000,0.150000,0.150000}%
\pgfsetstrokecolor{textcolor}%
\pgfsetfillcolor{textcolor}%
\pgftext[x=0.439258in,y=2.602239in,left,base]{\color{textcolor}\sffamily\fontsize{22.000000}{26.400000}\selectfont \(\displaystyle 10^{-1}\)}%
\end{pgfscope}%
\begin{pgfscope}%
\pgfsetbuttcap%
\pgfsetroundjoin%
\definecolor{currentfill}{rgb}{0.150000,0.150000,0.150000}%
\pgfsetfillcolor{currentfill}%
\pgfsetlinewidth{0.803000pt}%
\definecolor{currentstroke}{rgb}{0.150000,0.150000,0.150000}%
\pgfsetstrokecolor{currentstroke}%
\pgfsetdash{}{0pt}%
\pgfsys@defobject{currentmarker}{\pgfqpoint{-0.083333in}{0.000000in}}{\pgfqpoint{0.000000in}{0.000000in}}{%
\pgfpathmoveto{\pgfqpoint{0.000000in}{0.000000in}}%
\pgfpathlineto{\pgfqpoint{-0.083333in}{0.000000in}}%
\pgfusepath{stroke,fill}%
}%
\begin{pgfscope}%
\pgfsys@transformshift{1.125000in}{5.074091in}%
\pgfsys@useobject{currentmarker}{}%
\end{pgfscope}%
\end{pgfscope}%
\begin{pgfscope}%
\definecolor{textcolor}{rgb}{0.150000,0.150000,0.150000}%
\pgfsetstrokecolor{textcolor}%
\pgfsetfillcolor{textcolor}%
\pgftext[x=0.594814in,y=4.958016in,left,base]{\color{textcolor}\sffamily\fontsize{22.000000}{26.400000}\selectfont \(\displaystyle 10^{0}\)}%
\end{pgfscope}%
\begin{pgfscope}%
\pgfsetbuttcap%
\pgfsetroundjoin%
\definecolor{currentfill}{rgb}{0.150000,0.150000,0.150000}%
\pgfsetfillcolor{currentfill}%
\pgfsetlinewidth{0.401500pt}%
\definecolor{currentstroke}{rgb}{0.150000,0.150000,0.150000}%
\pgfsetstrokecolor{currentstroke}%
\pgfsetdash{}{0pt}%
\pgfsys@defobject{currentmarker}{\pgfqpoint{-0.041667in}{0.000000in}}{\pgfqpoint{0.000000in}{0.000000in}}{%
\pgfpathmoveto{\pgfqpoint{0.000000in}{0.000000in}}%
\pgfpathlineto{\pgfqpoint{-0.041667in}{0.000000in}}%
\pgfusepath{stroke,fill}%
}%
\begin{pgfscope}%
\pgfsys@transformshift{1.125000in}{1.071697in}%
\pgfsys@useobject{currentmarker}{}%
\end{pgfscope}%
\end{pgfscope}%
\begin{pgfscope}%
\pgfsetbuttcap%
\pgfsetroundjoin%
\definecolor{currentfill}{rgb}{0.150000,0.150000,0.150000}%
\pgfsetfillcolor{currentfill}%
\pgfsetlinewidth{0.401500pt}%
\definecolor{currentstroke}{rgb}{0.150000,0.150000,0.150000}%
\pgfsetstrokecolor{currentstroke}%
\pgfsetdash{}{0pt}%
\pgfsys@defobject{currentmarker}{\pgfqpoint{-0.041667in}{0.000000in}}{\pgfqpoint{0.000000in}{0.000000in}}{%
\pgfpathmoveto{\pgfqpoint{0.000000in}{0.000000in}}%
\pgfpathlineto{\pgfqpoint{-0.041667in}{0.000000in}}%
\pgfusepath{stroke,fill}%
}%
\begin{pgfscope}%
\pgfsys@transformshift{1.125000in}{1.486528in}%
\pgfsys@useobject{currentmarker}{}%
\end{pgfscope}%
\end{pgfscope}%
\begin{pgfscope}%
\pgfsetbuttcap%
\pgfsetroundjoin%
\definecolor{currentfill}{rgb}{0.150000,0.150000,0.150000}%
\pgfsetfillcolor{currentfill}%
\pgfsetlinewidth{0.401500pt}%
\definecolor{currentstroke}{rgb}{0.150000,0.150000,0.150000}%
\pgfsetstrokecolor{currentstroke}%
\pgfsetdash{}{0pt}%
\pgfsys@defobject{currentmarker}{\pgfqpoint{-0.041667in}{0.000000in}}{\pgfqpoint{0.000000in}{0.000000in}}{%
\pgfpathmoveto{\pgfqpoint{0.000000in}{0.000000in}}%
\pgfpathlineto{\pgfqpoint{-0.041667in}{0.000000in}}%
\pgfusepath{stroke,fill}%
}%
\begin{pgfscope}%
\pgfsys@transformshift{1.125000in}{1.780856in}%
\pgfsys@useobject{currentmarker}{}%
\end{pgfscope}%
\end{pgfscope}%
\begin{pgfscope}%
\pgfsetbuttcap%
\pgfsetroundjoin%
\definecolor{currentfill}{rgb}{0.150000,0.150000,0.150000}%
\pgfsetfillcolor{currentfill}%
\pgfsetlinewidth{0.401500pt}%
\definecolor{currentstroke}{rgb}{0.150000,0.150000,0.150000}%
\pgfsetstrokecolor{currentstroke}%
\pgfsetdash{}{0pt}%
\pgfsys@defobject{currentmarker}{\pgfqpoint{-0.041667in}{0.000000in}}{\pgfqpoint{0.000000in}{0.000000in}}{%
\pgfpathmoveto{\pgfqpoint{0.000000in}{0.000000in}}%
\pgfpathlineto{\pgfqpoint{-0.041667in}{0.000000in}}%
\pgfusepath{stroke,fill}%
}%
\begin{pgfscope}%
\pgfsys@transformshift{1.125000in}{2.009155in}%
\pgfsys@useobject{currentmarker}{}%
\end{pgfscope}%
\end{pgfscope}%
\begin{pgfscope}%
\pgfsetbuttcap%
\pgfsetroundjoin%
\definecolor{currentfill}{rgb}{0.150000,0.150000,0.150000}%
\pgfsetfillcolor{currentfill}%
\pgfsetlinewidth{0.401500pt}%
\definecolor{currentstroke}{rgb}{0.150000,0.150000,0.150000}%
\pgfsetstrokecolor{currentstroke}%
\pgfsetdash{}{0pt}%
\pgfsys@defobject{currentmarker}{\pgfqpoint{-0.041667in}{0.000000in}}{\pgfqpoint{0.000000in}{0.000000in}}{%
\pgfpathmoveto{\pgfqpoint{0.000000in}{0.000000in}}%
\pgfpathlineto{\pgfqpoint{-0.041667in}{0.000000in}}%
\pgfusepath{stroke,fill}%
}%
\begin{pgfscope}%
\pgfsys@transformshift{1.125000in}{2.195688in}%
\pgfsys@useobject{currentmarker}{}%
\end{pgfscope}%
\end{pgfscope}%
\begin{pgfscope}%
\pgfsetbuttcap%
\pgfsetroundjoin%
\definecolor{currentfill}{rgb}{0.150000,0.150000,0.150000}%
\pgfsetfillcolor{currentfill}%
\pgfsetlinewidth{0.401500pt}%
\definecolor{currentstroke}{rgb}{0.150000,0.150000,0.150000}%
\pgfsetstrokecolor{currentstroke}%
\pgfsetdash{}{0pt}%
\pgfsys@defobject{currentmarker}{\pgfqpoint{-0.041667in}{0.000000in}}{\pgfqpoint{0.000000in}{0.000000in}}{%
\pgfpathmoveto{\pgfqpoint{0.000000in}{0.000000in}}%
\pgfpathlineto{\pgfqpoint{-0.041667in}{0.000000in}}%
\pgfusepath{stroke,fill}%
}%
\begin{pgfscope}%
\pgfsys@transformshift{1.125000in}{2.353400in}%
\pgfsys@useobject{currentmarker}{}%
\end{pgfscope}%
\end{pgfscope}%
\begin{pgfscope}%
\pgfsetbuttcap%
\pgfsetroundjoin%
\definecolor{currentfill}{rgb}{0.150000,0.150000,0.150000}%
\pgfsetfillcolor{currentfill}%
\pgfsetlinewidth{0.401500pt}%
\definecolor{currentstroke}{rgb}{0.150000,0.150000,0.150000}%
\pgfsetstrokecolor{currentstroke}%
\pgfsetdash{}{0pt}%
\pgfsys@defobject{currentmarker}{\pgfqpoint{-0.041667in}{0.000000in}}{\pgfqpoint{0.000000in}{0.000000in}}{%
\pgfpathmoveto{\pgfqpoint{0.000000in}{0.000000in}}%
\pgfpathlineto{\pgfqpoint{-0.041667in}{0.000000in}}%
\pgfusepath{stroke,fill}%
}%
\begin{pgfscope}%
\pgfsys@transformshift{1.125000in}{2.490016in}%
\pgfsys@useobject{currentmarker}{}%
\end{pgfscope}%
\end{pgfscope}%
\begin{pgfscope}%
\pgfsetbuttcap%
\pgfsetroundjoin%
\definecolor{currentfill}{rgb}{0.150000,0.150000,0.150000}%
\pgfsetfillcolor{currentfill}%
\pgfsetlinewidth{0.401500pt}%
\definecolor{currentstroke}{rgb}{0.150000,0.150000,0.150000}%
\pgfsetstrokecolor{currentstroke}%
\pgfsetdash{}{0pt}%
\pgfsys@defobject{currentmarker}{\pgfqpoint{-0.041667in}{0.000000in}}{\pgfqpoint{0.000000in}{0.000000in}}{%
\pgfpathmoveto{\pgfqpoint{0.000000in}{0.000000in}}%
\pgfpathlineto{\pgfqpoint{-0.041667in}{0.000000in}}%
\pgfusepath{stroke,fill}%
}%
\begin{pgfscope}%
\pgfsys@transformshift{1.125000in}{2.610520in}%
\pgfsys@useobject{currentmarker}{}%
\end{pgfscope}%
\end{pgfscope}%
\begin{pgfscope}%
\pgfsetbuttcap%
\pgfsetroundjoin%
\definecolor{currentfill}{rgb}{0.150000,0.150000,0.150000}%
\pgfsetfillcolor{currentfill}%
\pgfsetlinewidth{0.401500pt}%
\definecolor{currentstroke}{rgb}{0.150000,0.150000,0.150000}%
\pgfsetstrokecolor{currentstroke}%
\pgfsetdash{}{0pt}%
\pgfsys@defobject{currentmarker}{\pgfqpoint{-0.041667in}{0.000000in}}{\pgfqpoint{0.000000in}{0.000000in}}{%
\pgfpathmoveto{\pgfqpoint{0.000000in}{0.000000in}}%
\pgfpathlineto{\pgfqpoint{-0.041667in}{0.000000in}}%
\pgfusepath{stroke,fill}%
}%
\begin{pgfscope}%
\pgfsys@transformshift{1.125000in}{3.427474in}%
\pgfsys@useobject{currentmarker}{}%
\end{pgfscope}%
\end{pgfscope}%
\begin{pgfscope}%
\pgfsetbuttcap%
\pgfsetroundjoin%
\definecolor{currentfill}{rgb}{0.150000,0.150000,0.150000}%
\pgfsetfillcolor{currentfill}%
\pgfsetlinewidth{0.401500pt}%
\definecolor{currentstroke}{rgb}{0.150000,0.150000,0.150000}%
\pgfsetstrokecolor{currentstroke}%
\pgfsetdash{}{0pt}%
\pgfsys@defobject{currentmarker}{\pgfqpoint{-0.041667in}{0.000000in}}{\pgfqpoint{0.000000in}{0.000000in}}{%
\pgfpathmoveto{\pgfqpoint{0.000000in}{0.000000in}}%
\pgfpathlineto{\pgfqpoint{-0.041667in}{0.000000in}}%
\pgfusepath{stroke,fill}%
}%
\begin{pgfscope}%
\pgfsys@transformshift{1.125000in}{3.842305in}%
\pgfsys@useobject{currentmarker}{}%
\end{pgfscope}%
\end{pgfscope}%
\begin{pgfscope}%
\pgfsetbuttcap%
\pgfsetroundjoin%
\definecolor{currentfill}{rgb}{0.150000,0.150000,0.150000}%
\pgfsetfillcolor{currentfill}%
\pgfsetlinewidth{0.401500pt}%
\definecolor{currentstroke}{rgb}{0.150000,0.150000,0.150000}%
\pgfsetstrokecolor{currentstroke}%
\pgfsetdash{}{0pt}%
\pgfsys@defobject{currentmarker}{\pgfqpoint{-0.041667in}{0.000000in}}{\pgfqpoint{0.000000in}{0.000000in}}{%
\pgfpathmoveto{\pgfqpoint{0.000000in}{0.000000in}}%
\pgfpathlineto{\pgfqpoint{-0.041667in}{0.000000in}}%
\pgfusepath{stroke,fill}%
}%
\begin{pgfscope}%
\pgfsys@transformshift{1.125000in}{4.136633in}%
\pgfsys@useobject{currentmarker}{}%
\end{pgfscope}%
\end{pgfscope}%
\begin{pgfscope}%
\pgfsetbuttcap%
\pgfsetroundjoin%
\definecolor{currentfill}{rgb}{0.150000,0.150000,0.150000}%
\pgfsetfillcolor{currentfill}%
\pgfsetlinewidth{0.401500pt}%
\definecolor{currentstroke}{rgb}{0.150000,0.150000,0.150000}%
\pgfsetstrokecolor{currentstroke}%
\pgfsetdash{}{0pt}%
\pgfsys@defobject{currentmarker}{\pgfqpoint{-0.041667in}{0.000000in}}{\pgfqpoint{0.000000in}{0.000000in}}{%
\pgfpathmoveto{\pgfqpoint{0.000000in}{0.000000in}}%
\pgfpathlineto{\pgfqpoint{-0.041667in}{0.000000in}}%
\pgfusepath{stroke,fill}%
}%
\begin{pgfscope}%
\pgfsys@transformshift{1.125000in}{4.364931in}%
\pgfsys@useobject{currentmarker}{}%
\end{pgfscope}%
\end{pgfscope}%
\begin{pgfscope}%
\pgfsetbuttcap%
\pgfsetroundjoin%
\definecolor{currentfill}{rgb}{0.150000,0.150000,0.150000}%
\pgfsetfillcolor{currentfill}%
\pgfsetlinewidth{0.401500pt}%
\definecolor{currentstroke}{rgb}{0.150000,0.150000,0.150000}%
\pgfsetstrokecolor{currentstroke}%
\pgfsetdash{}{0pt}%
\pgfsys@defobject{currentmarker}{\pgfqpoint{-0.041667in}{0.000000in}}{\pgfqpoint{0.000000in}{0.000000in}}{%
\pgfpathmoveto{\pgfqpoint{0.000000in}{0.000000in}}%
\pgfpathlineto{\pgfqpoint{-0.041667in}{0.000000in}}%
\pgfusepath{stroke,fill}%
}%
\begin{pgfscope}%
\pgfsys@transformshift{1.125000in}{4.551465in}%
\pgfsys@useobject{currentmarker}{}%
\end{pgfscope}%
\end{pgfscope}%
\begin{pgfscope}%
\pgfsetbuttcap%
\pgfsetroundjoin%
\definecolor{currentfill}{rgb}{0.150000,0.150000,0.150000}%
\pgfsetfillcolor{currentfill}%
\pgfsetlinewidth{0.401500pt}%
\definecolor{currentstroke}{rgb}{0.150000,0.150000,0.150000}%
\pgfsetstrokecolor{currentstroke}%
\pgfsetdash{}{0pt}%
\pgfsys@defobject{currentmarker}{\pgfqpoint{-0.041667in}{0.000000in}}{\pgfqpoint{0.000000in}{0.000000in}}{%
\pgfpathmoveto{\pgfqpoint{0.000000in}{0.000000in}}%
\pgfpathlineto{\pgfqpoint{-0.041667in}{0.000000in}}%
\pgfusepath{stroke,fill}%
}%
\begin{pgfscope}%
\pgfsys@transformshift{1.125000in}{4.709176in}%
\pgfsys@useobject{currentmarker}{}%
\end{pgfscope}%
\end{pgfscope}%
\begin{pgfscope}%
\pgfsetbuttcap%
\pgfsetroundjoin%
\definecolor{currentfill}{rgb}{0.150000,0.150000,0.150000}%
\pgfsetfillcolor{currentfill}%
\pgfsetlinewidth{0.401500pt}%
\definecolor{currentstroke}{rgb}{0.150000,0.150000,0.150000}%
\pgfsetstrokecolor{currentstroke}%
\pgfsetdash{}{0pt}%
\pgfsys@defobject{currentmarker}{\pgfqpoint{-0.041667in}{0.000000in}}{\pgfqpoint{0.000000in}{0.000000in}}{%
\pgfpathmoveto{\pgfqpoint{0.000000in}{0.000000in}}%
\pgfpathlineto{\pgfqpoint{-0.041667in}{0.000000in}}%
\pgfusepath{stroke,fill}%
}%
\begin{pgfscope}%
\pgfsys@transformshift{1.125000in}{4.845793in}%
\pgfsys@useobject{currentmarker}{}%
\end{pgfscope}%
\end{pgfscope}%
\begin{pgfscope}%
\pgfsetbuttcap%
\pgfsetroundjoin%
\definecolor{currentfill}{rgb}{0.150000,0.150000,0.150000}%
\pgfsetfillcolor{currentfill}%
\pgfsetlinewidth{0.401500pt}%
\definecolor{currentstroke}{rgb}{0.150000,0.150000,0.150000}%
\pgfsetstrokecolor{currentstroke}%
\pgfsetdash{}{0pt}%
\pgfsys@defobject{currentmarker}{\pgfqpoint{-0.041667in}{0.000000in}}{\pgfqpoint{0.000000in}{0.000000in}}{%
\pgfpathmoveto{\pgfqpoint{0.000000in}{0.000000in}}%
\pgfpathlineto{\pgfqpoint{-0.041667in}{0.000000in}}%
\pgfusepath{stroke,fill}%
}%
\begin{pgfscope}%
\pgfsys@transformshift{1.125000in}{4.966296in}%
\pgfsys@useobject{currentmarker}{}%
\end{pgfscope}%
\end{pgfscope}%
\begin{pgfscope}%
\definecolor{textcolor}{rgb}{0.150000,0.150000,0.150000}%
\pgfsetstrokecolor{textcolor}%
\pgfsetfillcolor{textcolor}%
\pgftext[x=0.383703in,y=3.015000in,,bottom,rotate=90.000000]{\color{textcolor}\sffamily\fontsize{24.000000}{28.800000}\selectfont Spectral radius, \(\displaystyle \rho\)}%
\end{pgfscope}%
\begin{pgfscope}%
\pgfpathrectangle{\pgfqpoint{1.125000in}{0.750000in}}{\pgfqpoint{6.975000in}{4.530000in}} %
\pgfusepath{clip}%
\pgfsetroundcap%
\pgfsetroundjoin%
\pgfsetlinewidth{1.405250pt}%
\definecolor{currentstroke}{rgb}{0.133333,0.133333,0.133333}%
\pgfsetstrokecolor{currentstroke}%
\pgfsetdash{}{0pt}%
\pgfpathmoveto{\pgfqpoint{1.442045in}{3.351104in}}%
\pgfpathlineto{\pgfqpoint{4.476888in}{2.216210in}}%
\pgfpathlineto{\pgfqpoint{5.466254in}{2.222248in}}%
\pgfpathlineto{\pgfqpoint{6.074705in}{2.223528in}}%
\pgfpathlineto{\pgfqpoint{6.515318in}{2.223996in}}%
\pgfpathlineto{\pgfqpoint{6.860945in}{2.224218in}}%
\pgfpathlineto{\pgfqpoint{7.145363in}{2.224340in}}%
\pgfpathlineto{\pgfqpoint{7.387025in}{2.224415in}}%
\pgfpathlineto{\pgfqpoint{7.597122in}{2.224463in}}%
\pgfpathlineto{\pgfqpoint{7.782955in}{2.224497in}}%
\pgfusepath{stroke}%
\end{pgfscope}%
\begin{pgfscope}%
\pgfpathrectangle{\pgfqpoint{1.125000in}{0.750000in}}{\pgfqpoint{6.975000in}{4.530000in}} %
\pgfusepath{clip}%
\pgfsetbuttcap%
\pgfsetmiterjoin%
\definecolor{currentfill}{rgb}{0.133333,0.133333,0.133333}%
\pgfsetfillcolor{currentfill}%
\pgfsetlinewidth{0.000000pt}%
\definecolor{currentstroke}{rgb}{0.133333,0.133333,0.133333}%
\pgfsetstrokecolor{currentstroke}%
\pgfsetdash{}{0pt}%
\pgfsys@defobject{currentmarker}{\pgfqpoint{-0.062500in}{-0.062500in}}{\pgfqpoint{0.062500in}{0.062500in}}{%
\pgfpathmoveto{\pgfqpoint{-0.000000in}{-0.062500in}}%
\pgfpathlineto{\pgfqpoint{0.062500in}{0.062500in}}%
\pgfpathlineto{\pgfqpoint{-0.062500in}{0.062500in}}%
\pgfpathclose%
\pgfusepath{fill}%
}%
\begin{pgfscope}%
\pgfsys@transformshift{1.442045in}{3.351104in}%
\pgfsys@useobject{currentmarker}{}%
\end{pgfscope}%
\begin{pgfscope}%
\pgfsys@transformshift{4.476888in}{2.216210in}%
\pgfsys@useobject{currentmarker}{}%
\end{pgfscope}%
\begin{pgfscope}%
\pgfsys@transformshift{5.466254in}{2.222248in}%
\pgfsys@useobject{currentmarker}{}%
\end{pgfscope}%
\begin{pgfscope}%
\pgfsys@transformshift{6.074705in}{2.223528in}%
\pgfsys@useobject{currentmarker}{}%
\end{pgfscope}%
\begin{pgfscope}%
\pgfsys@transformshift{6.515318in}{2.223996in}%
\pgfsys@useobject{currentmarker}{}%
\end{pgfscope}%
\begin{pgfscope}%
\pgfsys@transformshift{6.860945in}{2.224218in}%
\pgfsys@useobject{currentmarker}{}%
\end{pgfscope}%
\begin{pgfscope}%
\pgfsys@transformshift{7.145363in}{2.224340in}%
\pgfsys@useobject{currentmarker}{}%
\end{pgfscope}%
\begin{pgfscope}%
\pgfsys@transformshift{7.387025in}{2.224415in}%
\pgfsys@useobject{currentmarker}{}%
\end{pgfscope}%
\begin{pgfscope}%
\pgfsys@transformshift{7.597122in}{2.224463in}%
\pgfsys@useobject{currentmarker}{}%
\end{pgfscope}%
\begin{pgfscope}%
\pgfsys@transformshift{7.782955in}{2.224497in}%
\pgfsys@useobject{currentmarker}{}%
\end{pgfscope}%
\end{pgfscope}%
\begin{pgfscope}%
\pgfpathrectangle{\pgfqpoint{1.125000in}{0.750000in}}{\pgfqpoint{6.975000in}{4.530000in}} %
\pgfusepath{clip}%
\pgfsetroundcap%
\pgfsetroundjoin%
\pgfsetlinewidth{1.405250pt}%
\definecolor{currentstroke}{rgb}{0.318370,0.318370,0.318370}%
\pgfsetstrokecolor{currentstroke}%
\pgfsetdash{}{0pt}%
\pgfpathmoveto{\pgfqpoint{1.442045in}{0.955909in}}%
\pgfpathlineto{\pgfqpoint{4.476888in}{1.909681in}}%
\pgfpathlineto{\pgfqpoint{5.466254in}{1.968295in}}%
\pgfpathlineto{\pgfqpoint{6.074705in}{2.083604in}}%
\pgfpathlineto{\pgfqpoint{6.515318in}{2.140550in}}%
\pgfpathlineto{\pgfqpoint{6.860945in}{2.158229in}}%
\pgfpathlineto{\pgfqpoint{7.145363in}{2.161745in}}%
\pgfpathlineto{\pgfqpoint{7.387025in}{2.176897in}}%
\pgfpathlineto{\pgfqpoint{7.597122in}{2.174171in}}%
\pgfpathlineto{\pgfqpoint{7.782955in}{2.186772in}}%
\pgfusepath{stroke}%
\end{pgfscope}%
\begin{pgfscope}%
\pgfpathrectangle{\pgfqpoint{1.125000in}{0.750000in}}{\pgfqpoint{6.975000in}{4.530000in}} %
\pgfusepath{clip}%
\pgfsetbuttcap%
\pgfsetroundjoin%
\definecolor{currentfill}{rgb}{0.318370,0.318370,0.318370}%
\pgfsetfillcolor{currentfill}%
\pgfsetlinewidth{0.000000pt}%
\definecolor{currentstroke}{rgb}{0.318370,0.318370,0.318370}%
\pgfsetstrokecolor{currentstroke}%
\pgfsetdash{}{0pt}%
\pgfsys@defobject{currentmarker}{\pgfqpoint{-0.048611in}{-0.048611in}}{\pgfqpoint{0.048611in}{0.048611in}}{%
\pgfpathmoveto{\pgfqpoint{0.000000in}{-0.048611in}}%
\pgfpathcurveto{\pgfqpoint{0.012892in}{-0.048611in}}{\pgfqpoint{0.025257in}{-0.043489in}}{\pgfqpoint{0.034373in}{-0.034373in}}%
\pgfpathcurveto{\pgfqpoint{0.043489in}{-0.025257in}}{\pgfqpoint{0.048611in}{-0.012892in}}{\pgfqpoint{0.048611in}{0.000000in}}%
\pgfpathcurveto{\pgfqpoint{0.048611in}{0.012892in}}{\pgfqpoint{0.043489in}{0.025257in}}{\pgfqpoint{0.034373in}{0.034373in}}%
\pgfpathcurveto{\pgfqpoint{0.025257in}{0.043489in}}{\pgfqpoint{0.012892in}{0.048611in}}{\pgfqpoint{0.000000in}{0.048611in}}%
\pgfpathcurveto{\pgfqpoint{-0.012892in}{0.048611in}}{\pgfqpoint{-0.025257in}{0.043489in}}{\pgfqpoint{-0.034373in}{0.034373in}}%
\pgfpathcurveto{\pgfqpoint{-0.043489in}{0.025257in}}{\pgfqpoint{-0.048611in}{0.012892in}}{\pgfqpoint{-0.048611in}{0.000000in}}%
\pgfpathcurveto{\pgfqpoint{-0.048611in}{-0.012892in}}{\pgfqpoint{-0.043489in}{-0.025257in}}{\pgfqpoint{-0.034373in}{-0.034373in}}%
\pgfpathcurveto{\pgfqpoint{-0.025257in}{-0.043489in}}{\pgfqpoint{-0.012892in}{-0.048611in}}{\pgfqpoint{0.000000in}{-0.048611in}}%
\pgfpathclose%
\pgfusepath{fill}%
}%
\begin{pgfscope}%
\pgfsys@transformshift{1.442045in}{0.955909in}%
\pgfsys@useobject{currentmarker}{}%
\end{pgfscope}%
\begin{pgfscope}%
\pgfsys@transformshift{4.476888in}{1.909681in}%
\pgfsys@useobject{currentmarker}{}%
\end{pgfscope}%
\begin{pgfscope}%
\pgfsys@transformshift{5.466254in}{1.968295in}%
\pgfsys@useobject{currentmarker}{}%
\end{pgfscope}%
\begin{pgfscope}%
\pgfsys@transformshift{6.074705in}{2.083604in}%
\pgfsys@useobject{currentmarker}{}%
\end{pgfscope}%
\begin{pgfscope}%
\pgfsys@transformshift{6.515318in}{2.140550in}%
\pgfsys@useobject{currentmarker}{}%
\end{pgfscope}%
\begin{pgfscope}%
\pgfsys@transformshift{6.860945in}{2.158229in}%
\pgfsys@useobject{currentmarker}{}%
\end{pgfscope}%
\begin{pgfscope}%
\pgfsys@transformshift{7.145363in}{2.161745in}%
\pgfsys@useobject{currentmarker}{}%
\end{pgfscope}%
\begin{pgfscope}%
\pgfsys@transformshift{7.387025in}{2.176897in}%
\pgfsys@useobject{currentmarker}{}%
\end{pgfscope}%
\begin{pgfscope}%
\pgfsys@transformshift{7.597122in}{2.174171in}%
\pgfsys@useobject{currentmarker}{}%
\end{pgfscope}%
\begin{pgfscope}%
\pgfsys@transformshift{7.782955in}{2.186772in}%
\pgfsys@useobject{currentmarker}{}%
\end{pgfscope}%
\end{pgfscope}%
\begin{pgfscope}%
\pgfpathrectangle{\pgfqpoint{1.125000in}{0.750000in}}{\pgfqpoint{6.975000in}{4.530000in}} %
\pgfusepath{clip}%
\pgfsetroundcap%
\pgfsetroundjoin%
\pgfsetlinewidth{1.405250pt}%
\definecolor{currentstroke}{rgb}{0.501961,0.501961,0.501961}%
\pgfsetstrokecolor{currentstroke}%
\pgfsetdash{}{0pt}%
\pgfpathmoveto{\pgfqpoint{1.442045in}{4.049851in}}%
\pgfpathlineto{\pgfqpoint{4.476888in}{3.411725in}}%
\pgfpathlineto{\pgfqpoint{5.466254in}{3.411414in}}%
\pgfpathlineto{\pgfqpoint{6.074705in}{3.413126in}}%
\pgfpathlineto{\pgfqpoint{6.515318in}{3.410862in}}%
\pgfpathlineto{\pgfqpoint{6.860945in}{3.413494in}}%
\pgfpathlineto{\pgfqpoint{7.145363in}{3.410979in}}%
\pgfpathlineto{\pgfqpoint{7.387025in}{3.407534in}}%
\pgfpathlineto{\pgfqpoint{7.597122in}{3.408383in}}%
\pgfpathlineto{\pgfqpoint{7.782955in}{3.408087in}}%
\pgfusepath{stroke}%
\end{pgfscope}%
\begin{pgfscope}%
\pgfpathrectangle{\pgfqpoint{1.125000in}{0.750000in}}{\pgfqpoint{6.975000in}{4.530000in}} %
\pgfusepath{clip}%
\pgfsetbuttcap%
\pgfsetroundjoin%
\definecolor{currentfill}{rgb}{0.501961,0.501961,0.501961}%
\pgfsetfillcolor{currentfill}%
\pgfsetlinewidth{0.000000pt}%
\definecolor{currentstroke}{rgb}{0.501961,0.501961,0.501961}%
\pgfsetstrokecolor{currentstroke}%
\pgfsetdash{}{0pt}%
\pgfsys@defobject{currentmarker}{\pgfqpoint{-0.034722in}{-0.034722in}}{\pgfqpoint{0.034722in}{0.034722in}}{%
\pgfpathmoveto{\pgfqpoint{0.000000in}{-0.034722in}}%
\pgfpathcurveto{\pgfqpoint{0.009208in}{-0.034722in}}{\pgfqpoint{0.018041in}{-0.031064in}}{\pgfqpoint{0.024552in}{-0.024552in}}%
\pgfpathcurveto{\pgfqpoint{0.031064in}{-0.018041in}}{\pgfqpoint{0.034722in}{-0.009208in}}{\pgfqpoint{0.034722in}{0.000000in}}%
\pgfpathcurveto{\pgfqpoint{0.034722in}{0.009208in}}{\pgfqpoint{0.031064in}{0.018041in}}{\pgfqpoint{0.024552in}{0.024552in}}%
\pgfpathcurveto{\pgfqpoint{0.018041in}{0.031064in}}{\pgfqpoint{0.009208in}{0.034722in}}{\pgfqpoint{0.000000in}{0.034722in}}%
\pgfpathcurveto{\pgfqpoint{-0.009208in}{0.034722in}}{\pgfqpoint{-0.018041in}{0.031064in}}{\pgfqpoint{-0.024552in}{0.024552in}}%
\pgfpathcurveto{\pgfqpoint{-0.031064in}{0.018041in}}{\pgfqpoint{-0.034722in}{0.009208in}}{\pgfqpoint{-0.034722in}{0.000000in}}%
\pgfpathcurveto{\pgfqpoint{-0.034722in}{-0.009208in}}{\pgfqpoint{-0.031064in}{-0.018041in}}{\pgfqpoint{-0.024552in}{-0.024552in}}%
\pgfpathcurveto{\pgfqpoint{-0.018041in}{-0.031064in}}{\pgfqpoint{-0.009208in}{-0.034722in}}{\pgfqpoint{0.000000in}{-0.034722in}}%
\pgfpathclose%
\pgfusepath{fill}%
}%
\begin{pgfscope}%
\pgfsys@transformshift{1.442045in}{4.049851in}%
\pgfsys@useobject{currentmarker}{}%
\end{pgfscope}%
\begin{pgfscope}%
\pgfsys@transformshift{4.476888in}{3.411725in}%
\pgfsys@useobject{currentmarker}{}%
\end{pgfscope}%
\begin{pgfscope}%
\pgfsys@transformshift{5.466254in}{3.411414in}%
\pgfsys@useobject{currentmarker}{}%
\end{pgfscope}%
\begin{pgfscope}%
\pgfsys@transformshift{6.074705in}{3.413126in}%
\pgfsys@useobject{currentmarker}{}%
\end{pgfscope}%
\begin{pgfscope}%
\pgfsys@transformshift{6.515318in}{3.410862in}%
\pgfsys@useobject{currentmarker}{}%
\end{pgfscope}%
\begin{pgfscope}%
\pgfsys@transformshift{6.860945in}{3.413494in}%
\pgfsys@useobject{currentmarker}{}%
\end{pgfscope}%
\begin{pgfscope}%
\pgfsys@transformshift{7.145363in}{3.410979in}%
\pgfsys@useobject{currentmarker}{}%
\end{pgfscope}%
\begin{pgfscope}%
\pgfsys@transformshift{7.387025in}{3.407534in}%
\pgfsys@useobject{currentmarker}{}%
\end{pgfscope}%
\begin{pgfscope}%
\pgfsys@transformshift{7.597122in}{3.408383in}%
\pgfsys@useobject{currentmarker}{}%
\end{pgfscope}%
\begin{pgfscope}%
\pgfsys@transformshift{7.782955in}{3.408087in}%
\pgfsys@useobject{currentmarker}{}%
\end{pgfscope}%
\end{pgfscope}%
\begin{pgfscope}%
\pgfpathrectangle{\pgfqpoint{1.125000in}{0.750000in}}{\pgfqpoint{6.975000in}{4.530000in}} %
\pgfusepath{clip}%
\pgfsetbuttcap%
\pgfsetroundjoin%
\pgfsetlinewidth{1.405250pt}%
\definecolor{currentstroke}{rgb}{0.000000,0.000000,0.000000}%
\pgfsetstrokecolor{currentstroke}%
\pgfsetdash{{5.180000pt}{2.240000pt}}{0.000000pt}%
\pgfpathmoveto{\pgfqpoint{1.442045in}{5.074091in}}%
\pgfpathlineto{\pgfqpoint{4.476888in}{5.074091in}}%
\pgfpathlineto{\pgfqpoint{5.466254in}{5.074091in}}%
\pgfpathlineto{\pgfqpoint{6.074705in}{5.074091in}}%
\pgfpathlineto{\pgfqpoint{6.515318in}{5.074091in}}%
\pgfpathlineto{\pgfqpoint{6.860945in}{5.074091in}}%
\pgfpathlineto{\pgfqpoint{7.145363in}{5.074091in}}%
\pgfpathlineto{\pgfqpoint{7.387025in}{5.074091in}}%
\pgfpathlineto{\pgfqpoint{7.597122in}{5.074091in}}%
\pgfpathlineto{\pgfqpoint{7.782955in}{5.074091in}}%
\pgfusepath{stroke}%
\end{pgfscope}%
\begin{pgfscope}%
\pgfsetrectcap%
\pgfsetmiterjoin%
\pgfsetlinewidth{1.254687pt}%
\definecolor{currentstroke}{rgb}{0.150000,0.150000,0.150000}%
\pgfsetstrokecolor{currentstroke}%
\pgfsetdash{}{0pt}%
\pgfpathmoveto{\pgfqpoint{1.125000in}{0.750000in}}%
\pgfpathlineto{\pgfqpoint{1.125000in}{5.280000in}}%
\pgfusepath{stroke}%
\end{pgfscope}%
\begin{pgfscope}%
\pgfsetrectcap%
\pgfsetmiterjoin%
\pgfsetlinewidth{1.254687pt}%
\definecolor{currentstroke}{rgb}{0.150000,0.150000,0.150000}%
\pgfsetstrokecolor{currentstroke}%
\pgfsetdash{}{0pt}%
\pgfpathmoveto{\pgfqpoint{8.100000in}{0.750000in}}%
\pgfpathlineto{\pgfqpoint{8.100000in}{5.280000in}}%
\pgfusepath{stroke}%
\end{pgfscope}%
\begin{pgfscope}%
\pgfsetrectcap%
\pgfsetmiterjoin%
\pgfsetlinewidth{1.254687pt}%
\definecolor{currentstroke}{rgb}{0.150000,0.150000,0.150000}%
\pgfsetstrokecolor{currentstroke}%
\pgfsetdash{}{0pt}%
\pgfpathmoveto{\pgfqpoint{1.125000in}{0.750000in}}%
\pgfpathlineto{\pgfqpoint{8.100000in}{0.750000in}}%
\pgfusepath{stroke}%
\end{pgfscope}%
\begin{pgfscope}%
\pgfsetrectcap%
\pgfsetmiterjoin%
\pgfsetlinewidth{1.254687pt}%
\definecolor{currentstroke}{rgb}{0.150000,0.150000,0.150000}%
\pgfsetstrokecolor{currentstroke}%
\pgfsetdash{}{0pt}%
\pgfpathmoveto{\pgfqpoint{1.125000in}{5.280000in}}%
\pgfpathlineto{\pgfqpoint{8.100000in}{5.280000in}}%
\pgfusepath{stroke}%
\end{pgfscope}%
\begin{pgfscope}%
\pgfsetroundcap%
\pgfsetroundjoin%
\pgfsetlinewidth{1.405250pt}%
\definecolor{currentstroke}{rgb}{0.133333,0.133333,0.133333}%
\pgfsetstrokecolor{currentstroke}%
\pgfsetdash{}{0pt}%
\pgfpathmoveto{\pgfqpoint{5.797283in}{1.836132in}}%
\pgfpathlineto{\pgfqpoint{6.213950in}{1.836132in}}%
\pgfusepath{stroke}%
\end{pgfscope}%
\begin{pgfscope}%
\pgfsetbuttcap%
\pgfsetmiterjoin%
\definecolor{currentfill}{rgb}{0.133333,0.133333,0.133333}%
\pgfsetfillcolor{currentfill}%
\pgfsetlinewidth{0.000000pt}%
\definecolor{currentstroke}{rgb}{0.133333,0.133333,0.133333}%
\pgfsetstrokecolor{currentstroke}%
\pgfsetdash{}{0pt}%
\pgfsys@defobject{currentmarker}{\pgfqpoint{-0.062500in}{-0.062500in}}{\pgfqpoint{0.062500in}{0.062500in}}{%
\pgfpathmoveto{\pgfqpoint{-0.000000in}{-0.062500in}}%
\pgfpathlineto{\pgfqpoint{0.062500in}{0.062500in}}%
\pgfpathlineto{\pgfqpoint{-0.062500in}{0.062500in}}%
\pgfpathclose%
\pgfusepath{fill}%
}%
\begin{pgfscope}%
\pgfsys@transformshift{6.005617in}{1.836132in}%
\pgfsys@useobject{currentmarker}{}%
\end{pgfscope}%
\end{pgfscope}%
\begin{pgfscope}%
\definecolor{textcolor}{rgb}{0.150000,0.150000,0.150000}%
\pgfsetstrokecolor{textcolor}%
\pgfsetfillcolor{textcolor}%
\pgftext[x=6.380617in,y=1.763216in,left,base]{\color{textcolor}\sffamily\fontsize{15.000000}{18.000000}\selectfont Linear}%
\end{pgfscope}%
\begin{pgfscope}%
\pgfsetroundcap%
\pgfsetroundjoin%
\pgfsetlinewidth{1.405250pt}%
\definecolor{currentstroke}{rgb}{0.318370,0.318370,0.318370}%
\pgfsetstrokecolor{currentstroke}%
\pgfsetdash{}{0pt}%
\pgfpathmoveto{\pgfqpoint{5.797283in}{1.530346in}}%
\pgfpathlineto{\pgfqpoint{6.213950in}{1.530346in}}%
\pgfusepath{stroke}%
\end{pgfscope}%
\begin{pgfscope}%
\pgfsetbuttcap%
\pgfsetroundjoin%
\definecolor{currentfill}{rgb}{0.318370,0.318370,0.318370}%
\pgfsetfillcolor{currentfill}%
\pgfsetlinewidth{0.000000pt}%
\definecolor{currentstroke}{rgb}{0.318370,0.318370,0.318370}%
\pgfsetstrokecolor{currentstroke}%
\pgfsetdash{}{0pt}%
\pgfsys@defobject{currentmarker}{\pgfqpoint{-0.048611in}{-0.048611in}}{\pgfqpoint{0.048611in}{0.048611in}}{%
\pgfpathmoveto{\pgfqpoint{0.000000in}{-0.048611in}}%
\pgfpathcurveto{\pgfqpoint{0.012892in}{-0.048611in}}{\pgfqpoint{0.025257in}{-0.043489in}}{\pgfqpoint{0.034373in}{-0.034373in}}%
\pgfpathcurveto{\pgfqpoint{0.043489in}{-0.025257in}}{\pgfqpoint{0.048611in}{-0.012892in}}{\pgfqpoint{0.048611in}{0.000000in}}%
\pgfpathcurveto{\pgfqpoint{0.048611in}{0.012892in}}{\pgfqpoint{0.043489in}{0.025257in}}{\pgfqpoint{0.034373in}{0.034373in}}%
\pgfpathcurveto{\pgfqpoint{0.025257in}{0.043489in}}{\pgfqpoint{0.012892in}{0.048611in}}{\pgfqpoint{0.000000in}{0.048611in}}%
\pgfpathcurveto{\pgfqpoint{-0.012892in}{0.048611in}}{\pgfqpoint{-0.025257in}{0.043489in}}{\pgfqpoint{-0.034373in}{0.034373in}}%
\pgfpathcurveto{\pgfqpoint{-0.043489in}{0.025257in}}{\pgfqpoint{-0.048611in}{0.012892in}}{\pgfqpoint{-0.048611in}{0.000000in}}%
\pgfpathcurveto{\pgfqpoint{-0.048611in}{-0.012892in}}{\pgfqpoint{-0.043489in}{-0.025257in}}{\pgfqpoint{-0.034373in}{-0.034373in}}%
\pgfpathcurveto{\pgfqpoint{-0.025257in}{-0.043489in}}{\pgfqpoint{-0.012892in}{-0.048611in}}{\pgfqpoint{0.000000in}{-0.048611in}}%
\pgfpathclose%
\pgfusepath{fill}%
}%
\begin{pgfscope}%
\pgfsys@transformshift{6.005617in}{1.530346in}%
\pgfsys@useobject{currentmarker}{}%
\end{pgfscope}%
\end{pgfscope}%
\begin{pgfscope}%
\definecolor{textcolor}{rgb}{0.150000,0.150000,0.150000}%
\pgfsetstrokecolor{textcolor}%
\pgfsetfillcolor{textcolor}%
\pgftext[x=6.380617in,y=1.457430in,left,base]{\color{textcolor}\sffamily\fontsize{15.000000}{18.000000}\selectfont DMG}%
\end{pgfscope}%
\begin{pgfscope}%
\pgfsetroundcap%
\pgfsetroundjoin%
\pgfsetlinewidth{1.405250pt}%
\definecolor{currentstroke}{rgb}{0.501961,0.501961,0.501961}%
\pgfsetstrokecolor{currentstroke}%
\pgfsetdash{}{0pt}%
\pgfpathmoveto{\pgfqpoint{5.797283in}{1.224561in}}%
\pgfpathlineto{\pgfqpoint{6.213950in}{1.224561in}}%
\pgfusepath{stroke}%
\end{pgfscope}%
\begin{pgfscope}%
\pgfsetbuttcap%
\pgfsetroundjoin%
\definecolor{currentfill}{rgb}{0.501961,0.501961,0.501961}%
\pgfsetfillcolor{currentfill}%
\pgfsetlinewidth{0.000000pt}%
\definecolor{currentstroke}{rgb}{0.501961,0.501961,0.501961}%
\pgfsetstrokecolor{currentstroke}%
\pgfsetdash{}{0pt}%
\pgfsys@defobject{currentmarker}{\pgfqpoint{-0.034722in}{-0.034722in}}{\pgfqpoint{0.034722in}{0.034722in}}{%
\pgfpathmoveto{\pgfqpoint{0.000000in}{-0.034722in}}%
\pgfpathcurveto{\pgfqpoint{0.009208in}{-0.034722in}}{\pgfqpoint{0.018041in}{-0.031064in}}{\pgfqpoint{0.024552in}{-0.024552in}}%
\pgfpathcurveto{\pgfqpoint{0.031064in}{-0.018041in}}{\pgfqpoint{0.034722in}{-0.009208in}}{\pgfqpoint{0.034722in}{0.000000in}}%
\pgfpathcurveto{\pgfqpoint{0.034722in}{0.009208in}}{\pgfqpoint{0.031064in}{0.018041in}}{\pgfqpoint{0.024552in}{0.024552in}}%
\pgfpathcurveto{\pgfqpoint{0.018041in}{0.031064in}}{\pgfqpoint{0.009208in}{0.034722in}}{\pgfqpoint{0.000000in}{0.034722in}}%
\pgfpathcurveto{\pgfqpoint{-0.009208in}{0.034722in}}{\pgfqpoint{-0.018041in}{0.031064in}}{\pgfqpoint{-0.024552in}{0.024552in}}%
\pgfpathcurveto{\pgfqpoint{-0.031064in}{0.018041in}}{\pgfqpoint{-0.034722in}{0.009208in}}{\pgfqpoint{-0.034722in}{0.000000in}}%
\pgfpathcurveto{\pgfqpoint{-0.034722in}{-0.009208in}}{\pgfqpoint{-0.031064in}{-0.018041in}}{\pgfqpoint{-0.024552in}{-0.024552in}}%
\pgfpathcurveto{\pgfqpoint{-0.018041in}{-0.031064in}}{\pgfqpoint{-0.009208in}{-0.034722in}}{\pgfqpoint{0.000000in}{-0.034722in}}%
\pgfpathclose%
\pgfusepath{fill}%
}%
\begin{pgfscope}%
\pgfsys@transformshift{6.005617in}{1.224561in}%
\pgfsys@useobject{currentmarker}{}%
\end{pgfscope}%
\end{pgfscope}%
\begin{pgfscope}%
\definecolor{textcolor}{rgb}{0.150000,0.150000,0.150000}%
\pgfsetstrokecolor{textcolor}%
\pgfsetfillcolor{textcolor}%
\pgftext[x=6.380617in,y=1.151644in,left,base]{\color{textcolor}\sffamily\fontsize{15.000000}{18.000000}\selectfont AMG}%
\end{pgfscope}%
\begin{pgfscope}%
\pgfsetbuttcap%
\pgfsetroundjoin%
\pgfsetlinewidth{1.405250pt}%
\definecolor{currentstroke}{rgb}{0.000000,0.000000,0.000000}%
\pgfsetstrokecolor{currentstroke}%
\pgfsetdash{{5.180000pt}{2.240000pt}}{0.000000pt}%
\pgfpathmoveto{\pgfqpoint{5.797283in}{0.918775in}}%
\pgfpathlineto{\pgfqpoint{6.213950in}{0.918775in}}%
\pgfusepath{stroke}%
\end{pgfscope}%
\begin{pgfscope}%
\definecolor{textcolor}{rgb}{0.150000,0.150000,0.150000}%
\pgfsetstrokecolor{textcolor}%
\pgfsetfillcolor{textcolor}%
\pgftext[x=6.380617in,y=0.845858in,left,base]{\color{textcolor}\sffamily\fontsize{15.000000}{18.000000}\selectfont Boundary \(\displaystyle \rho = 1\)}%
\end{pgfscope}%
\end{pgfpicture}%
\makeatother%
\endgroup%

%% file: draft.bbl
\begin{thebibliography}{10}

\bibitem{abadi2016tensorflow}
{\sc M.~Abadi, P.~Barham, J.~Chen, Z.~Chen, A.~Davis, J.~Dean, M.~Devin,
  S.~Ghemawat, G.~Irving, M.~Isard, et~al.}, {\em Tensorflow: A system for
  large-scale machine learning.}, in OSDI, vol.~16, 2016, pp.~265--283.

\bibitem{avron2011randomized}
{\sc H.~Avron and S.~Toledo}, {\em Randomized algorithms for estimating the
  trace of an implicit symmetric positive semi-definite matrix}, Journal of the
  ACM (JACM), 58 (2011), p.~8.

\bibitem{baur1983complexity}
{\sc W.~Baur and V.~Strassen}, {\em The complexity of partial derivatives},
  Theoretical computer science, 22 (1983), pp.~317--330.

\bibitem{pyamg}
{\sc W.~N. Bell, L.~N. Olson, and J.~B. Schroder}, {\em {PyAMG}: Algebraic
  multigrid solvers in {Python} v3.0}, 2015,
  \url{https://github.com/pyamg/pyamg}.
\newblock Release 3.2.

\bibitem{briggs2000multigrid}
{\sc W.~L. Briggs, V.~E. Henson, and S.~F. McCormick}, {\em A multigrid
  tutorial}, SIAM, 2000.

\bibitem{dai2013fast}
{\sc R.~Dai, Y.~Wang, and J.~Zhang}, {\em Fast and high accuracy multiscale
  multigrid method with multiple coarse grid updating strategy for the 3d
  convection--diffusion equation}, Computers \& Mathematics with Applications,
  66 (2013), pp.~542--559.

\bibitem{de1990matrix}
{\sc P.~M. De~Zeeuw}, {\em Matrix-dependent prolongations and restrictions in a
  blackbox multigrid solver}, Journal of computational and applied mathematics,
  33 (1990), pp.~1--27.

\bibitem{dendy1982black}
{\sc J.~Dendy}, {\em Black box multigrid}, Journal of Computational Physics, 48
  (1982), pp.~366--386.

\bibitem{goodfellow2016deep}
{\sc I.~Goodfellow, Y.~Bengio, and A.~Courville}, {\em Deep Learning}, MIT
  Press, 2016.
\newblock \url{http://www.deeplearningbook.org}.

\bibitem{grasedyck2016nearly}
{\sc L.~Grasedyck, L.~Wang, and J.~Xu}, {\em A nearly optimal multigrid method
  for general unstructured grids}, Numerische Mathematik, 134 (2016),
  pp.~637--666.

\bibitem{hackbusch2013multi}
{\sc W.~Hackbusch}, {\em Multi-grid methods and applications}, vol.~4, Springer
  Science \& Business Media, 2013.

\bibitem{he2016deep}
{\sc K.~He, X.~Zhang, S.~Ren, and J.~Sun}, {\em Deep residual learning for
  image recognition}, in Proceedings of the IEEE conference on computer vision
  and pattern recognition, 2016, pp.~770--778.

\bibitem{hecht1988theory}
{\sc R.~Hecht-Nielsen et~al.}, {\em Theory of the backpropagation neural
  network.}, Neural Networks, 1 (1988), pp.~445--448.

\bibitem{hemker1990order}
{\sc P.~Hemker}, {\em On the order of prolongations and restrictions in
  multigrid procedures}, Journal of Computational and Applied Mathematics, 32
  (1990), pp.~423--429.

\bibitem{hutchinson1990stochastic}
{\sc M.~F. Hutchinson}, {\em A stochastic estimator of the trace of the
  influence matrix for laplacian smoothing splines}, Communications in
  Statistics-Simulation and Computation, 19 (1990), pp.~433--450.

\bibitem{kawulok2016advances}
{\sc M.~Kawulok, E.~Celebi, and B.~Smolka}, {\em Advances in Face Detection and
  Facial Image Analysis}, Springer, 2016.

\bibitem{khelifi2014hybrid}
{\sc S.~C. Khelifi, N.~M{\'e}chitoua, F.~H{\"u}lsemann, and F.~Magoul{\`e}s},
  {\em A hybrid multigrid method for convection--diffusion problems}, Journal
  of Computational and Applied Mathematics, 259 (2014), pp.~711--719.

\bibitem{kingma2014adam}
{\sc D.~Kingma and J.~Ba}, {\em Adam: A method for stochastic optimization},
  arXiv preprint arXiv:1412.6980,  (2014).

\bibitem{kozyakin2009accuracy}
{\sc V.~Kozyakin}, {\em On accuracy of approximation of the spectral radius by
  the gelfand formula}, Linear Algebra and its Applications, 431 (2009),
  pp.~2134--2141.

\bibitem{kressner2014subspace}
{\sc D.~Kressner and B.~Vandereycken}, {\em Subspace methods for computing the
  pseudospectral abscissa and the stability radius}, SIAM Journal on Matrix
  Analysis and Applications, 35 (2014), pp.~292--313.

\bibitem{krizhevsky2012imagenet}
{\sc A.~Krizhevsky, I.~Sutskever, and G.~E. Hinton}, {\em Imagenet
  classification with deep convolutional neural networks}, in Advances in
  neural information processing systems, 2012, pp.~1097--1105.

\bibitem{lecun2015deep}
{\sc Y.~LeCun, Y.~Bengio, and G.~Hinton}, {\em Deep learning}, Nature, 521
  (2015), pp.~436--444.

\bibitem{livshits2014scalable}
{\sc I.~Livshits}, {\em A scalable multigrid method for solving indefinite
  helmholtz equations with constant wave numbers}, Numerical Linear Algebra
  with Applications, 21 (2014), pp.~177--193.

\bibitem{maclaurin2015autograd}
{\sc D.~Maclaurin, D.~Duvenaud, and R.~P. Adams}, {\em Autograd: Effortless
  gradients in numpy}, in ICML 2015 AutoML Workshop, 2015.

\bibitem{mengi2014numerical}
{\sc E.~Mengi, E.~A. Yildirim, and M.~Kili{\c{c}}}, {\em Numerical optimization
  of eigenvalues of hermitian matrix functions}, SIAM Journal on Matrix
  Analysis and Applications, 35 (2014), pp.~699--724.

\bibitem{nesterov2007smoothing}
{\sc Y.~Nesterov}, {\em Smoothing technique and its applications in
  semidefinite optimization}, Mathematical Programming, 110 (2007),
  pp.~245--259.

\bibitem{nesterov2013optimizing}
{\sc Y.~Nesterov and V.~Y. Protasov}, {\em Optimizing the spectral radius},
  SIAM Journal on Matrix Analysis and Applications, 34 (2013), pp.~999--1013.

\bibitem{olshanskii2014iterative}
{\sc M.~A. Olshanskii and E.~E. Tyrtyshnikov}, {\em Iterative methods for
  linear systems: theory and applications}, SIAM, 2014.

\bibitem{park2015deep}
{\sc Y.~Park and M.~Kellis}, {\em Deep learning for regulatory genomics},
  Nature biotechnology, 33 (2015), pp.~825--826.

\bibitem{ruge1987algebraic}
{\sc J.~W. Ruge and K.~St{\"u}ben}, {\em Algebraic multigrid}, Multigrid
  methods, 3 (1987), pp.~73--130.

\bibitem{socher2012convolutional}
{\sc R.~Socher, B.~Huval, B.~Bath, C.~D. Manning, and A.~Y. Ng}, {\em
  Convolutional-recursive deep learning for 3d object classification}, in
  Advances in Neural Information Processing Systems, 2012, pp.~656--664.

\bibitem{socher2013recursive}
{\sc R.~Socher, A.~Perelygin, J.~Wu, J.~Chuang, C.~D. Manning, A.~Ng, and
  C.~Potts}, {\em Recursive deep models for semantic compositionality over a
  sentiment treebank}, in Proceedings of the 2013 conference on empirical
  methods in natural language processing, 2013, pp.~1631--1642.

\bibitem{sosnovik2017neural}
{\sc I.~Sosnovik and I.~Oseledets}, {\em Neural networks for topology
  optimization}, arXiv preprint arXiv:1709.09578,  (2017).

\bibitem{stolk2014multigrid}
{\sc C.~C. Stolk, M.~Ahmed, and S.~K. Bhowmik}, {\em A multigrid method for the
  helmholtz equation with optimized coarse grid corrections}, SIAM Journal on
  Scientific Computing, 36 (2014), pp.~A2819--A2841.

\bibitem{stuben2001review}
{\sc K.~St{\"u}ben}, {\em A review of algebraic multigrid}, Journal of
  Computational and Applied Mathematics, 128 (2001), pp.~281--309.

\bibitem{stynes2007convection}
{\sc M.~Stynes}, {\em Convection-diffusion-reaction problems, sdfem/supg and a
  priori meshes}, International Journal of Computing Science and Mathematics, 1
  (2007), pp.~412--431.

\bibitem{sun2014deep}
{\sc Y.~Sun, Y.~Chen, X.~Wang, and X.~Tang}, {\em Deep learning face
  representation by joint identification-verification}, in Advances in neural
  information processing systems, 2014, pp.~1988--1996.

\bibitem{theano}
{\sc {Theano Development Team}}, {\em {Theano: A {Python} framework for fast
  computation of mathematical expressions}}, arXiv e-prints, abs/1605.02688
  (2016), \url{http://arxiv.org/abs/1605.02688}.

\bibitem{trottenberg2000multigrid}
{\sc U.~Trottenberg, C.~W. Oosterlee, and A.~Schuller}, {\em Multigrid},
  Academic press, 2000.

\bibitem{watson1989modern}
{\sc L.~T. Watson and R.~T. Haftka}, {\em Modern homotopy methods in
  optimization}, Computer Methods in Applied Mechanics and Engineering, 74
  (1989), pp.~289--305.

\bibitem{xu1996auxiliary}
{\sc J.~Xu}, {\em The auxiliary space method and optimal multigrid
  preconditioning techniques for unstructured grids}, Computing, 56 (1996),
  pp.~215--235.

\bibitem{yi2014deep}
{\sc D.~Yi, Z.~Lei, S.~Liao, and S.~Z. Li}, {\em Deep metric learning for
  person re-identification}, in Pattern Recognition (ICPR), 2014 22nd
  International Conference on, IEEE, 2014, pp.~34--39.

\bibitem{yserentant1993old}
{\sc H.~Yserentant}, {\em Old and new convergence proofs for multigrid
  methods}, Acta numerica, 2 (1993), pp.~285--326.

\end{thebibliography}
